\documentclass[twocolumn]{svjour3}          

\usepackage[colorlinks,citecolor=black,linkcolor=black,urlcolor=black]{hyperref}
\usepackage[figuresright]{rotating}

\usepackage{amssymb}
\usepackage{amsmath,amssymb}
\usepackage{subfig}
\usepackage{stmaryrd}

\usepackage{booktabs}
\usepackage{algorithm}
\usepackage{algorithmicx}
\usepackage{algpseudocode}
\usepackage{tabularx}
\usepackage{xcolor}
\usepackage{caption}
\usepackage[authoryear,round]{natbib}
\DeclareMathOperator{\buffer}{buffer}

\DeclareMathOperator{\supp}{supp}
\DeclareMathOperator{\trunc}{trunc}

\DeclareMathOperator{\N}{N}
\DeclareMathOperator{\G}{G}
\DeclareMathOperator{\U}{U}
\DeclareMathOperator{\T}{T}

\newcolumntype{P}[1]{>{\centering\arraybackslash}p{#1}}
\newcolumntype{M}[1]{>{\centering\arraybackslash}m{#1}}

\journalname{}

\begin{document}
	
	\title{Extended isogeometric analysis of multi-material and multi-physics problems using hierarchical B-splines}
	
	\author{Mathias Schmidt$^1$, Lise No\"el$^2$, Keenan Doble$^3$, John A. Evans$^3$, Kurt Maute$^3$}

		\institute{{$^1$Computational Engineering Division, Lawrence Livermore National Laboratory}, 
        {7000 East Ave}, {Livermore}, {94550}, {California}, {USA} \\ 
        {$^2$Department of Precision and Microsystems Engineering, Faculty of Mechanical, Maritime and Materials Engineering},
        {Delft University of Technology}, {Mekelweg 2}, {2628 CD Delft} , {The Netherlands}\\
        {$^3$Aerospace Mechanics Research Center, Department of Aerospace Engineering Sciences}, 
		{University of Colorado Boulder}, {3775 Discovery Dr}, {Boulder}, {80309-0429}, {Colorado}, {USA} \\ 
  		\email{schmidt43@llnl.gov} }
	
	\date{Received: date / Accepted: date}
	
	\maketitle
	
\begin{abstract}
This paper presents an immersed, isogeometric finite element framework to predict the response of multi-material, multi-physics problems with complex geometries using locally refined discretizations. To circumvent the need to generate conformal meshes, this work uses an eXtended Finite Element Method (XFEM) to discretize the governing equations on non-conforming, embedding meshes. A flexible approach to create truncated hierarchical B-splines discretizations is presented. This approach enables the refinement of each state variable field individually to meet field-specific accuracy requirements. To obtain an immersed geometry representation that is consistent across all hierarchically refined B-spline discretizations, the geometry is immersed into a single mesh, the XFEM background mesh, which is constructed from the union of all hierarchical B-spline meshes. An extraction operator is introduced to represent the truncated hierarchical B-spline bases in terms of Lagrange shape functions on the XFEM background mesh without loss of accuracy. The truncated hierarchical B-spline bases are enriched using a generalized Heaviside enrichment strategy to accommodate small geometric features and multi-material problems. The governing equations are augmented by a formulation of the face-oriented ghost stabilization enhanced for locally refined B-spline bases. We present examples for two- and three-dimensional linear elastic and thermo-elastic problems. The numerical results validate the accuracy of our framework.  The results also demonstrate the applicability of the proposed framework to large, geometrically complex problems.
\end{abstract}

\keywords{Immersed Finite Element Method, eXtended Isogeometric Analysis, Multi-material Problems, Multi-physics Problems, Truncated Hierarchical B-splines, Lagrange Extraction}
\section{Introduction}\label{sec:intro}

Finite element analysis is frequently used to predict the response of systems described by partial differential equations defined over a spatial domain. In classical finite element methods, the domain is discretized using a single mesh that conforms to the external boundaries and internal material interfaces. For problems with complex shapes and multiple material phases, the construction of this conformal mesh is often a major bottleneck in the analysis process, see \cite{Bazilevs2010}. Conformal mesh generation may also hamper the automation of the finite element analysis for problems with changing geometry as encountered in, for example, phase-change problems with dynamically evolving interfaces or shape and topology optimization.

The state variable fields may exhibit large spatial gradients at boundaries and material interfaces, as well as in the vicinity of small geometric features. To resolve these spatial gradients, a sufficiently fine discretization is needed. While uniformly fine discretization may yield an accurate approximation, the associated computational cost may exceed practical limits. A locally refined mesh balances discretization needs and computational cost. However, generating locally refined meshes for standard finite element methods further increases the complexity of mesh generation. 
 
\begin{figure*}
	\centering
	\subfloat[Thermo-elastic multi-material problem subject to a spatially varying heat load\label{subfig-2:1}]
	{\includegraphics[width=0.7\columnwidth]{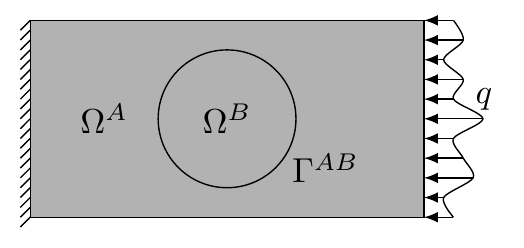}}
	\hfill
	\subfloat[Norm of the diffusive flux vector \label{subfig-3:1}]
	{\includegraphics[width=0.6\columnwidth]{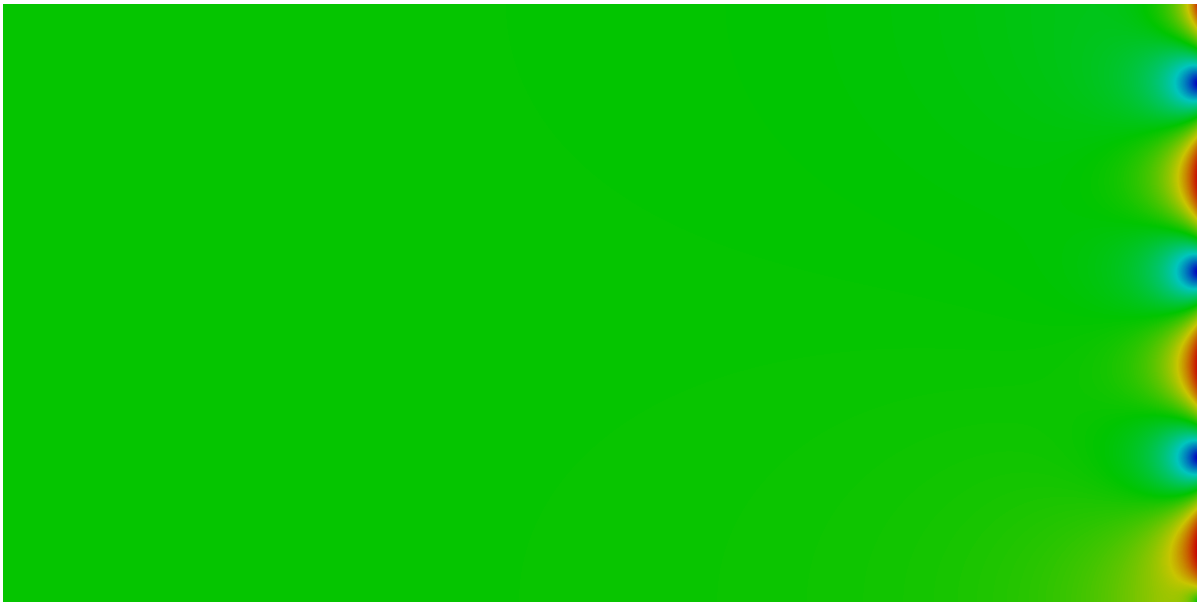}}
	\hfill
	\subfloat[Von Mises stress field  \label{subfig-4:1}]
	{\includegraphics[width=0.6\columnwidth]{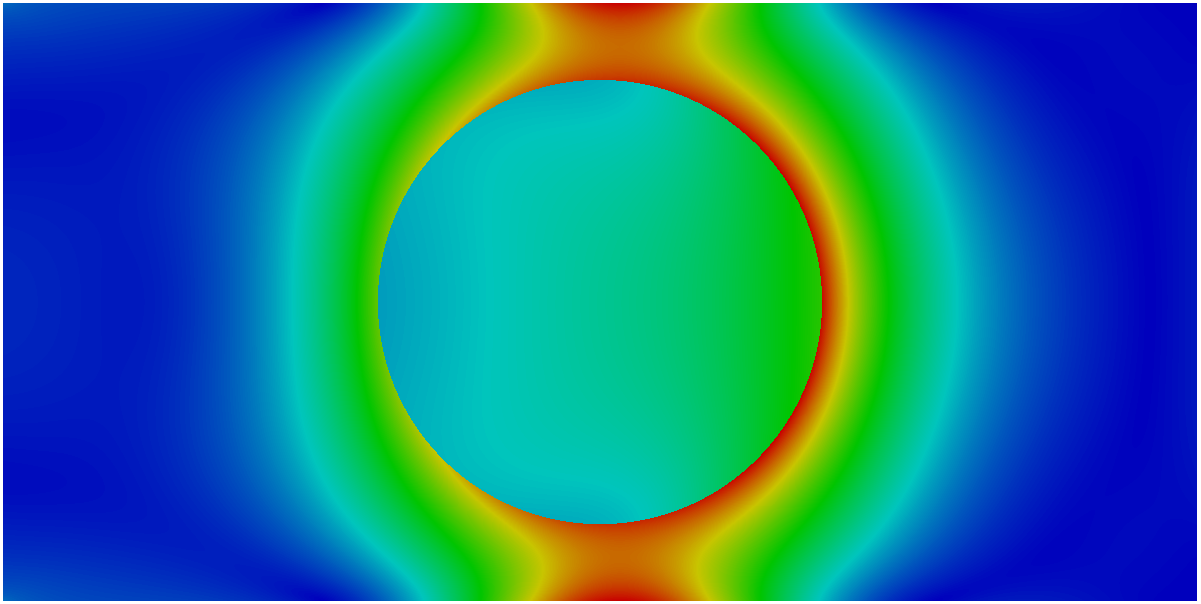}}
	\caption{Illustration of an example thermo-elastic problem}
	\label{Fig_introduction}
\end{figure*}

\begin{figure*}
	\centering
	\subfloat[Locally refined mesh for the temperature field\label{subfig-2:dummy}]
	{\includegraphics[width=0.3\textwidth]{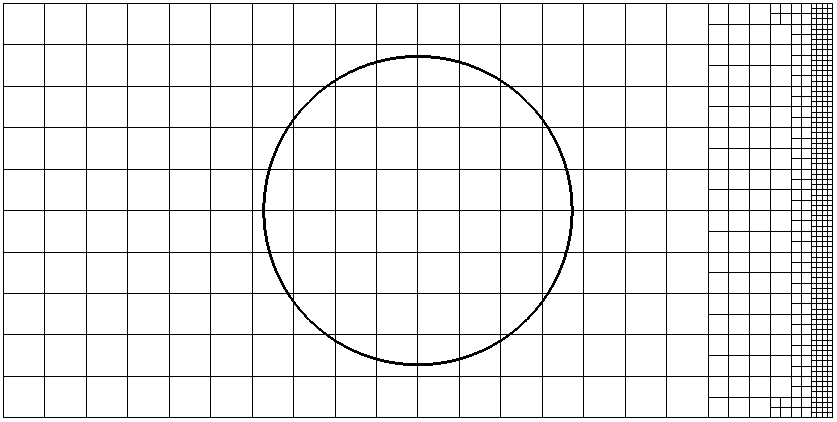}}
	\hspace{1.0cm}
	\subfloat[Locally refined mesh for the displacement field\label{subfig-3:dummy}]
	{\includegraphics[width=0.3\textwidth]{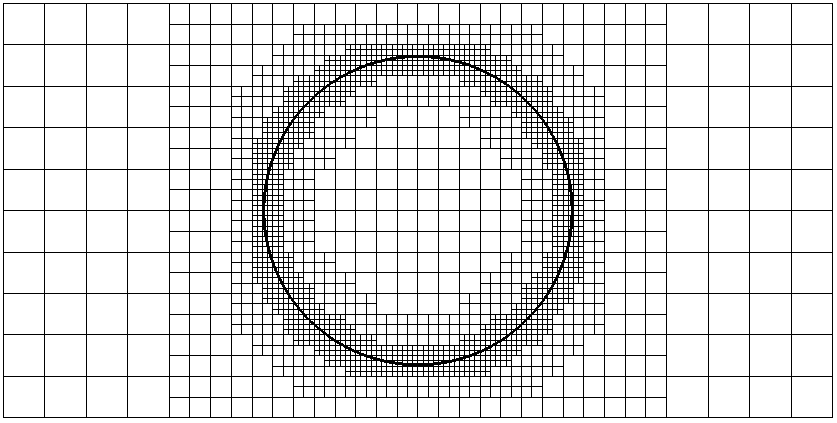}}
	\caption{Example of locally refined discretizations for the thermo-elastic problem depicted in Fig.~\ref{subfig-2:1}}
	\label{Fig_introduction_discretization}
\end{figure*}

Many engineering and science applications involve multiple, often coupled state variable fields. For such problems, each field may require differently refined discretizations. For example, consider the two-material, thermo-elastic problem of a circle embedded in a rectangular plate shown in Fig.~\ref{subfig-2:1}. We assume that the circle is made of a material with a finite coefficient of thermal expansion (CTE), while the plate's CTE is zero. The plate is subject to a spatially varying heat flux along its right edge and clamped at its left edge. The heat flux increases the temperature in the system, causing the circle to expand which generates stresses in both the plate and circle due to the CTE mismatch. To predict the stress field at the material interface with high accuracy, accurate approximations of both the temperature and displacement fields are required. The contour of the norm of the diffusive flux is shown in Fig.~\ref{subfig-3:1} and the contour of the Von Mises stress field in Fig.~\ref{subfig-4:1}. The regions with large spatial gradients differ between the displacement and the temperature field. Separately adapting the discretizations associated with each field enables a sufficient resulution of each field while minimizing the overall computational cost. An example of such discretizations for the temperature and displacement fields is depicted in Fig.~\ref{subfig-2:dummy} and \ref{subfig-3:dummy}, respectively.


Lagrange polynomials are the standard choice for basis functions in finite elements. The element-local nature of this class of basis functions and their interpolatory property facilitate a convenient and efficient implementation and enforcement of Dirichlet boundary and interface conditions. However, Lagrange basis functions limit the inter-element continuity to $C^{0}$, irrespective of the polynomial order of the approximation within the element. The lack of higher-order inter-element continuity affects the efficiency of Lagrange bases, measured by the number of degrees of freedom needed to achieve a desired accuracy; see for example \cite{Evans2009}.

In this work, we propose a computational framework that addresses the issues encountered with standard finite element approaches as outlined above. This framework synthesizes the following concepts which will be further discussed in detail below. An immersed finite element approach is adopted to eliminate the need for conformal mesh generation. Higher order B-spline basis functions are used to increase the discretization efficiency over Lagrange bases. For each state variable, a truncated hierarchically refined B-spline basis is generated to meet field-specific accuracy requirements. 

Immersed Boundary Methods (IBMs) have gained in popularity in recent years, see \cite{Babuska1997,Peskin2002,Mittal2005}. The general idea of these methods is to immerse the geometry of the physical domain into a computational domain with a much simpler geometry. Thus, the geometries of the physical and computational domains are decoupled. The geometric simplicity of the computational domain allows for the convenient generation of IBM background meshes, such as tensor grids. The geometry is immersed into the IBM background mesh, which simplifies the discretization of the state variable fields.

In this work, we focus on finite element formulations of IBMs. Immersed Finite Element Methods (IFEMs), also referred to as geometrically unfitted or embedded domain finite element methods, include the Finite Cell Method , see \cite{Parvizian2007a,Duster2008,Schillinger2015}, the CutFEM (\cite{Burman2015}), the Generalized Finite Element Method (GFEM) (\cite{Duarte2000,Strouboulis2000}), and the eXtended Finite Element Method (XFEM), see \cite{Belytschko1999,Belytschko2009}. In this paper, we consider specifically the latter approach. The XFEM augments the standard finite element basis with additional basis functions to represent discontinuities of the state variable field within an XFEM background element that is intersected by a boundary or an interface. The augmented finite element basis satisfies the partition of unity (PU) property. In this work, we adopt a Heaviside enrichment strategy for its flexibility in modeling interface and boundary conditions of multi-material problems with complex geometries, see \cite{Noel2022}. 

Traditionally, Heaviside enriched XFEM approaches discretize state variable fields by Lagrange basis functions which are defined on the XFEM background mesh. In this paper, we adopt higher order B-splines for discretizing state variable fields. The advantages of B-splines basis functions for finite element methods have been demonstrated in the context of Isogeometric Analysis (IGA), see \cite{Hughes2005}. While IGA was originally developed to eliminate the discrepancy between CAD geometry representation and finite element analysis, \cite{Evans2009} showed that B-spline basis functions in general improve accuracy, robustness, and computational efficiency. In this work, we further utilize the refineability property of B-splines, see \cite{Garau2016}. 

B-spline basis functions were studied with IFEMs by \cite{Methods2012,Schillinger2011,Schillinger2015,Verhoosel2015,Elfverson2018,Divi2020,Noel2020}. \cite{Nguyen2012} introduced B-spline approximations to the XFEM and coined the term \textit{X-IGA}. These works demonstrated that integrating higher order B-splines discretizations into IFEM approaches yields improved accuracy and computational efficiency compared to the use of Lagrange basis functions. 

The ability to locally refine B-splines enables the convenient construction of locally refined approximation spaces, see \cite{Giannelli2012a,Buffa2017,Bracco2019}. While standard hierarchical B-spline bases do not fulfill the PU property, truncating the bases restores the PU property \cite{Giannelli2012a}. Truncated Hierarchical B-spline (THB) basis functions form a sparse, strongly stable basis \cite{Giannelli2012a} and are employed in this work.

\begin{figure*}
  \centering
  \includegraphics[width=0.3\textwidth]{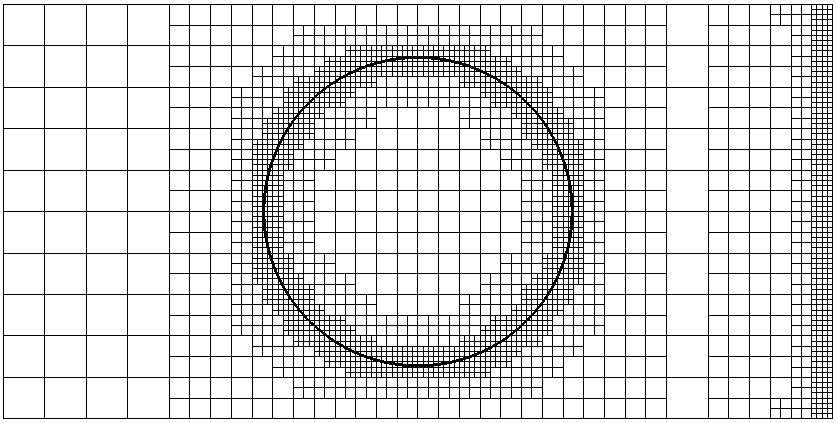}
  \caption{Union of the locally refined meshes depicted in Figures~\ref{subfig-2:dummy} and \ref{subfig-3:dummy} \label{subfig-4:dummy}}
  \label{Fig_introduction_discretization_union}
\end{figure*}

This paper contributes to the work on IFEMs as follows. We introduce a versatile discretization approach based on THBs defined on locally refined tensor meshes. Each state variable can be discretized individually by choosing the polynomial order and the local refinement independently from the discretization of other fields. This approach leads to an explicit control over the accuracy of the field approximations and the computational cost of the analysis. Computationally efficient algorithms and data structures are introduced to enable the application of this discretization approach to large problems in two and three dimensions using parallel computing. To integrate the locally refined B-spline discretization approach into the XFEM, a single union mesh is generated from the individual B-spline meshes. Fig.~\ref{subfig-4:dummy} shows the union mesh generated for the two locally refined meshes presented in Fig.~\ref{Fig_introduction_discretization}. The union mesh serves as the XFEM background mesh in which the geometry is immersed. To increase the geometry resolution, the XFEM background mesh can be additionally refined. We introduce an exact extraction operator that expresses THBs in terms of Lagrange basis functions over each element of the XFEM background mesh.

The enrichment strategy of \cite{Makhija2014,Villanueva2014,Noel2022} is generalized to enrich the locally refined THB basis functions considering their support for a given intersection geometry. The XFEM problem is augmented by an X-IGA formulation of the face-oriented ghost stabilization from \cite{Noel2022}. The governing equations are integrated by standard quadrature rules on a geometry-conforming integration mesh that is constructed by cutting XFEM background elements that are intersected by the external boundaries or internal material interfaces. Boundary and interface conditions are enforced weakly by Nitsche's method. We illustrate the main characteristics and features of the proposed XFEM framework by numerical examples considering linear elastic and thermo-elastic problems. We study the convergence of geometric and state variable discretization errors with uniform and local mesh refinement for different B-spline orders. To demonstrate the applicability of our XFEM framework to complex multi-material, multi-physics problems, we apply our XFEM framework to the thermo-elastic analysis of a 3D polycrystalline micro-structure.

This paper is organized as follows: Section \ref{sec:xfem} outlines the Heaviside-enriched XFEM framework. Section \ref{sec:bsplines} recalls the fundamentals of hierarchical B-splines as well as their truncation to restore the PU property. Section \ref{sec:details} details the meshing algorithms and data structures. Sections \ref{sec:model} summarizes the governing equations, Nitsche's formulation, and face-oriented ghost stabilization for linear thermo-elasticity. Numerical two and three dimensional examples are presented in Section \ref{numericalExamples}. The main findings are summarized in Section 7, together with recommendations for future work.

\section{The eXtended Finite Element Method}\label{sec:xfem}

In this work, we adopt the XFEM to perform analysis on non-conforming background meshes. We follow the basic concepts of a generalized Heaviside enrichment strategy, introduced by \cite{Terada2003}. In this section, the basic concepts of the XFEM as relevant for this work are briefly described. 

\begin{figure*}
	\centering
	\includegraphics[width=0.24\textwidth]{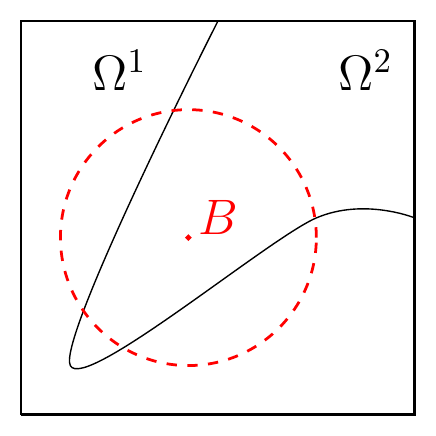}
	\includegraphics[width=0.24\textwidth]{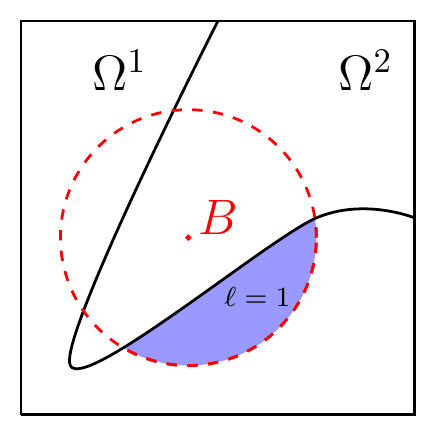}
	\includegraphics[width=0.24\textwidth]{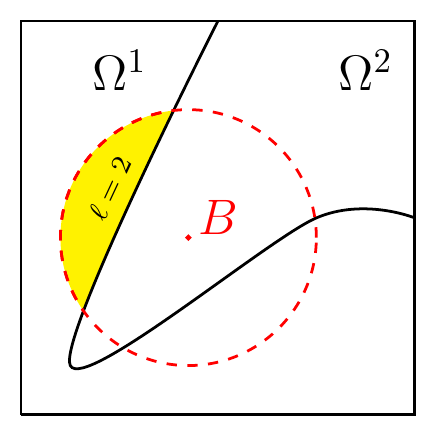}
	\includegraphics[width=0.24\textwidth]{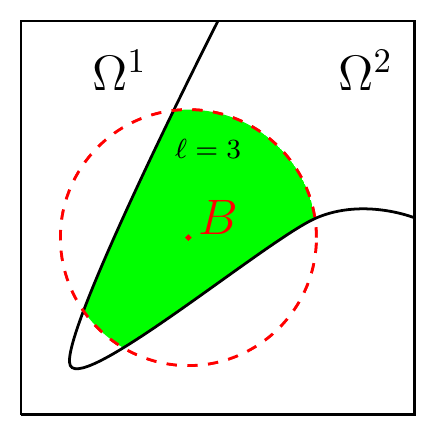}
	\caption{Basis function enrichment for a two-material problem}
	\label{Fig_basis_function_enrichment}
\end{figure*}

We start from a non-conforming finite element approximation space with THB basis functions $\left\{ B_j \right\}_{j=1}^{\tilde{n}}$ that are defined on a locally refined discretization. To approximate the state variable fields in different material phases, each basis function is enriched. The enrichment of a particular basis function depends on the number of topologically disconnected regions of all material phases in the support of this basis function. The approach is illustrated in Fig.~\ref{Fig_basis_function_enrichment} for a configuration with two material phases, $\Omega^1$ and $\Omega^2$. The support of the basis function $B$ is depicted by the dashed red line. It spans three topologically disconnected regions, each occupied by one  of the two material phases. Therefore, the basis function is enriched three times.

This approach is applicable to any number of material phases, intersection configurations, and basis function support. To define the enriched basis, we introduce the indicator function $\psi_{j}^{\ell}$. This function is equal to 1 for points located within the portion of the support of the $j^\textup{th}$ basis function corresponding to the enrichment level $\ell$ and equal to zero elsewhere. The enriched basis then is simply $\left\{B_j \psi_{j}^{\ell}: j \in \left\{1, \ldots, \tilde{n}\right\} \textup{ and } \ell \in \left\{1, \ldots, L_{j}\right\} \right\}$ where $L_{j}$ is the maximum number of enrichment levels for basis function $B_j$.  Since the original THB basis satisfies the PU principle, so does the enriched THB basis.  

The enriched finite element approximation of a vector-valued state variable $\mathbf{u}^h(\mathbf{x})$ can be written in terms of the enriched THB basis as follows: 
\begin{equation}
\mathbf{u}^h(\mathbf{x})
=\sum_{j = 1}^{\tilde{n}} \left( \sum_{l=1}^{L_j} \psi_{j}^{\ell}(\mathbf{x}) \, B_j(\mathbf{x}) \, \mathbf{c}_{j}^{\ell} \right),
\label{Eq_EnrichedApprox}
\end{equation}
where $\mathbf{c}_{j}^{\ell}$ is the coefficient associated to the $j^\text{th}$ original THB basis function and the $j^\text{th}$ enrichment level.  The indicator function $\psi_{j}^{\ell}$ enforces that only one set of enriched basis functions is used to describe the state variable at each point in the domain. A more detailed description of the enrichment strategy can be found in \cite{Noel2022}.
 
The Heaviside enriched XFEM formulation outlined above enables the modeling of $C^{-1}$ intra-element discontinuities of state variables within a non-conforming background element. Essential boundary conditions can be enforced weakly by, for example, Nitsche's method (\cite{Nitsche1971,Burman2012a}) or the stabilized Lagrange multiplier method (\cite{Gerstenberger2008}).

Immersing geometry into the XFEM background mesh can result in basis functions with small support within the geometric domain, leading to poorly conditioned systems of discretized governing equations. Various strategies to mitigate this issue have been studied in the literature, such as the face-oriented ghost stabilization (\cite{Burman2010,Burman2014a,Noel2022}), basis function removal (\cite{Embar2010,Elfverson2018}), and pre-conditioning (\cite{Lang2014,DePrenter2017}). In this work, we extend the face-oriented ghost stabilization to hierarchically refined B-spline discretizations; as further described in Subsection \ref{subsectionStabilization}.

\section{Hierarchical B-splines}\label{sec:bsplines}

This section focuses on hierarchical B-splines for a locally refined discretization of state variable fields. First, the basic concepts of B-splines in one and multiple dimensions are recalled. Then the B-spline refinement and the construction of non-truncated (HB) and truncated hierarchical B-spline (THB) bases are described.

\subsection{B-spline Basis Functions}\label{BSplineBasisFunctions}

In 1D, we define a knot vector $\Xi = \{ \xi_{1},\xi_{2}, \dots ,\xi_{n+p+1} \}$, for which $\xi \in \mathbb{R}$ and $\xi_{1} \leq \xi_{2} \leq \dots \leq \xi_{n+p+1}$. A univariate B-spline basis function $N_{i,p}(\xi)$ of degree $p$ is constructed recursively starting from the piecewise constant basis function:
\begin{equation}
N_{i,0}(\xi) =
\begin{cases}
1, & \text{if}\ \xi_{i} \leq \xi \leq\xi_{i+1},\\
0, & \text{otherwise},
\end{cases}
\end{equation}
and using the Cox de Boor recursion formula (\cite{Calculating1972}) for higher degrees, $p > 0$:
\begin{equation}
\begin{split}
N_{i,p}(\xi) = & \frac{\xi - \xi_{i}}{\xi_{i+p} - \xi_{i}}\ N_{i,p-1}(\xi)\\
& \hspace{1.25cm} + \frac{\xi_{i+p+1} - \xi}{\xi_{i+p+1} - \xi_{i+1}}\ N_{i+1,p-1}(\xi).
\end{split}
\end{equation}
To guarantee a $C^{p-1}$ continuity over the entire computational domain, none of the interior knots should be repeated.  The corresponding B-spline basis exhibits a $C^{p-1}$ continuity at every knot in the interior of the domain, while it is $C^{\infty}$ continuous in between the knots. A knot span is defined as the half open interval $[\xi_i, \xi_{i+1})$ and a B-spline element is defined as a non-empty knot span.  

In 2D and 3D, tensor-product B-spline basis functions $B_{i,p}(\xi)$ are constructed by applying the tensor-product operation to univariate B-spline basis functions in each parametric direction. Denoting the parametric space dimension as $d_{p}$, a tensor-product B-spline basis is constructed starting from $d_{p}$ knot vectors $\Xi^{m} = \{\xi_{1}^{m},\xi_{2}^{m},\-\dots,\-\xi_{n_{m}+p_{m}+1}^{m}\}$ with $p_{m}$ being the polynomial degree and $n_{m}$ the number of basis functions in the parametric direction $m = 1, \dots, d_{p}$. A tensor-product B-spline basis function is generated from $d_{p}$ univariate B-spline basis functions $N_{i_{m},p_{m}}^{m}(\xi^{m})$ in each parametric direction $m$ as follows:
\begin{equation}
B_{\mathbf{i},\mathbf{p}}(\boldsymbol{\xi}) = \prod_{m=1}^{d_{p}} N_{i_{m},p_{m}}^{m}(\xi^{m}),
\label{Eq_multiTensorProd}
\end{equation}
where the position in the tensor-product structure is given by the index vector $\mathbf{i} =\{ i_{1}, \dots, i_{d_{p}} \}$, and the vector $\mathbf{p} = \{ p_{1}, \dots, p_{d_{p}} \}$ defines the polynomial degree in each direction. Similarly to the univariate case, an element is defined as the tensor-product of $d_{p}$ non-empty knot spans. Additionally, a B-spline space $\mathcal{V}$ is defined as the span of B-spline basis functions.

\subsection{B-spline Refinement}\label{BSplinesRefinement}

\begin{figure}[h]\centering
	\includegraphics[width=0.90\columnwidth]{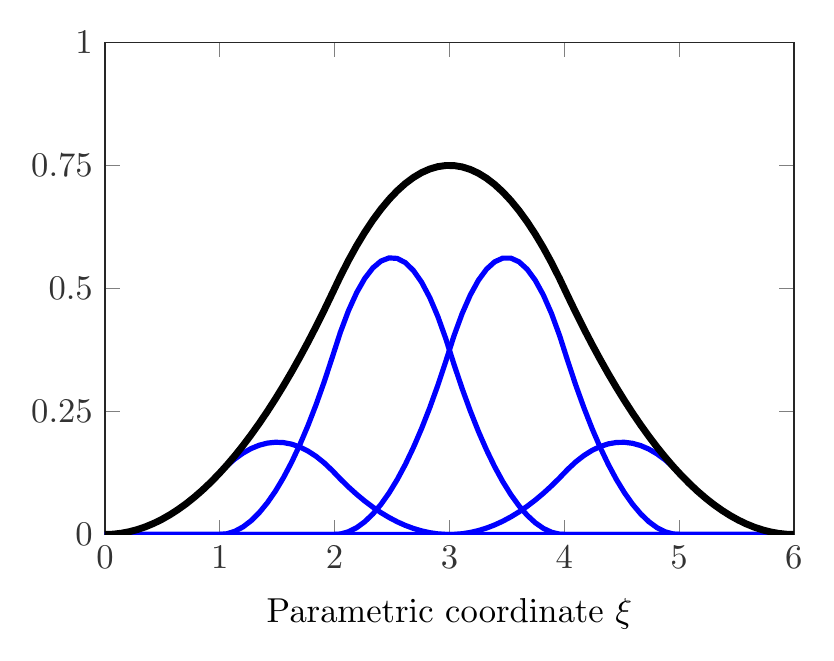}
	\caption{Subdivision of a quadratic B-spline basis function (black) into $p+2$ contracted B-spline basis functions (blue) of half the knot span width}
	\label{fig_BSplineSubdivision}
\end{figure}

Hierarchical refinement of uniform B-splines is achieved by subdivision. A univariate B-spline basis function is expressed as a linear combination of $p+2$ contracted, translated, and scaled copies of itself:
\begin{equation}
N_{p}(\xi) = 2^{-p}\sum_{j=0}^{p+1}\binom{p+1}{j}N_{p}(2\xi-j),
\end{equation}
where the binomial coefficient is defined as:
\begin{equation}
\binom{p+1}{j} = \frac{(p+1)!}{j!(p+1-j)!}.
\end{equation}
Fig.~\ref{fig_BSplineSubdivision} shows the refinement of a quadratic univariate B-spline basis function obtained by subdivision.

The subdivision in Eq.~\eqref{Eq_multiTensorProd} for a univariate B-spline basis can be extended to tensor-product B-spline basis functions $B_{p}$ as follows, see \cite{Methods2012}:
\begin{equation}
B_{\mathbf{p}}(\boldsymbol{\xi}) = \sum_{\mathbf{j}} \left( \prod_{m=1}^{d} 2^{-p_{m}}\ \binom{p_{m}+1}{j_{m}}\ N_{p_{m}}(2\xi^{m}-j_{m}) \right),
\end{equation}
where the index vector $\mathbf{j} = \{i_{1}, \dots, i_{d_{p}}\}$ collects the positions in the tensor-product structure.

\subsection{Hierarchical B-splines}\label{hierarchicalBSplinesSub}

To define a hierarchical mesh of depth $n$, a sequence of subdomains $\Omega^{l}$ is introduced:
\begin{equation}
\Omega^{n-1} \subseteq \Omega^{n-2} \subseteq \dots \subseteq \Omega^{0} = \Omega,
\end{equation}
where each subdomain $\Omega^l$ is a refined sub-region of $\Omega^{l-1}$. Consequently, $\Omega$ is equal to the union of all the subdomains $\Omega^l$.

\begin{figure}[b!]\center
	\includegraphics[width=0.90\columnwidth]{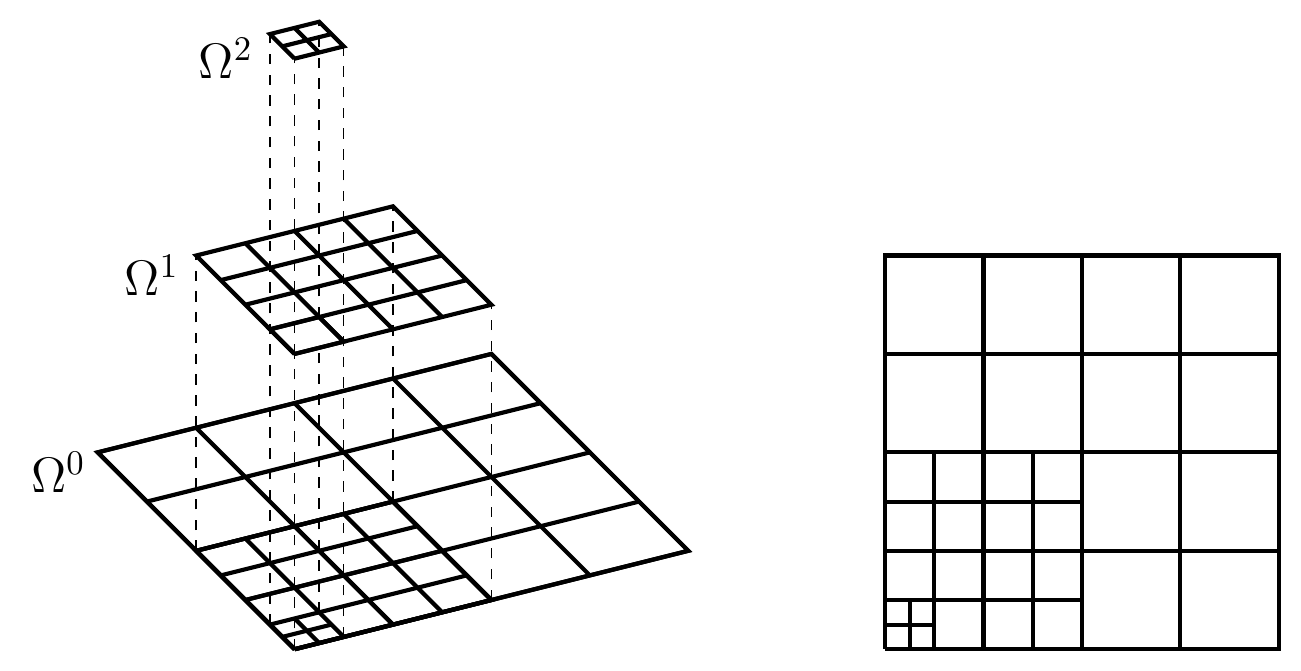}\hspace{0.25cm}
	\caption{Hierarchically refined mesh}
	\label{fig_hierarchical_Mesh}
\end{figure}

To build a hierarchical B-spline basis, a sequence of tensor-product B-spline spaces is introduced:
\begin{equation}
\mathcal{V}^{0} \subset \mathcal{V}^{1} \subset \mathcal{V}^{2} \subset \mathcal{V}^{3} \subset \dots .
\end{equation}
Each B-spline space $\mathcal{V}^{l}$ has a corresponding basis $\mathcal{B}^{l}$ and a corresponding mesh $\mathcal{K}^l$ of tensor-product elements.

A hierarchical B-spline basis $\mathcal{H}$ is constructed recursively based on the sequence of B-spline bases $\mathcal{B}^{l}$ that span the domains $\Omega^{l}$. In an initial step, the basis functions defined on the coarsest level, $l=0$, are collected and assigned to $\mathcal{H}^{0}$. The hierarchical B-spline basis $\mathcal{H}^{l+1}$ is constructed by taking the union of all basis functions $B$ in $\mathcal{H}^{l}$ whose support is not fully enclosed in $\Omega^{l+1}$ and all basis functions $B$ in $\mathcal{B}^{l+1}$ whose support lies in $\Omega^{l+1}$. The recursive algorithm reads (\cite{Garau2016}):
\begin{equation}
\left\{
\begin{array}{lll}
\mathcal{H}^0      & =: & \mathcal{B}^0\\[5pt]
\mathcal{H}^{l+1} & =: & \{ B \in \mathcal{H}^{l}\ |\ \supp(B) \not\subseteq \Omega^{l+1}\}\ \cup \\[5pt]
&&\{ B \in \mathcal{B}^{l+1}\ |\ \supp(B) \subseteq \Omega^{l+1}\},\\[5pt]
&&\mbox{for}\ l = 0, \dots, n-2,\\
\end{array}
\right.
\label{Eq_recursiveAlgorithm}
\end{equation}
where the index $l$ denotes the level of refinement. Basis functions collected in $\mathcal{H}$, where $\mathcal{H} := \mathcal{H}^{n-1}$, are called active, while basis functions in $\mathcal{B}^{l}$ not present in $\mathcal{H}$ are said to be inactive.

Associated with a hierarchical B-spline basis is a hierarchically refined mesh
\begin{equation}
\mathcal{K} := \cup_{l = 0}^{n-1} \left\{ K \in \mathcal{K}^l: K \in \Omega^{l} \text{ and } K \notin \Omega^{l+1} \right\}
\end{equation}
wherein $\Omega^n$ is taken to be the empty set.  An example of a hierarchically refined mesh associated with a two-dimensional hierarchical B-spline basis is displayed in Figure \ref{fig_hierarchical_Mesh}.  A hierarchical B-spline basis is smooth over each element of its associated hierarchically refined mesh.

\begin{figure}
\center
	\includegraphics[width=1.0\columnwidth]{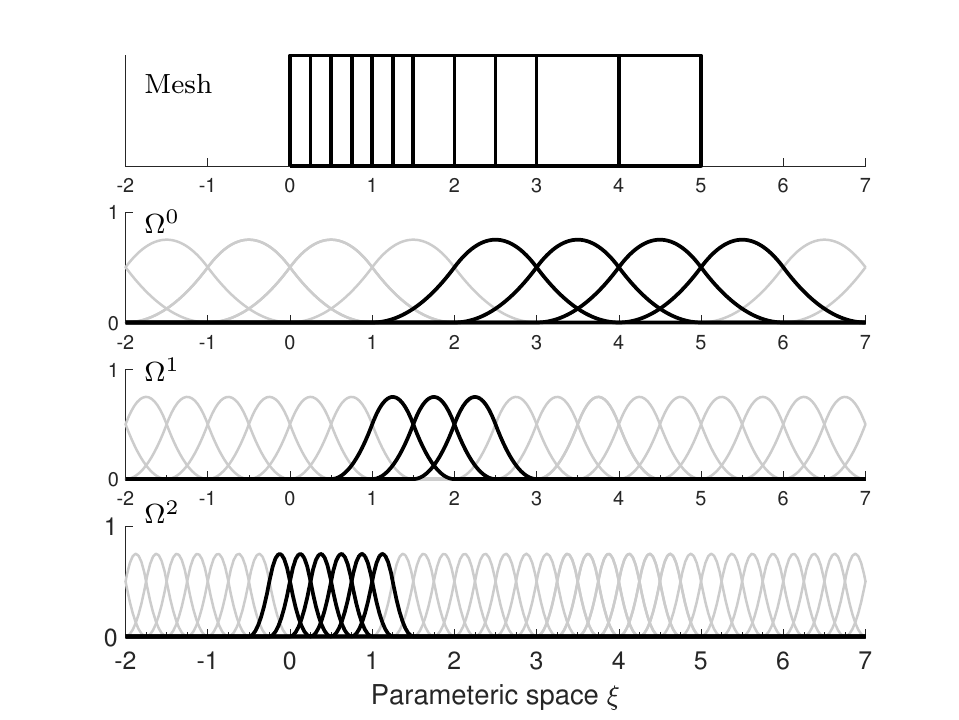}
	\caption{One-dimensional hierarchical B-spline basis attained using three levels of refinement }
	\label{fig_MeshRefinement1D}
\end{figure}

A hierarchical B-spline basis $\mathcal{H}$ is illustrated for a one-dimensional example in Fig.~\ref{fig_MeshRefinement1D}. The top row shows a one-dimensional hierarchically refined mesh. Below the mesh, the basis functions for three refinement levels are shown. Following the recursion rule of Eq.~\eqref{Eq_recursiveAlgorithm}, a B-spline basis $\mathcal{H}$ is created through an initialization step with all bases in the subdomain $\Omega^{0}$ refined to a level $l=0$. All bases in the subdomain $\Omega^{l+1}$ with higher refinement level $l+1$ are added recursively, while existing basis functions of level $l$ fully enclosed in $\Omega^{l+1}$ are discarded.
The active B-spline basis functions $\mathcal{H}$ are shown in black, while the inactive B-spline basis functions are shown in gray.

Hanging nodes are a byproduct of \textit{h}-refinement in hierarchical refined meshes and are naturally handled by the B-spline bases. In contrast, hanging nodes in classic finite elements require extra treatments, such as the introduction of multi-point constraints.

\subsection{Truncated B-splines}\label{truncatedBSplines}

By construction, the hierarchical B-spline basis presented above does not fulfill the PU property. The truncated hierarchical B-spline basis constitutes an alternative to the hierarchical B-spline basis that does satisfy the PU property. Truncation also reduces the number of overlapping functions on adjacent hierarchical levels, see \cite{Giannelli2012a}. Considering a basis function $B^{l}$, part of $\mathcal{B}^{l}$ and defined on the domain $\Omega^{l}$, its representation in terms of the finer basis of level $l+1$ is given as:
\begin{equation}
B^{l} = \sum_{B^{l+1} \in \, \mathcal{B}^{l+1}} c_{B^{l+1}}^{l+1} \left( B^{l} \right)\ B^{l+1},
\end{equation}
where $c_{B^{l+1}}^{l+1}$ is the coefficient associated to a basis function $B^{l+1}$.

The truncation of this basis function $B^{l}$, whose support overlaps with the support of finer basis functions $B^{l+1}$, part of $\mathcal{B}^{l+1}$ and defined on $\Omega^{l+1}$, is attained as follows (\cite{Giannelli2012a, Garau2016}):
\begin{equation}
\begin{array}{lll}
\trunc^{l+1}(B^{l}) & = &\displaystyle \sum_{
	\begin{array}{c}
	\scalebox{0.75}{$B^{l+1} \in \mathcal{B}^{l+1},$}\\ 
	\scalebox{0.75}{$\supp(B^{l+1}) \not\subseteq \Omega^{l+1}$}
	\end{array}} c_{B^{l+1}}^{l+1} \left( B^{l} \right)\ B^{l+1},\\[40pt]
& = & \displaystyle B^{l} - \sum_{\supp(B^{l+1}) \subseteq \Omega^{l+1}} c_{B^{l+1}}^{l+1} \left( B^{l} \right)\ B^{l+1}.\\
\end{array}
\label{Eq_truncatedBSpline}
\end{equation}

Following the creation of a hierarchical B-spline basis $\mathcal{H}$, a THB basis $\mathcal{T}$ is constructed recursively by considering the truncation in Eq.~\eqref{Eq_truncatedBSpline}, see \cite{Giannelli2012a,Garau2016}:
\begin{equation}
\left\{
\begin{array}{lll}
\mathcal{T}^0      & =: & \mathcal{B}^0\\[5pt]
\mathcal{T}^{l+1} & =: & \{ \trunc^{l+1}(B) \ |\ B \ \mbox{in} \ \mathcal{T}^{l} \wedge \supp(B) \not\subseteq \Omega^{l+1}\}\\[5pt]
&& \cup\ \{ B \in \mathcal{B}^{l+1} \ |\ \supp(B) \subseteq \Omega^{l+1}\}, \\[5pt]
&& \mbox{for}\ l = 0, \dots, n-2.
\end{array}
\right.
\label{Eq_truncatedRecursiveAlgorithm}
\end{equation}

The truncated basis $\mathcal{T}$ spans the same space as the non-truncated basis $\mathcal{H}$ and it admits a strong stability property (\cite{Giannelli2014}). Moreover, the smaller support results in a reduction in the number of nonzero basis functions per element and consequently a sparser system of linear equations in a finite element analysis. 

\begin{figure*}[ht]\center
	\includegraphics[width=.45\textwidth]{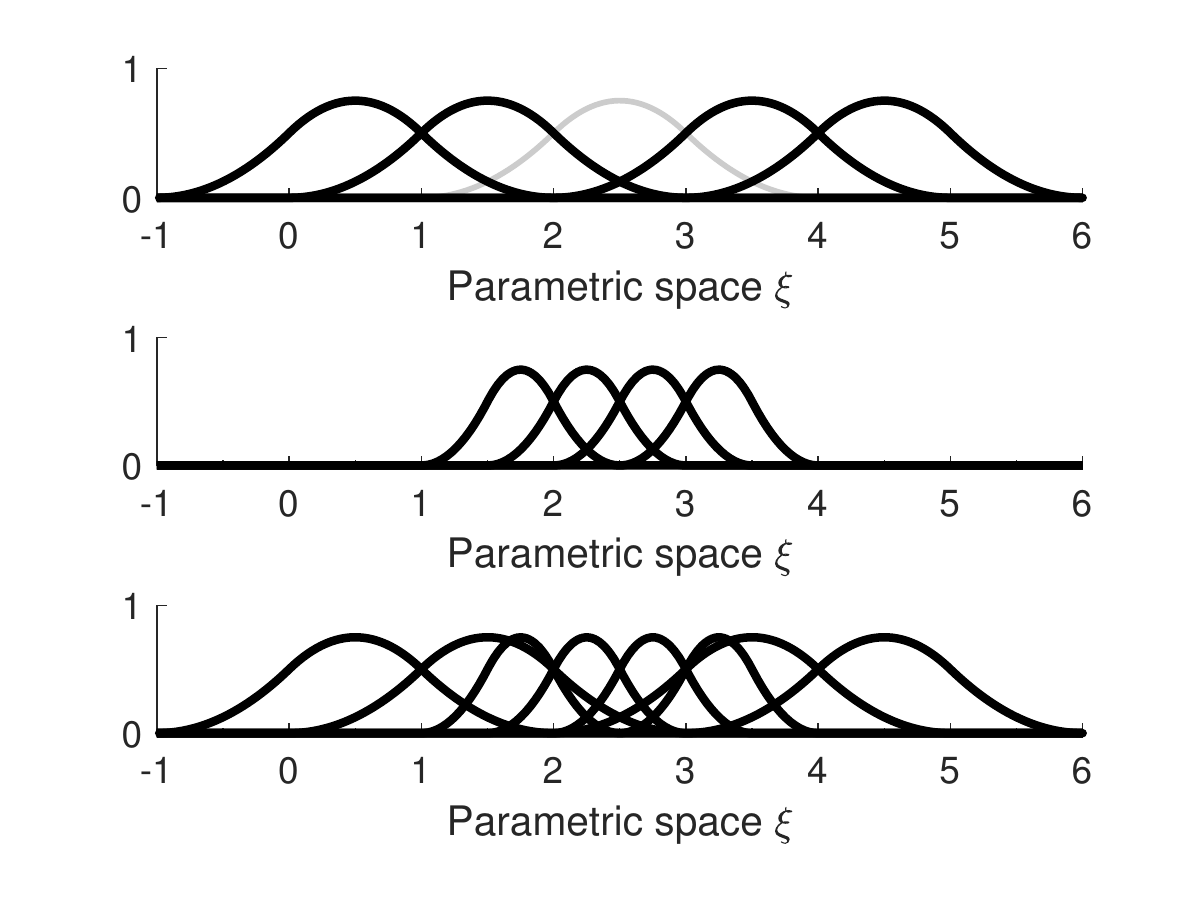}\hspace{0.5cm}
	\includegraphics[width=.45\textwidth]{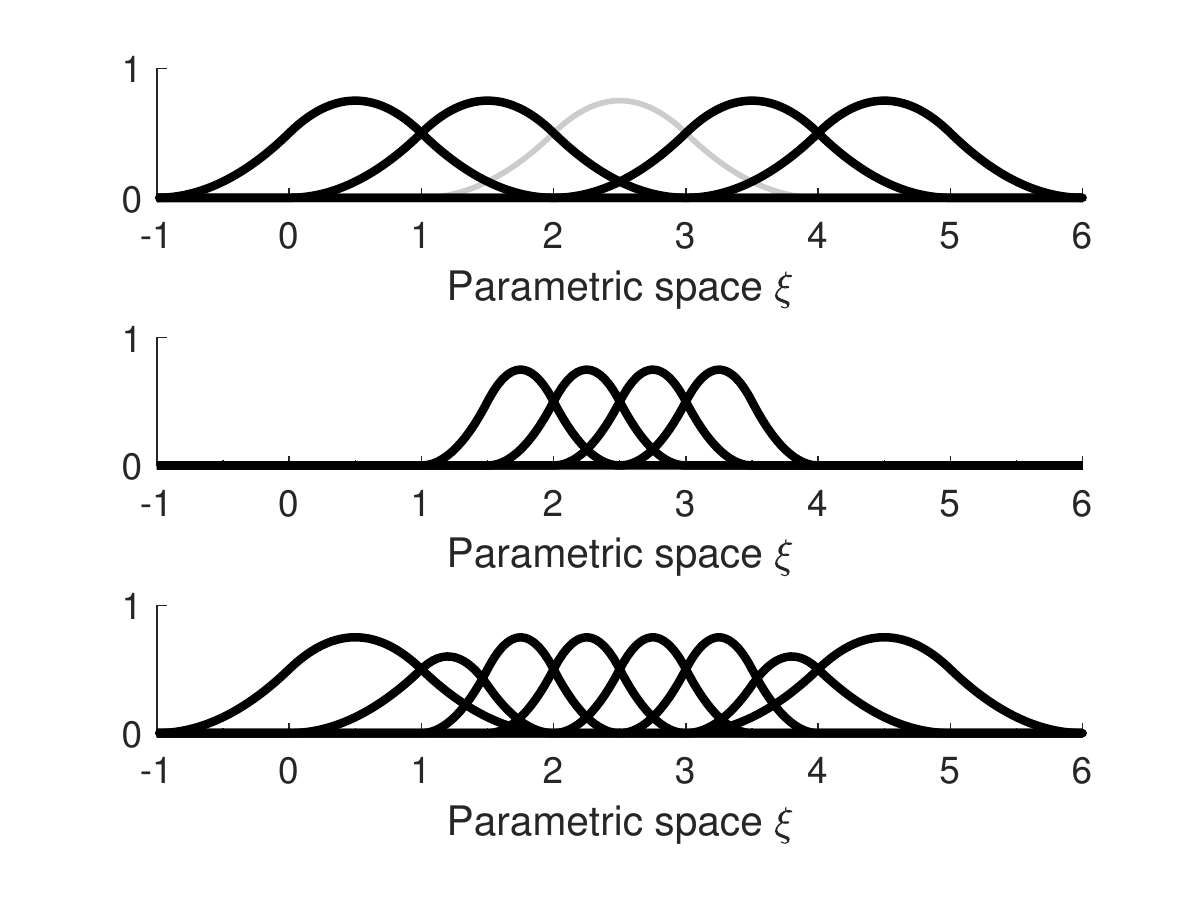}
	\caption{Comparison of univariate HB-spline (left) and THB-spline basis functions (right)}
	\label{fig_truncation1D}
\end{figure*}

The effect of the truncation is illustrated in Fig.~\ref{fig_truncation1D}. A univariate truncated and non-truncated basis $\mathcal{T}$ and $\mathcal{H}$ are juxtaposed. The first and second levels correspond to $\Omega^0$ and $\Omega^1$ respectively, while the bottom level represents the combination of the functions on these two levels. The comparison shows the reduced support of the truncated B-spline basis functions. 

\section{Implementation Details}\label{sec:details}

This section details the implementation of the THB discretization described above. A computationally efficient approach is introduced to represent hierarchically refined meshes via a poly-tree data structure. Local mesh refinement strategies are outlined. A methodology for expressing THB basis functions associated with multiple, differently refined meshes in terms of Lagrange shape functions over each element in a union background mesh is presented.

\subsection{Poly-tree Data Structure for Hierarchically Refined B-spline Discretizations}\label{Refinement pattern}

The generation of a large hierarchical refined mesh can be expensive due to the high memory consumption, especially when creating all possible elements for each refinement level. The cost is compounded when creating multiple, differently refined hierarchical meshes. In this paper, we present a computationally and memory efficient approach to build differently refined meshes.

We first construct a Poly-Tree (PT) data structure, i.e., a quadtree in 2D and an octree in 3D. The PT data structure represents a set of hierarchically refined meshes. The tree depth corresponds to a refinement level $l$ with the base level, $0$, representing a coarse uniform tensor mesh. The nodes of this PT data structure are labeled PT cells and used to construct discretizations of state variable fields. Note that the PT cells are not associated with any specific set of basis functions.

PT cells with higher refinement levels $l > 0$ are created recursively. Starting from level $0$, the PT cells are recursively subdivided into 4 and 8 PT cells in 2D and 3D, respectively. Considering a PT cell at a refinement level $l$, the PT cell at refinement level $l-1$ from which the PT cell is created is referred to as its parent. The PT cells at refinement level $l+1$ created by subdivision of a PT cell at level $l$ are called its children. The PT cells are only created once and only as needed.

 To efficiently represent multiple, differently refined meshes with the same PT data structure, we introduce the concept of PT cell activation states. Each PT cell has multiple activation states which are represented by Activation Indices (AI).
 Possible activation states for a particular AI are either \textit{active}, \textit{refined}, or \textit{inactive}. \textit{Active} PT cells for a specific AI have \textit{refined} parent PT cells and \textit{inactive} children PT cells. 

\begin{figure*}[ht]\centering
	\includegraphics[width=0.66\columnwidth]{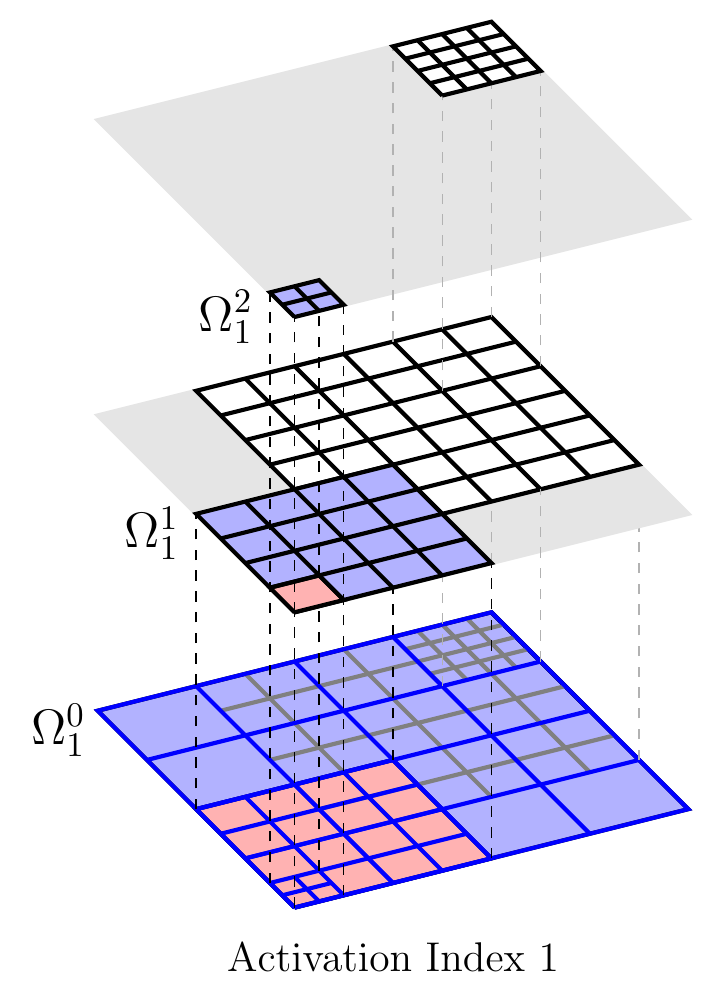}\hspace{0.25cm}
	\includegraphics[width=0.60\columnwidth]{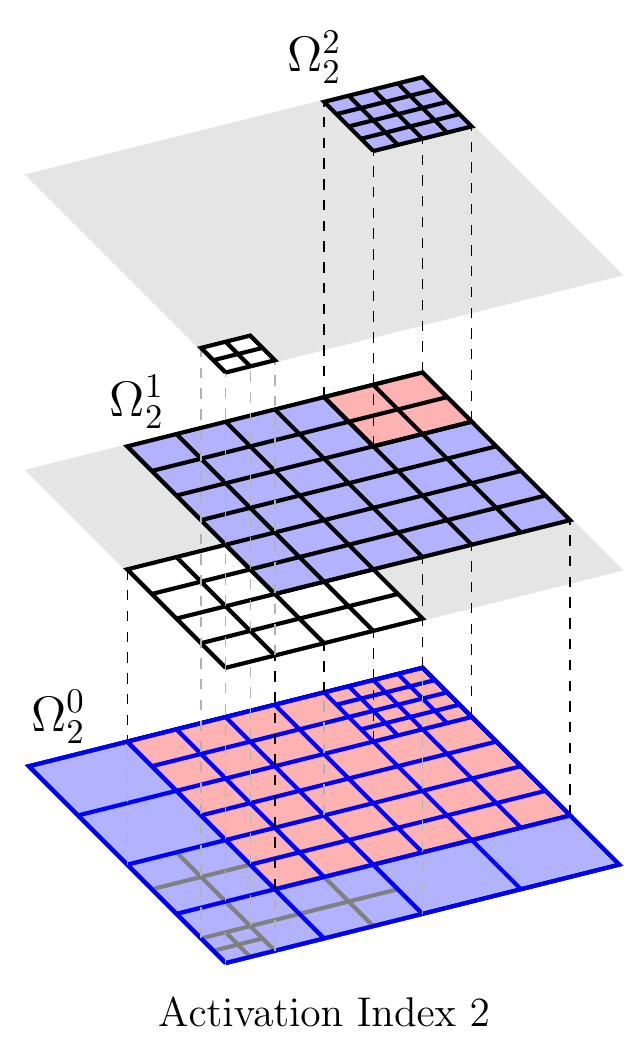}
	\caption{One hierarchical refined quad-tree data structure with two sets of AIs	}
	\label{fig_hierarchicalRefinedMesh}
\end{figure*}

Although PT cells are only constructed once, the activation state concept allows for the creation of multiple, differently refined meshes. A single activation state can be stored as a binary number with only two bits. Thus, the activation states for all AIs of a PT cell are stored in $2 \cdot \# \mathrm{AIs}$ bits. An illustration of a single quad-tree data structure with two AIs is presented in Fig.~\ref{fig_hierarchicalRefinedMesh}. The coloring of the PT cells indicates the activation state for a specific AI. For a particular AI, PT cells in blue, red, and white are \textit{active}, \textit{refined}, and \textit{inactive}, respectively. Blue framed PT cells on the zero level illustrate the resulting mesh of active PT cells on a given AI. This enables the efficient construction of differently refined hierarchical B-spline discretizations based on a single PT data structure and a set of AIs.  To see this, let $\Omega^l_I$ be the union of active and refined cells associated with level $l$ and AI $I$.  Then, we can create a hierarchical B-spline basis $\mathcal{H}_I$ and a corresponding hierarchically refined mesh $\mathcal{K}_I$ for AI $I$ from the sequence
\begin{equation}
\Omega^{n-1}_I \subseteq \Omega^{n-2}_I \subseteq \dots \subseteq \Omega^{0}_I = \Omega.
\end{equation}
The resulting hierarchically refined mesh $\mathcal{K}_I$ is precisely the set of active PT cells for AI $I$.

In our PT data structure, we only store PT cells that are either active or refined for at least one AI $I$.  The element edge length $h_{\mathrm{PT cell},m}^0$ for parametric direction $m$ and on the coarsest level $0$ is computed based on the given mesh size and the number of coarsest elements in each dimension, both of them are predefined by the user. The size of mesh elements corresponding to PT cells on levels $l > 0$ may be computed based on the coarsest element edge length and the refinement level via $h_{\mathrm{PT cells},m}^l = 2^{-l} \ h_{\mathrm{PT cell},m}^0$. Basis functions are assigned to mesh elements based on a pre-defined elemental basis function ordering. The elemental basis function ordering in this work follows the Exodus II standard (\cite{Shemon2014}).

\subsection{Local Refinement Strategies}\label{refStrategies}

To generate a PT data structure that supports different hierarchically refined meshes, the PT is recursively refined for each AI and corresponding refinement criteria. The refinement procedure is described in Algorithm~\ref{alg_refinment}. For a specific AI, \textit{active} PT cells are flagged for refinement based on chosen refinement criteria. Adjacent PT cells can also be flagged for refinement to increase the size of the refined region. Moreover, additional active PT cells for the current AI may be flagged considering mesh regularity requirements as discussed in Subsection~\ref{meshRegularity}. In case the PT data structure is generated in parallel using a domain decomposition strategy, flagged PT cells are communicated across adjacent subdomains to guarantee consistent refinement. Refinement is then performed by creating new children PT cells through subdividing all flagged PT cells unless the children PT cells do already exist. New PT cells are initialized with an \textit{inactive} activation state for all AIs. For the AI currently considered for refinement, the activation state of the children PT cells is set to \textit{active} while the parent PT cell is set to \textit{refined}.

\begin{algorithm}[H]
	\begin{algorithmic}[1] 
		\For{Activation Index (AIs)}
		\State{Flag PT cells based on refinement criteria for AI} 
		\State{Collect flagged PT cells into queue for refinement} 
		\While {changes in refinement queue }
		\For{all PT cells in queue for refinement}
		\State{Apply refinement buffer Algorithm~\ref{alg_refinementBuffer}}
		\State{Collect flagged PT cells in queue for refinement} 
		\EndFor
		\EndWhile
		\State{Communicate refinement queue}
		\If{Refined children PT cell does not exist}
		\State{Create refined children PT cells}
		\State{Initialize new PT cells as \textit{inactive} for all AIs}
		\EndIf
		\State{Flag children PT cells as \textit{active} for AI }
		\State{Flag parent PT cell as \textit{refined} for AI }
		\EndFor
	\end{algorithmic}	
	\caption{Refinement algorithm \label{alg_refinment}}
\end{algorithm}

\subsubsection{Mesh Regularity Requirements}\label{meshRegularity}

For construction of THB bases as discussed in Section~\ref{sec:bsplines}, mesh regularity requirements need to be considered when constructing the PT data structure. The difference in refinement level between adjacent PT cells in the refined PT is limited to one. Furthermore, all active neighbor PT cells inside a so-called buffer zone of an active PT cell on level $l$ must be of level greater or equal $l-1$, see Eq. (\ref{Eq_truncatedRecursiveAlgorithm}).

The buffer range $d_{\buffer,m}^l$ for a particular active PT cells is calculated by multiplying the PT cells size with a user-defined buffer parameter $b_{\buffer}$, i.e., $d_{\buffer,m}^l = b_{\buffer} \ h_{\mathrm{PT cells},m}^l$. When creating a B-spline basis, the width of the buffer zone in parametric direction $m$ must be larger than or equal to the width of the basis function supports in parametric dimension $m$. In this work each state variable field is interpolated with an individual interpolation order $p$. To satisfy the mesh regularity requirement for all interpolation functions, the buffer parameter must be chosen as $b_{\buffer} \geq p_{\mathrm{max}}$, where $p_{\mathrm{max}}$ is the maximal polynomial degree of all used bases. 

\begin{figure*}[ht]\centering
	\includegraphics[width=2.0\columnwidth]{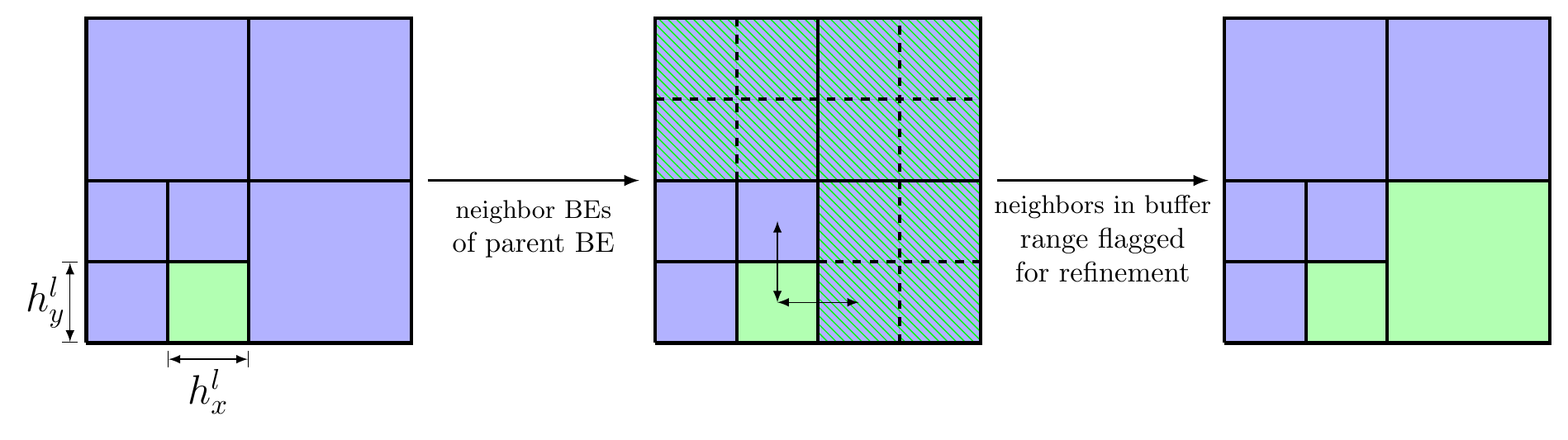}
	\caption{Visualization of refinement buffer Algorithm \ref{alg_refinementBuffer}}
	\label{fig_refBuffer}
\end{figure*} 

The refinement procedure for enforcing a buffer zone is summarized in Algorithm~\ref{alg_refinementBuffer}. The algorithm is applied to each PT cell initially flagged for refinement in Step 2 of Algorithm~\ref{alg_refinment} and starts by determining its parent PT cell. The refinement status of the parent's neighbors, i.e., cells within the buffer range of the considered parent, are checked. If these neighbors are neither refined nor flagged for refinement, the distance $d_{m}$ between the centers of the considered PT cell and its neighbor's children cells is calculated. This operation is trivial even if these neighboring children cells do not exist due to the PT data structure information. If any distance $d_{m}$ is smaller than the buffer range $d_{\buffer,m}^l$, the particular neighbor cell is flagged for refinement. The algorithm is then applied recursively to all newly flagged neighbor cells, until no further cells are flagged for refinement. An efficient access to hierarchical mesh information, such as neighborhood relationships, is provided by the PT data structure discussed in Subsection~\ref{Paralle_Considerations}. A visual representation of Algorithm~\ref{alg_refinementBuffer} is provided in Fig.~\ref{fig_refBuffer} for $b_{\buffer} = 1$. PT cells with a green fill are cells flagged for refinement while PT cells with a green pattern fill are neighboring PT cells of the parent cell.

\begin{algorithm}[ht]\center
	\begin{algorithmic}[1]
		\State{PT cell flagged for refinement in Algorithm~\ref{alg_refinment} }
		\State{Get parent}
		\State{Get parent's neighbors in buffer range $ d_{\buffer}^l$}
		\If{neighbor is active and not flagged for refinement}
		\State{Calculate distance $d_{\mathrm{ma}}$ considered cell and neighbors non-existing children PT cells}
		\If{$d_{\mathrm{m}} < d_{\buffer}$}
		\State{Flag neighbor for refinement}
		\State{Apply refinement buffer Algorithm~\ref{alg_refinementBuffer} for neighbor}
		\EndIf
		\EndIf
	\end{algorithmic}        
	\caption{Refinement buffer algorithm}
	\label{alg_refinementBuffer}
\end{algorithm}

\subsection{Union Background Mesh and Extraction Operators}\label{Projection}

The framework presented above allows for a different hierarchically refined B-spline discretization for each state variable field. These discretizations are defined on locally refined background meshes that are not aligned with the external boundaries and material interfaces.  To simplify the finite element formation and assembly process, a union\footnote{Note that the union background mesh is not attained via a set union of the separate hierarchically refined meshes.  Instead, elements in this set union that fully contain smaller elements in the set union do not belong to the union background mesh.} background mesh $\mathcal{K}_{union}$ is constructed from the hierarchically refined meshes for each state variable field as follows:
\begin{align}
\mathcal{K}_{union} := \left\{ K: \right. & K \in \mathcal{K}_I \text{ for some } I \text{ and there is no } \nonumber \\ & K' \in \mathcal{K}_J \text{ for } J \neq I \text{ such that } \left. K' \subsetneq K \right\}.
\end{align}
That is, an element of a hierarchically refined mesh for one of the state variable fields is an element of the union background mesh if it contains no finer element belonging to a hierarchically refined mesh for one of the other state variable fields.  An example of a union background mesh is displayed in Fig.~\ref{subfig-4:dummy}.  The union background mesh is specially defined so that the hierarchical B-spline basis functions associated with each state variable field are smooth over each element of the union background mesh.  The union background mesh  can be conveniently created via the PT data structure by combining all AIs used to generate THB background meshes.  In particular, a PT cell corresponds to an element of the union background mesh if it is active for at least one activation index and active or inactive for all other activation indices.  The union background mesh serves as the XFEM background mesh in which the geometry is immersed. Intersected XFEM background elements of the union background mesh are cut to generate a body-fitted integration mesh. Each THB discretization is enriched separately considering the B-spline interpolation and the immersed geometry represented on the union background mesh. For details on the enrichment strategy, the reader is referred to \cite{Noel2022}. To increase the geometry resolution, the union mesh can be further refined, either locally or globally.

\begin{figure}[ht]\centering
	\includegraphics[width=1.0\columnwidth]{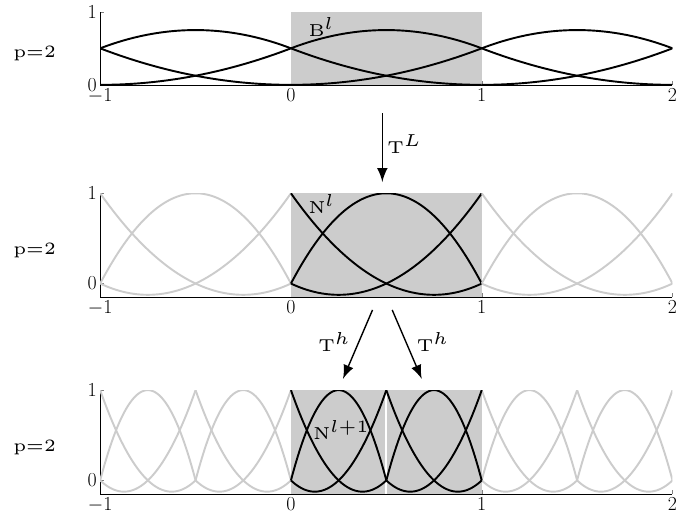}
	\caption{Illustration of the extraction process}
    \label{fig_Extraction}
\end{figure}

To facilitate operations performed on the union background mesh and on the integration mesh in the XFEM analysis, the THB basis functions for each locally refined background mesh are represented in terms of Lagrange shape functions on each element of the union background mesh via Lagrange extraction, see \cite{Remacle2012c,DAngella2018}.  Namely, the THB basis functions for a locally refined background mesh can be represented element-wise in terms of Lagrange shape functions as
\begin{equation}
B_k(\boldsymbol{\xi}) = \sum_{j} B_k(\boldsymbol{\xi}_j) N^L_j(\boldsymbol{\xi}) = \sum_{j} T^L_{jk} N^L_j(\boldsymbol{\xi}) 
\end{equation}
where $B_k$ is the $k^\text{th}$ THB basis function, $\left\{ N^L_j \right\}_j$ are the Lagrange shape functions over the element, $\left\{ \boldsymbol{\xi}_j \right\}_j$ are the locations at which the Lagrange shape functions $\left\{ N^L_j \right\}_j$ are interpolatory, and $T^L_{jk} = B_k(\boldsymbol{\xi}_j)$.  We refer to $T^L$ as a Lagrange extraction operator.  The Lagrange shape functions over each element can in turn be represented in terms of Lagrange shape functions over a child element as
\begin{equation}
N^L_j(\boldsymbol{\xi}) = \sum_{i} N^L_j(\hat{\boldsymbol{\xi}}_i) \hat{N}^L_i(\hat{\boldsymbol{\xi}}_i) = \sum_{i} T^h_{ij} \hat{N}^L_i(\hat{\boldsymbol{\xi}}_i)
\end{equation}
where $\left\{ \hat{N}^L_i \right\}_i$ are the Lagrange shape functions over the child element, $\hat{\boldsymbol{\xi}}_i$ are the locations at which the Lagrange shape functions $\left\{ \hat{N}^L_i \right\}_i$ are interpolatory, and $T^h_{ij} = N^L_j(\hat{\boldsymbol{\xi}}_i)$.  We refer to $T^h$ as an $h$-refinement extraction operator.  It follows that the THB basis functions can be expressed in terms of the Lagrange shape functions over the child element as
\begin{equation}
B_k(\boldsymbol{\xi}) = \sum_{i} T_{ik} \hat{N}^L_i(\boldsymbol{\xi}) 
\end{equation}
where
\begin{equation}
T_{ik} = \sum_{j} T^L_{jk} T^h_{ij} 
\end{equation}
is a Lagrange extraction operator that can be computed using the aforementioned Lagrange and $h$-refinement extraction operators.  This process can be repeated to represent THB basis functions in terms of Lagrange shape functions and easily computable Lagrange extraction operators over descendent elements of the child element as well, thus enabling us to represent the THB basis functions for each state variable field in terms of Lagrange shape functions and Lagrange extraction operators over each element of the union background mesh. The extraction process from a quadratic B-spline basis to a quadratic, once refined Lagrange basis is illustrated in Fig.~\ref{fig_Extraction}.  Our framework limits the background element refinement to a factor of two, as presented in Subsection~\ref{meshRegularity}. This consequently limits the number of \textit{h}-refinement extraction matrices to four in 2D and eight in 3D. These matrices can be precomputed and efficiently selected exploiting the PT data structure.

\subsection{Background Mesh Data Structure}\label{Paralle_Considerations}

In this subsection, we discuss computational aspects of the PT data structure for storing and managing hierarchically refined meshes. This data structure aims at improving the overall computational efficiency and minimize inter-process communication for parallel computations. The presented implementation is limited to quadtrees in 2D and octrees in 3D, which is sufficient for IFEMs. 

Given the dimensions of the computational domain and the number of elements in each spatial direction, the base level $0$ of the PT data structure is generated. Adopting a domain decomposition approach, the PT cells are grouped into subdomains. The mesh is hierarchically refined, and the union mesh is created in parallel on each processor. To facilitate the refinement process and the construction of the Lagrange extraction operators, we create overlapping subdomains to build an efficient PT data structure in parallel. In many cases a decomposition strategy assigning approximately the same number of PT cells to each subdomain is sufficient. However, as refinement can lead to a significant imbalance in number of PT cells across subdomains, it may be beneficial to choose a decomposition strategy that accounts for refinement.

\begin{figure}[ht]\center
	\includegraphics[width=0.9\columnwidth]{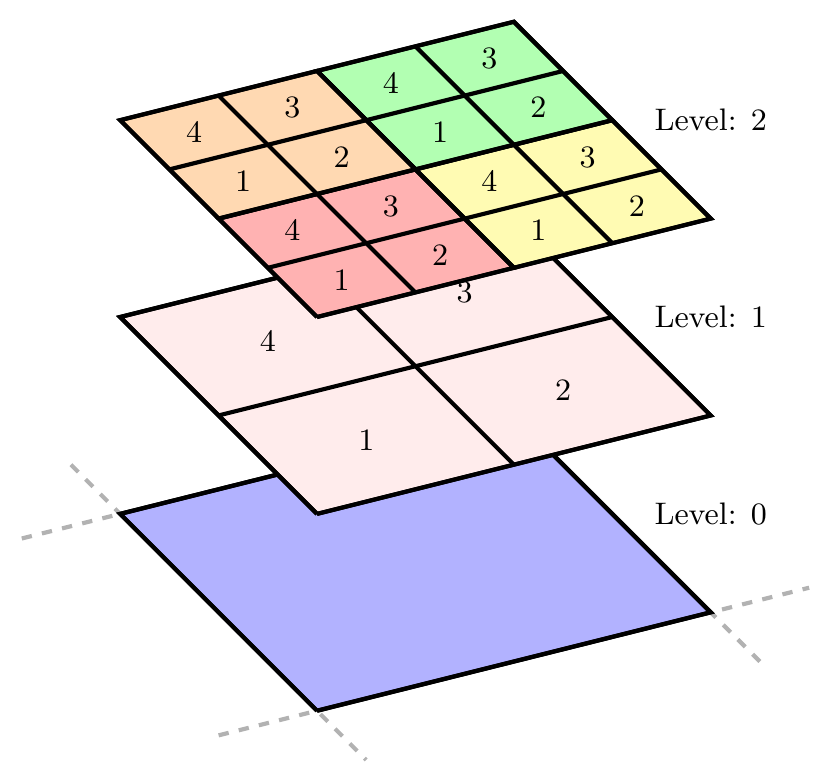}
	\caption{A coarse two dimensional PT cell refined through recursive subdivision and resulting quadtree data structure}
	\label{fig_Octree}
\end{figure}

\begin{figure}[ht]\center
	\includegraphics[width=1.0\columnwidth]{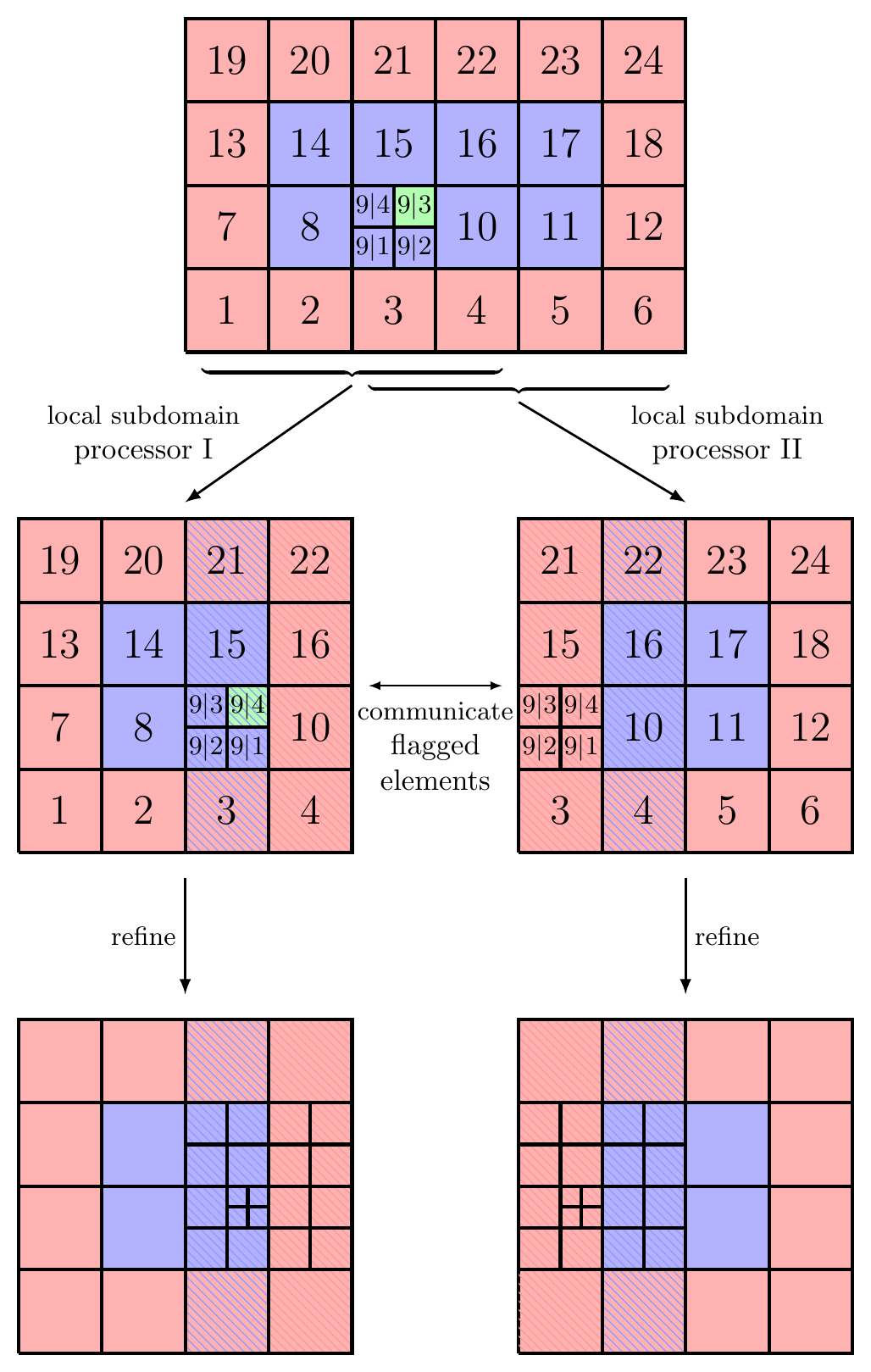}
	\caption{Parallel local refinement of a global $2 \times 4$ domain}
	\label{fig_parallel_refinement}
\end{figure}

At the coarsest level, $0$, a unique PT cell ID can be calculated based on the global PT cell location as illustrated in the uppermost graphic of Fig.~\ref{fig_parallel_refinement}. With the ID of the parent PT cell at refinement level $0$ and the position in the PT structure, unique IDs of refined PT cells can be determined directly, without inter-process communication and building local-to-global ID maps. The positions of the PT structure for a two dimensional twice refined PT cells are illustrated in Fig.~\ref{fig_Octree}. 

Furthermore, the subdomain-local PT data structure allows for efficient access to cell neighborhood and hierarchy relationships. The PT data structure also speeds up the construction of extraction operators and may be used for the identification of basis functions. Basis functions can be uniquely identified by utilizing the PT structure in combination with the elemental Exodus II basis function index. 

To reduce the communication across subdomains, an aura of PT cells is constructed around the set of PT cells owned by a specific processor, leading to overlapping subdomains. The width of the aura is chosen to be $b_{\buffer}$ PT cells on the coarsest refinement level with $b_{\buffer}$ defined as described in Subsection~\ref{meshRegularity}.  

The PT cells in an aura are needed to define all THB basis functions that are nonzero over non-aura elements. Therefore, most computations can be performed on a subdomain level, interlaced with only a few inter-subdomain communications. This includes mesh refinement along the processor boundaries, the construction of Lagrange extraction operators, and the evaluation of face-oriented ghost stabilization residuals and Jacobians, see Section~\ref{sec:model}. To illustrate the parallel concept, Fig. \ref{fig_parallel_refinement} presents a global two dimensional domain of size $2 \times 4$ elements in blue. PT cells in red indicate aura cells. The global domain is decomposed into two processor local subdomains. Unique global PT cell IDs are calculated based on IDs of the coarsest refinement level. Red striped PT cells are aura cells shared with the neighboring processor. Blue striped PT cells are cells of the inverse aura. The green cells are flagged for refinement. Parallel consistent refinement is guaranteed through communication of the flagged PT cell IDs on aura and inverse aura.

\section{Thermo-Elastic Model}\label{sec:model}

The XIGA framework presented above is applicable to a broad range of physical problems that are described by partial differential equations. In this paper, we demonstrate the functionality of the developed framework with static elastic and thermo-elastic problems. The elastic model assumes infinitesimal strains and a linear elastic, isotropic material response. The thermal model accounts for linear diffusive heat transfer. The structural response depends on the temperature field through inelastic, isotropic thermal expansion. In this section, we present the variational form of the stabilized governing equations of the thermo-elastic model.

\subsection{Governing Equations}\label{governingEquation}

The weak form of the residual is decomposed into the following four contributions:
\begin{equation}
\mathcal{R} = \mathcal{R}^{\T}+\mathcal{R}^{\U} + \mathcal{R}^{\N} + \mathcal{R}^{\G} = 0,
\label{Eq_RFull}
\end{equation} 
where $\mathcal{R}^{\T}$ and $\mathcal{R}^{\U}$ combine the volumetric contributions and contributions of Neumann boundary conditions to the thermal and structural residuals, respectively. The contributions of Nitsche's formulation for Dirichlet boundary conditions are collected in $\mathcal{R}^{\N}$ and the contributions of the face-oriented ghost stabilization in $\mathcal{R}^{\G}$.

In the absence of body loads, the thermal and structural residual, $\mathcal{R}^{\T}$ and $\mathcal{R}^{\U}$, are defined over all material domains $\Omega^m$ with $m=1, \ldots, N_m$ as follows:
\begin{equation}
\displaystyle \mathcal{R}^{\T} = 
\sum_{m=1}^{N_m} {
\int_{\Omega^m}\ \nabla \delta \mathbf{T} \ \boldsymbol{\kappa} \ \nabla \mathbf{T} \ d\Omega - \int_{\Gamma_{N,T}^m} \delta \mathbf{T} \ q_N \ d\Gamma },
\label{Eq_TLin}
\end{equation}
and
\begin{equation}
\displaystyle \mathcal{R}^{\U} = 
\sum_{m=1}^{N_m} {
\int_{\Omega^m}\  \delta \boldsymbol{\varepsilon}_u : \boldsymbol{\sigma} \ d\Omega 
- \int_{\Gamma_{N,u}^m} \delta \mathbf{u} \cdot \mathbf{f}_{N} \ d\Gamma},
\label{Eq_ULin}
\end{equation}
where $\mathbf{T}$ and $\delta \mathbf{T}$ are the temperature trial and test function, respectively. The isotropic thermal conductivity tensor is denoted by $\boldsymbol{\kappa}$. A heat flux ${q}_{N}$ is applied at the boundaries $\Gamma_{N,T}^m$. The displacement trial and test functions are denoted by $\mathbf{u}$ and $\delta \mathbf{u}$, respectively. Traction forces, $\mathbf{f}_{N}$, are applied on the boundaries $\Gamma_{N,u}^m$. The Cauchy stress tensor is denoted by $\boldsymbol{\sigma} = \mathbf{D}\ \boldsymbol{\varepsilon}_e$ where $\mathbf{D}$ is the isotropic constitutive tensor and $\boldsymbol{\varepsilon}_e$ is the elastic strain tensor with $\boldsymbol{\varepsilon}_e = \boldsymbol{\varepsilon}_{u} - \boldsymbol{\varepsilon}_{T}$.  The total infinitesimal strain  $\boldsymbol{\varepsilon}_{u}$ is computed by $\boldsymbol{\varepsilon}_{u} = \frac{1}{2} \left( \nabla \mathbf{u} + \nabla \mathbf{u}^T \right)$. The thermal strain is denoted by $\boldsymbol{\varepsilon}_{T} =  \alpha  \left( \mathbf{T} -\mathbf{T}_0 \right) \mathbf{I}$ where $\alpha$ is the CTE and $T_0$ is the reference temperature.

To weakly enforce Dirichlet boundary and interface conditions, we use an unsymmetric formulation of Nitsche's method. The associated residual, $\mathcal{R}^{\N}$, is decomposed into boundary and interface terms, separately for the thermal and elastic subproblems as follows:
\begin{equation}
\mathcal{R}^{\N}=\mathcal{R}^{\N,D}_{\mathbf{T}} + \mathcal{R}^{\N,D}_{\mathbf{u}} + \mathcal{R}^{\N,I}_{\mathbf{T}} + \mathcal{R}^{\N,I}_{\mathbf{u}},
\end{equation}
where the thermal and elastic contributions from Dirichlet boundary conditions are denoted by $\mathcal{R}^{\N,D}_{\mathbf{T}}$ and $\mathcal{R}^{\N,D}_{\mathbf{u}}$, respectively. The thermal and elastic contributions from interface conditions are $\mathcal{R}^{\N,I}_{\mathbf{T}}$ and $\mathcal{R}^{\N,I}_{\mathbf{u}}$, respectively. These Dirichlet boundary residual contributions are defined as follows:
\begin{equation}
\begin{array}{ll}
\mathcal{R}^{\N,D}_T =  \sum \limits_{m=1}^{N_m} 
& - \int_{\Gamma^{D,T}} \delta \mathbf{T} \  (\boldsymbol{\kappa} \ \nabla \mathbf{T}) \cdot  \mathbf{n}_{\Gamma} \ d\Gamma \\[8pt]
& + \int_{\Gamma^{D,T}} (\boldsymbol{\kappa} \ \nabla \delta \mathbf{T}) \cdot \mathbf{n}_{\Gamma} \ ( \mathbf{T} - \mathbf{T}_D ) \ d\Gamma \\[0pt]
& + \sum \limits_{e=1}^{N_{e}^{B}} \gamma_{D,T}^e \ \int_{\Gamma^{D,T}_e} \delta \mathbf{T} \ ( \mathbf{T} - \mathbf{T}_D ) \ d\Gamma ,
\end{array}
\label{Eq_RNitscheTD}
\end{equation}
and
\begin{equation}
\begin{array}{ll}
\mathcal{R}^{\N,D}_{u} =  \sum \limits_{m=1}^{N_m} 
& - \int_{\Gamma^{D,u}} \delta \mathbf{u} \cdot  \boldsymbol{\sigma}(\mathbf{u}) \ \mathbf{n}_{\Gamma} \ d\Gamma \\[8pt]
& + \int_{\Gamma^{D,u}} \boldsymbol{\sigma}(\mathbf{\delta u}) \ \mathbf{n}_{\Gamma} \cdot ( \mathbf{u} - \mathbf{u}_D ) \ d\Gamma \\[0pt]
& + \sum \limits_{e=1}^{N_{e}^{B}} \gamma_{D,u}^e \ \int_{\Gamma^{D,u}_e} \delta \mathbf{u} \cdot ( \mathbf{u} - \mathbf{u}_D ) \ d\Gamma .
\end{array}
\label{Eq_RNitscheD}
\end{equation}
The temperature $T_D$ is prescribed on the boundary $\Gamma^{D,T}$ and the displacement $\mathbf{u}_D$ is prescribed on $\Gamma^{D,u}$. The vector $\mathbf{n}_{\Gamma}$ denotes the outward pointing normal on the boundary. The summations over all $N_e^B$ elements in the union background mesh in Eqs.~\eqref{Eq_RNitscheTD} and \eqref{Eq_RNitscheD} penalize constraint violations along the boundaries $\Gamma_{D,T}^e$ and $\Gamma_{D,u}^e$ within the elements of the union background mesh.  We henceforth refer to the elements of the union background mesh as background elements. The penalty parameters $\gamma_{D,T}^e$ and $\gamma_{D,u}^e$ depend on the size of the background element, $h_e^B$, and are defined as follows:
\begin{equation}
\gamma_{D,T} = c_{D,T} \frac{\kappa}{h_e^B} \quad \textrm{and} \quad \gamma_{D,u} = c_{D,u} \frac{E}{h_e^B} ,
\end{equation}
where $\kappa$ is the isotropic material conductive and $E$ is the Young's modulus of the linear elastic material. The parameters $c_{D,T} \geq 0$ and $c_{D,u} \geq 0$ control the accuracy of enforcing the Dirichlet boundary conditions. 

Continuity of temperature and displacement fields and balance of heat flux and traction must be satisfied at all interfaces $\Gamma^{m,n} = \Omega^m \cap \Omega^n \neq \emptyset$. 
Nitsche's method for the thermal and structural interface conditions yields the following residual contributions: 
\begin{equation}
\begin{array}{ll}
\mathcal{R}^{\N,I}_T = \sum \limits_{\Gamma^{m,n}} \sum \limits_{e=1}^{N_e^B} 
& - \int_{\Gamma^{m,n}_e} \llbracket \delta \mathbf{T} \ \rrbracket \ \{ \boldsymbol{\kappa} \ \nabla \mathbf{T}  \} \cdot \mathbf{n}_{m,n} \ d\Gamma \\[6pt]
& + \int_{\Gamma^{m,n}_e} \{ \boldsymbol{\kappa} \ \nabla \delta \mathbf{T} \} \cdot \mathbf{n}_{m,n} \ \llbracket \mathbf{T} \rrbracket \ d\Gamma \\[6pt]
& + \gamma_{I,T}^e \ \int_{\Gamma^{m,n}_e} \llbracket \delta \mathbf{T} \rrbracket \ \llbracket \mathbf{T} \rrbracket \ d\Gamma,
\end{array}
\label{Eq_RNitscheTI}
\end{equation}
\begin{equation}
\begin{array}{ll}
\mathcal{R}^{\N,I}_u  = \sum \limits_{\Gamma^{m,n}} \sum \limits_{e=1}^{N_{e}^{B}} 
 & - \int_{\Gamma^{m,n}_e} \llbracket \delta \mathbf{u} \rrbracket \cdot \{ \boldsymbol{\sigma}(\mathbf{u}) \} \ \mathbf{n}_{m,n} \ d\Gamma \\[6pt]
 & + \int_{\Gamma^{m,n}_e} \{ \boldsymbol{\sigma}(\delta \mathbf{u}) \}  \ \mathbf{n}_{m,n}  \cdot \llbracket \mathbf{u} \rrbracket \ d\Gamma  \\[6pt]
 & + \gamma_{I,u}^e \ \int_{\Gamma^{m,n}_e} \llbracket \delta \mathbf{u}\rrbracket \cdot \llbracket \mathbf{u} \rrbracket d\Gamma,
\end{array}
\label{Eq_RNitscheI}
\end{equation}
where the jump operator is defined as $\llbracket \cdot \rrbracket = ( \cdot )^{m}-( \cdot )^{n}$. The numerical interface flux and traction are defined by the averaging operator as $\{ \cdot \} =w^{m} ( \cdot )^{m}+ w^{n}( \cdot )^{n}$, where $w^m$ and $w^n$ are weights. The vector $\mathbf{n}_{m,n}$ denotes the normal vector on the interface pointing from phase $m$ to phase $n$. The accuracy of enforcing the interface condition is controlled by the penalty terms with $\gamma_{I,T}^e$ and $\gamma_{I,u}^e$ being the elemental penalty factors. 

We follow the work of \cite{Annavarapu2012} and define weights for the numerical heat flux as follows:
\begin{equation}
\begin{array}{ll}
w_T^{m}  = & \frac{\text{\text{meas}}(\Omega^{m}) / \kappa^{m} }{\text{meas}(\Omega^{m})/\kappa^{m}+\text{meas}(\Omega^{n})/\kappa^{n}} , \\[6pt]
w_T^{n}  = & \frac{\text{\text{meas}}(\Omega^{n}) / \kappa^{n} }{\text{meas}(\Omega^{m})/\kappa^{m}+\text{meas}(\Omega^{n})/\kappa^{n}} ,
\end{array}
\end{equation}
and for the numerical traction as follows:
\begin{equation}
\begin{array}{ll}
w_u^{m}  = & \frac{\text{\text{meas}}(\Omega^{m}) / E^{m} }{\text{meas}(\Omega^{m})/E^{m}+\text{meas}(\Omega^{n})/E^{n}} , \\[6pt]
w_u^{n}  = & \frac{\text{\text{meas}}(\Omega^{n}) / E^{n} }{\text{meas}(\Omega^{m})/E^{m}+\text{meas}(\Omega^{n})/E^{n}} ,
\end{array}
\end{equation}
where $\text{meas}( \Omega^{k} )$ is the surface in 2D or the volume in 3D of the domain occupied by the phase $k$ within the background element. The elemental penalty factors in Eqs.~\eqref{Eq_RNitscheTI} and \eqref{Eq_RNitscheI} are defined by:
\begin{equation}
\begin{array}{ll}
\gamma_{I,T}^e = & 2 \ c_{I,T} \ \frac{\text{meas}(\Gamma^{m,n})}{\text{meas}(\Omega^{m}) /\kappa^{m} +\text{meas}(\Omega^{n})/\kappa^{n}} \\[6pt]
\gamma_{I,T}^e = & 2 \ c_{I,u} \ \frac{\text{meas}(\Gamma^{m,n})}{\text{meas}(\Omega^{m}) /E^{m} +\text{meas}(\Omega^{n})/E^{n}},
\end{array}
\end{equation}
where the operator $\text{meas}(\Gamma^{m,n} )$ measures the length in 2D or the area in 3D of the interface within the background element. The parameters $c_{I,T} \geq 0$ and $c_{I,u} \geq 0$ control the accuracy of enforcing the interface conditions.

\subsection{Face-oriented ghost stabilization}\label{subsectionStabilization}

Face-oriented ghost stabilization is used to mitigate numerical instabilities caused by basis functions with small support within the geometric domain. This may occur when an interface moves close to the boundary of the support of a basis function. Such configurations may result in ill-conditioning of the system of linear equations and in imprecise spatial gradients of the state variable field, see \cite{DePrenter2017}. This work adopts the face-oriented ghost stabilization approach presented by \cite{Burman2014a} and adapted by \cite{Noel2022} to fit the basis function enrichment strategy described in Section \ref{sec:xfem}.
 
Let $\Omega^m$ be a material subdomain phase $m$ and $K_{\Omega^m}$ the set of background elements that have a non-empty intersection with $\Omega^m$: 
\begin{equation}
	K_{\Omega^m} := \left\{ K \in \mathcal{K}_{union} : K \cap \Omega^m \neq \emptyset \right\}.
	\label{Eq_ghost_K}
\end{equation}
We define $\mathcal{F}^m_{int}$ as the set of interior facets of $K_{\Omega^m}$, i.e., the facets $F$ shared between two background elements $\Omega_{F}^{m,+}$ and $\Omega_{F}^{m,-}$ of $K_{\Omega^m}$. Let $\tilde{\Gamma}$ be the union of all material interfaces and geometric boundaries, here defined as material interfaces between void and non-void regions. The set of ghost facets for phase $m$ is:
\begin{equation}
	\mathcal{F}^m_{ghost} := \left\{ F \in \mathcal{F}^m_{int} : \Omega_{F}^{m,+} \cap \tilde{\Gamma} \neq \emptyset\ \mbox{or} \ \Omega_{F}^{m,-} \cap \tilde{\Gamma} \neq \emptyset \right\}.
	\label{Eq_ghost_F}
\end{equation}
Consider a ghost facet $F^m$ shared between two adjacent background elements $\Omega_{F}^{m,+}$ and $\Omega_F^{m,-}$. The normal to the facet is $\mathbf{n}_{F}$ and is chosen as $\mathbf{n}_{F} = \mathbf{n}_{F}^{m,+} = - \mathbf{n}_{F}^{m,-}$. The material layout subdivides the element $\Omega_{F}^{m,+}$ into $N_{F}^{m,+}$ connected subdomains $\Omega_{F,i}^{m,+}$ with $i=1, \dots, N_{F}^{m,+}$ and the element $\Omega_{F}^{m,-}$ into $N_{F}^{m,-}$ connected subdomains $\Omega_{F,j}^{m,-}$ with $j=1, \dots, N_{F}^{m,-}$.
 
\begin{figure}[t!]\centering
	\includegraphics[width=0.49\columnwidth]{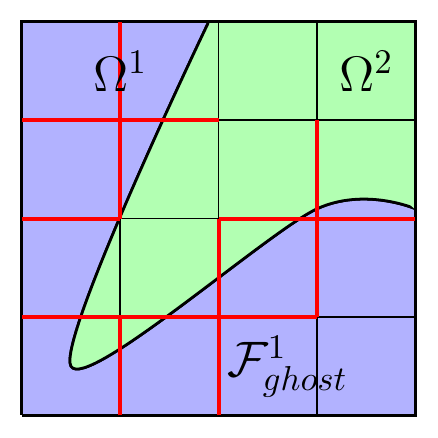}
	\includegraphics[width=0.49\columnwidth]{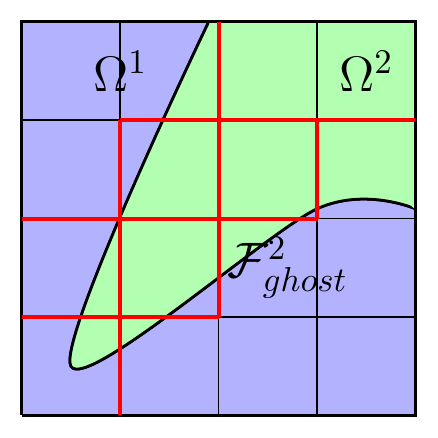}
	\caption{Set of ghost facets $\mathcal{F}_{ghost}^1$ and $\mathcal{F}_{ghost}^2$ used for the face-oriented ghost stabilization for a two-material problem}
	\label{fig_Ghost1}
\end{figure}

We define $\mathbf{u}_{F,i}^{m,+}$ as the polynomial extension of the field $\mathbf{u}|_{\Omega_{F,i}^{m,+}}$ to all of $\mathbb{R}^{d}$ and let $\mathbf{u}_{F,j}^{m,-}$ be the polynomial extension of the field $\mathbf{u}|_{\Omega_{F,j}^{m,-}}$ to all of $\mathbb{R}^{d}$. 

Additionally, our formulation requires that:
$| \partial \Omega_{F,i}^{m,+}\, \cap\, \partial \Omega_{F,j}^{m,-} | \neq 0,$
and the ghost stabilization is only applied between $\mathbf{u}_{F,i}^{m,+}$ and $\mathbf{u}_{F,j}^{m,-}$ when the boundaries of $\Omega_{F,i}^{m,+}$ and $\Omega_{F,j}^{m,-}$, $\partial \Omega_{F,i}^{m,+}$ and $\partial \Omega_{F,j}^{m,-}$, respectively, meet along a portion of the facet $F$ with a non-zero measure, e.g., the boundaries meet along more than a point in two dimensions and along more than a line in three dimensions.

With the above terminology defined, the contribution of the ghost stabilization for the displacement field to the residual equations is:
%
\begin{equation}
	\begin{array}{ll}
	    \mathcal{R}^{\G}_{\mathbf{u}} = 
	    \sum \limits_{m=1}^{N_m} \sum \limits_{F \in \mathcal{F}_{ghost}^m} &
		\sum \limits_{i=1}^{N_{F}^{m,+}} 
		\sum \limits_{\mathcal{J}_{\mathcal{F},i}^{m}}\\ [6pt]
		& \left[ \int_{F}  
		\gamma^{\mathbf{u}}_{G} \ h^{\tilde{k}} 
		\llbracket  \partial^{p}_{n}\mathbf{u} \rrbracket 
		\llbracket  \partial^{p}_{n}\mathbf{u} \rrbracket 
		d\Gamma \right],
	\end{array}
	\label{Eq_RGhostU2}
\end{equation}
where the set $\mathcal{I}_{\mathcal{F},i}^{m}$ is defined as:
\begin{equation}
	\mathcal{J}_{\mathcal{F},i}^m:= \{k \in \{1 \ldots N_{F}^{m,-} \} : | \partial \Omega_{F,i}^{m,+}\, \cap\, \partial \Omega_{F,j}^{m,-} | \neq 0 \},
	\label{Eq_RGhostU_Set}
\end{equation}
the jump operator $\llbracket \cdot \rrbracket$ is defined as:
\begin{equation}
	\llbracket  \partial^{p}_{n}\mathbf{u} \rrbracket =  \left( \partial^{p}_{n} \mathbf{u}_{F,i}^{m,+} - \partial^{p}_{n} \mathbf{u}_{F,j}^{m,-} \right),
	\label{Eq_JumpU}
\end{equation}
and
\begin{equation}
	\llbracket \partial^{p}_{n} \delta \mathbf{u} \rrbracket= \left( \partial^{p}_{n} \delta \mathbf{u}_{F,i}^{m,+} - \partial^{p}_{n} \delta \mathbf{u}_{F,j}^{m,-} \right),
	\label{Eq_JumpDU}
\end{equation}
where $\partial_n^k(\cdot)$ is the $k^{th}$ order normal derivative operator $\partial_n^k(\cdot) = \nabla^k(\cdot) \cdot \mathbf{n}_F$ where $\nabla^{k}(\cdot)$ is the $k^{th}$ order spatial gradient. The parameter $p$ is the polynomial order of the approximation, i.e., the THB discretization. It should be noted that only the $p^{th}$ contribution is nontrivial for the THB discretization used in this work. The ghost penalty parameter is denoted $\gamma_{G}^{\mathbf{u}}$ and is defined as a multiple of the Young's modulus $E$ of the considered material. 

The ghost penalization for the temperature field is defined similarly as:
\begin{equation}
\begin{array}{ll}
	 \mathcal{R}^{\G}_{\mathbf{T}} = 
	 \sum \limits_{m=1}^{N_m}
		\sum \limits_{F \in \mathcal{F}_{ghost}^m} &
	\sum \limits_{i=1}^{N_{F}^{m,+}} 
	\sum \limits_{\mathcal{J}_{\mathcal{F},i}^{m}}\\ [6pt] 
	&\left[ \int_{F}  
	\gamma^{\mathbf{T}}_{G}\ h^{\tilde{k}} 
	\llbracket  \partial^{p}_{n}\mathbf{T} \rrbracket 
	\llbracket  \partial^{p}_{n}\mathbf{T} \rrbracket 
	d\Gamma \right],
\end{array}
	\label{Eq_RGhostTheta}
\end{equation}
where the jump operator is defined as:
\begin{equation}
	\llbracket \partial^{p}_{n} \delta \mathbf{T} \rrbracket= \left( \partial^{p}_{n} \delta \mathbf{T}_{F,i}^{m,+} - \partial^{p}_{n} \delta \mathbf{T}_{F,j}^{m,-} \right),
	\label{Eq_JumpT}
\end{equation}
and
\begin{equation}
	\llbracket  \partial^{p}_{n}\mathbf{T} \rrbracket =  \left( \partial^{p}_{n} \mathbf{T}_{F,i}^{m,+} - \partial^{p}_{n} \mathbf{T}_{F,j}^{m,-} \right).
	\label{Eq_JumpDT}
\end{equation}
The ghost penalty parameters, $\gamma^{\mathbf{u}}_{G}$ and $\gamma_{G}^{\mathbf{T}}$, are defined as multiples of the Young's modulus $E$ and the conductivity $\kappa$ of the considered material, respectively. The parameter $\tilde{k}$ is defined as $\tilde{k} = 2(p-1)+1$ and enables control over the influence of the ghost penalty term. Due to the application of $C^{p-1}$ continuous THB bases, only jumps in gradients of order $p$ must be penalized. Further details on how to choose the ghost penalty term can be found in \cite{Burman2014a}. In this work we commonly choose a penalty parameter $\gamma_{G} = 0.001$.

\section{Numerical Examples}\label{numericalExamples}

In this section, we present 2D and 3D examples which illustrate the basic concepts and computational performance of the proposed discretization framework. Canonical examples show that our approach yields optimal convergence rates for linear, quadratic, and cubic B-spline approximations. An example with a stress singularity illustrates the computational advantages of locally  refined discretizations. 
Finally, the analysis of a polycrystalline micro-structure  demonstrates the applicability of our framework to complex multi-material problems.

To quantify the accuracy of the XFEM analyses, we define the error of a generic vector state variable field $\mathbf{v}$ in the $L^{2}$ norm and $H^1$ semi-norm as:
\begin{equation}
||\mathbf{v} - \mathbf{v}^h ||_{L^2} = \sqrt{ \sum_{m=1}^{N_m} \int_{\Omega^m} \left| \mathbf{v} - \mathbf{v}^h \right|^2 d \Omega },
\label{Eq::L2Norm}
\end{equation}
and
\begin{equation}
| \mathbf{v} - \mathbf{v}^h |_{H^1} = \sqrt{ \sum_{m=1}^{N_m} \int_{\Omega^m} \left| \nabla \mathbf{v} - \nabla \mathbf{v}^h \right|^2 d \Omega },
\label{Eq::H1SemiNorm}
\end{equation}
where $\mathbf{v}^h$ is the approximate field and $\mathbf{v}$ is the reference solution which is either the analytical solution if available or a solution computed on a sufficiently refined discretization.

The systems of discretized governing equations are solved by the direct solver PARDISO for 2D problems, see \cite{Kourounis2018}. A Generalized Minimal Residual (GMRES) method in combination with a dual threshold incomplete LU factorization with a degree of fill of 5.0 is used for 3D problems, see \cite{Saad2003}. The GMRES iterations are terminated if a relative drop of $10^{-10}$ of the normalized linear residual is achieved. All geometric and material parameters are given with each example in self-consistent units unless stated otherwise.

In the following subsections, we first present examples of single- and multi-material problems considering either a thermal or mechanical response. These examples characterize the fundamental features of the proposed analysis framework within single-physics settings. The last two examples consider coupled thermo-elastic single- and multi-material problems.

\subsection{Two-Material Elastic Bar}\label{hBar}

\begin{figure*}[h]\centering
	\includegraphics[width=0.98\columnwidth]{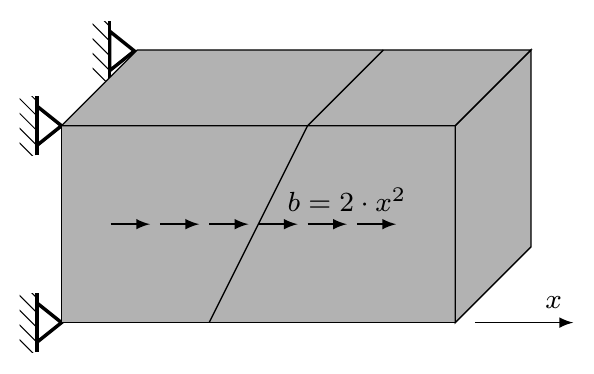}
	\includegraphics[width=0.98\columnwidth]{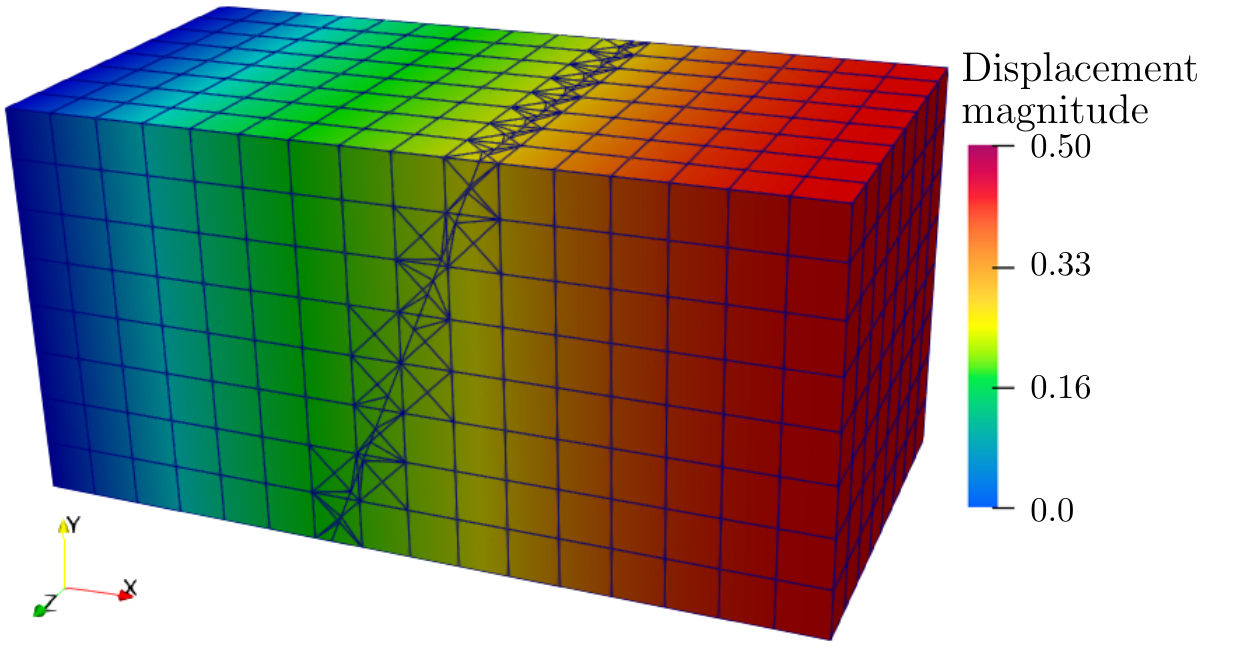}
	\caption{Three dimensional bar problem: problem setup (left); contours of displacement magnitude evaluated on coarsest discretization (right)}
	\label{fig_Sketch_Bar}
\end{figure*}

This example studies the convergence rates of the proposed immersed B-spline discretization framework. We consider the 3D bar shown in Fig.~\ref{fig_Sketch_Bar}. The bar has the dimension $1.0 \times 0.5 \times 0.5$. The left face of the bar is clamped, and the bar experiences a body load in axial direction $b_x=2 x^2$. The bar is composed of two linear elastic, isotropic materials, separated by an interface which is inclined against the x-axis. To facilitate comparison against an analytical solution, both materials are assigned the same properties: a Young's modulus $E = 1.0$ and a Poisson's ratio $\nu = 0.0$. Note that this setup allows for a 1D analytical model. The analytical displacement in x-direction is: 
\begin{equation}
u(x)=\frac{L}{6EA} \left(4 L^3 x - x^4 \right),
\end{equation}
where $L=1.0$ is the bar length and $A=0.25$ the cross-sectional area. 

To study the influence of the intersection configuration on the finite element solution with B-spline background mesh refinement, we rotate the interface in ten steps between $\frac{\pi}{4}$ and $\frac{\pi}{2}$ degrees. For each orientation, we increase the number of background mesh elements from $8 \times 4 \times 4$ to $ 128 \times 64 \times 64$ B-spline elements through uniform mesh refinement. We repeat this study for linear, quadratic, and cubic B-spline discretizations.

\begin{figure}[h]\centering
	\includegraphics[width=0.98\columnwidth]{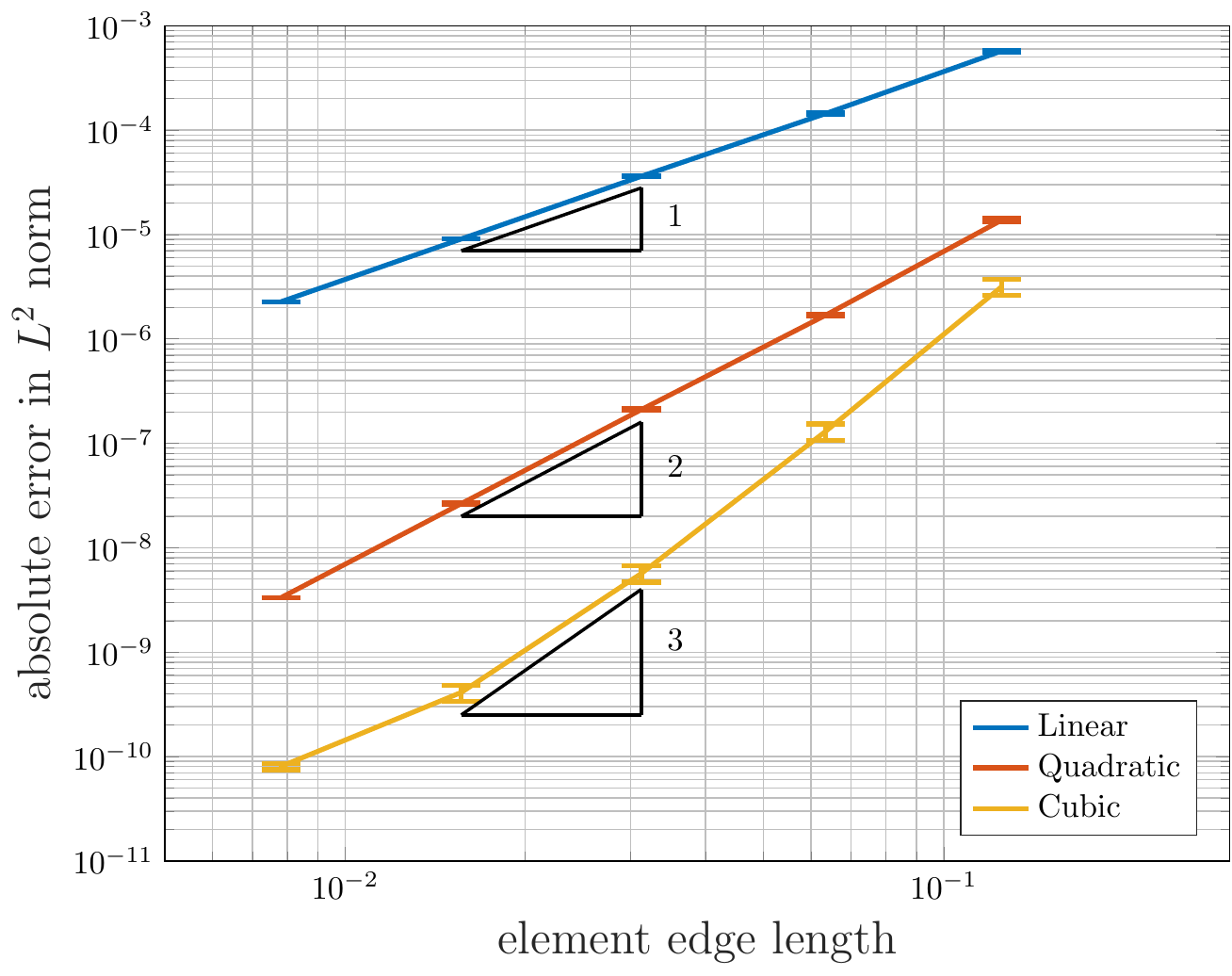}
	\includegraphics[width=0.98\columnwidth]{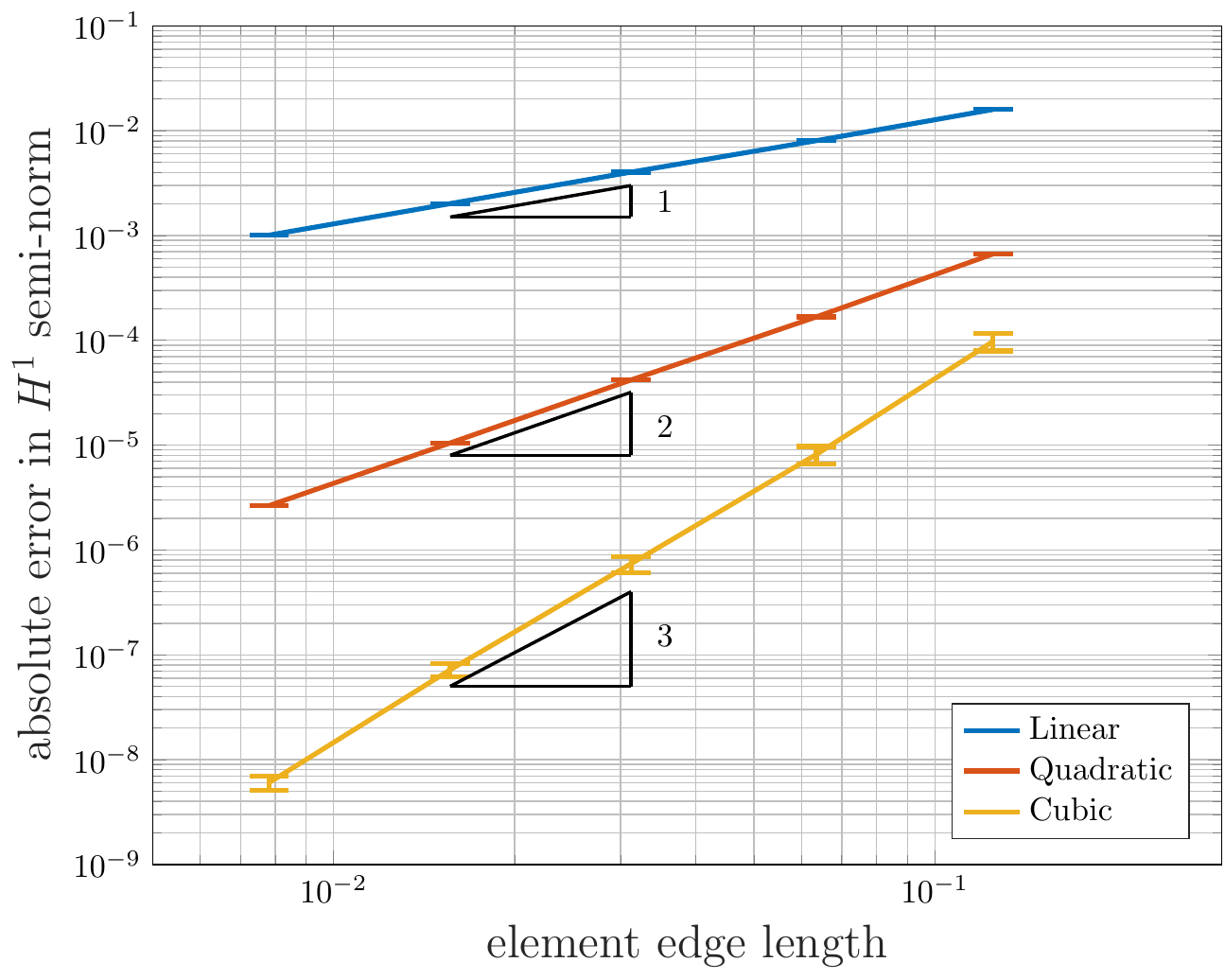}
	\caption{Convergence rates in the $L^2$ error norm and $H^1$ semi-norm for the presented two-material bar problem}
	\label{fig_bar_convergence}
\end{figure}

In Fig.~\ref{fig_bar_convergence}, we plot the mean error for all intersection configurations and for each approximation order over the B-spline element edge length as described in Eqs.(\ref{Eq::L2Norm}) and (\ref{Eq::H1SemiNorm}). The error bars represent the standard deviation for all intersection configurations per order and mesh refinement. For each B-spline order, we visualize the convergence rate by the triangles inserted in Fig.~\ref{fig_bar_convergence}. The numerical results show that the convergence rates of the $L^2$ error norm agree with the theoretical, optimal convergence rates of $p+1$, where $p$ is the polynomial B-spline order. Similarly, the theoretical, optimal convergence rate of $p$, as determined by \cite{Evans2009,Burman2015}, is achieved in the $H^1$ semi-norm for linear, quadratic, and cubic discretizations. Our results suggest that the proposed immersed finite element approach recovers optimal convergence rates for sufficiently smooth state variable fields with uniform mesh refinement.

\subsection{Single-Material Thermal Diffusion Problem}\label{L-Domain}

Analysis problems often include regions where the state variable fields exhibit large spatial gradients. To reduce the local and global approximation errors, a fine discretization is needed in these regions. The following thermal diffusion example demonstrates the ability of the proposed immersed finite element framework to construct and locally refine discretizations.

\begin{figure}[h]\centering
	\includegraphics[width=0.9\columnwidth]{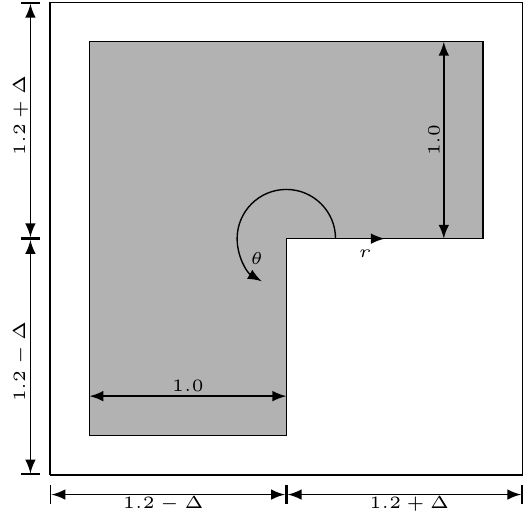}
	\includegraphics[width=0.9\columnwidth]{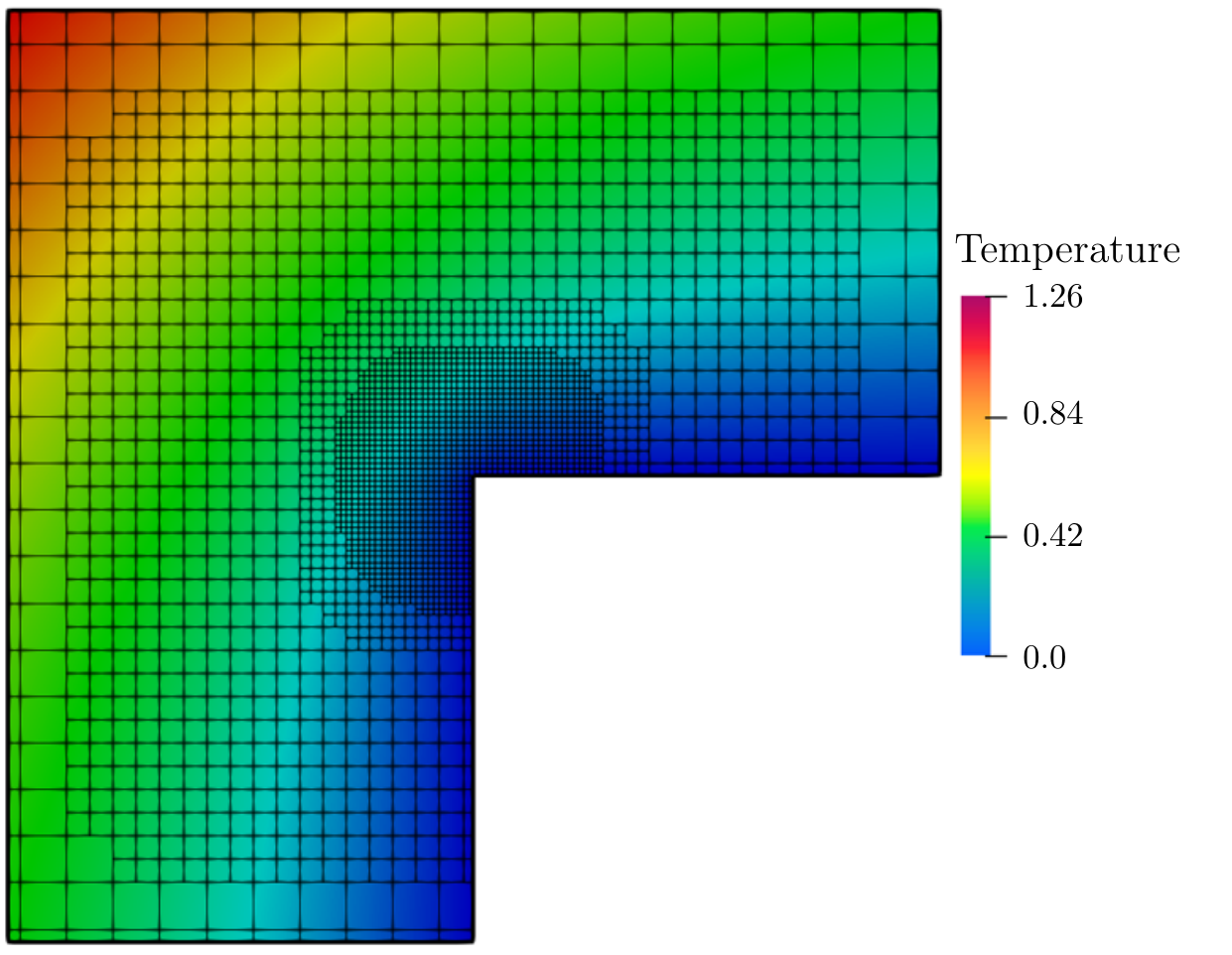}
	\caption{Single-material thermal diffusion problem: L-shaped domain immersed into a rectangular computational domain (top); locally  refined B-Spline discretization around the reentrant corner and temperature field (bottom) }
	\label{fig_Scetch_Example_4}
\end{figure}

We solve a thermal diffusion problem as described in Eq.~\eqref{Eq_TLin} in an L-shaped domain, assuming a single-material with an isotropic conductivity $\kappa = 1$, see Fig~\ref{fig_Scetch_Example_4}. For this problem, the exact solution in polar coordinates is given as (\cite{DAngella2016}):
\begin{equation}
T(r, \theta) = r^{ \frac{2}{3}} \mathrm{sin}( \frac{2}{3} \theta )
\end{equation}
The exact solution is enforced weakly to the entire boundary of the physical domain. The origin of the coordinate system is located at the reentrant corner. Note that the spatial gradients increase for $r \rightarrow 0$ and are infinite at $r = 0$.

Using the proposed analysis framework, the physical domain is immersed into a rectangular computational domain as shown in Fig.~\ref{fig_Scetch_Example_4}. We compare the convergence rates for uniform \textit{h}-refinement and for local \textit{h}-refinement around the reentrant corner.

The coarsest THB background mesh has $6 \times 6$ elements. Sequences of THB discretizations are constructed by uniform or local refinement with up to 6 refinement levels, see Section~\ref{sec:details}. The offset $\Delta$ between the computational and physical domains is chosen as $\Delta = - 0.127$. This offset guarantees that the L-shaped domain boundary does not align with the THB background mesh for any mesh refinement configuration. 

\begin{figure}[t!]\centering
	\includegraphics[width=0.96\columnwidth]{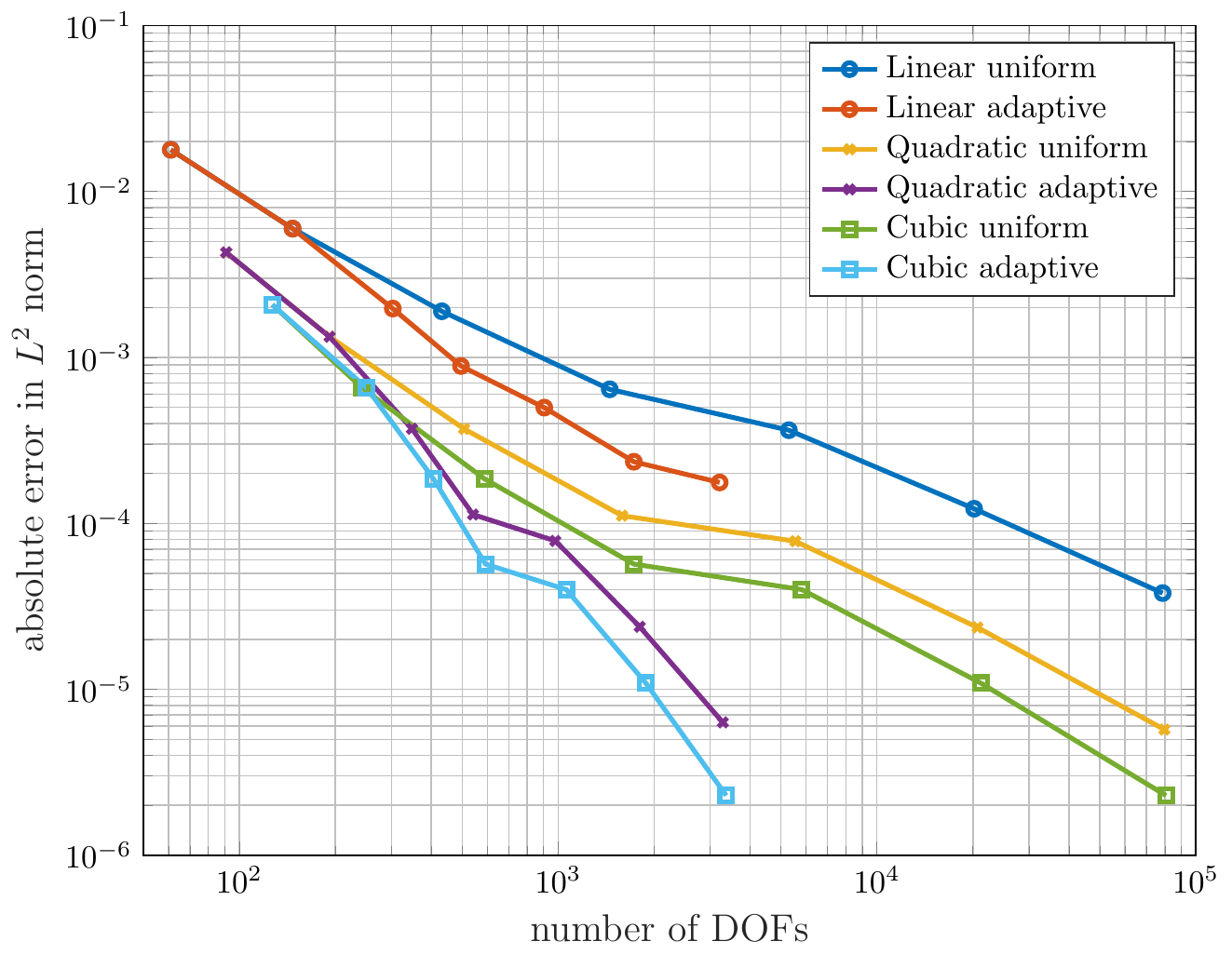}
	\caption{Error convergence rate of temperature field with mesh refinement in the $L^2$ norm}
	\label{fig_Example4_conv}
\end{figure}

For different B-spline order and mesh refinements, we plot the error in the $L^2$ norm versus the number of Degrees Of Freedom (DOFs) in Fig.~\ref{fig_Example4_conv}. Less DOFs are required to meet a specific error requirement for local refinement when compared to uniform refinement. This suggests that  local refinement may lead to a significant reduction in computational cost. Due to the singularity at the reentrant corner, optimal convergence rates with mesh refinement cannot be recovered. Moreover, this example demonstrates that the proposed framework allows for combining local \textit{h}-refinement with higher order B-spline discretizations. For example, using cubic B-spline and local refinement leads to the lowest error for any number of DOFs. We do not present $H^1$ semi-norm error plots as they are not meaningful due to the singularity at the reentrant corner.

\subsection{Polycrystalline Micro-structure}\label{grainStructure}

The proposed analysis framework is suited to model complex multi-material problems. While the previous examples established accuracy for single-material problems, this example demonstrates the multi-material capabilities of our approach and studies the computational cost associated with the generation of THB and union background meshes.

We consider a representative volume element of a polycrystalline micro-structure and analyze its structural response. The edge length of the volume element is $150 \, \mu m$. The grain geometries as well as the grain material are defined through a 3D image file with $149 \times 149 \times 149$ voxels; see, for example, \cite{Rodgers2015}. The voxels define 471 individual grains. To demonstrate the multi-material capabilities, unique material properties are assigned to each grain. The Young's modulus $E$ varies in the range of $[50 \cdot 10^3 \   \ldots \ 500 \cdot 10^3]\mathrm{MPa}$ and the Poisson ratio $\nu$ in a range of $[0.25 \ldots 0.35]$.

\begin{figure}[t!]\centering
	\includegraphics[width=0.96\columnwidth]{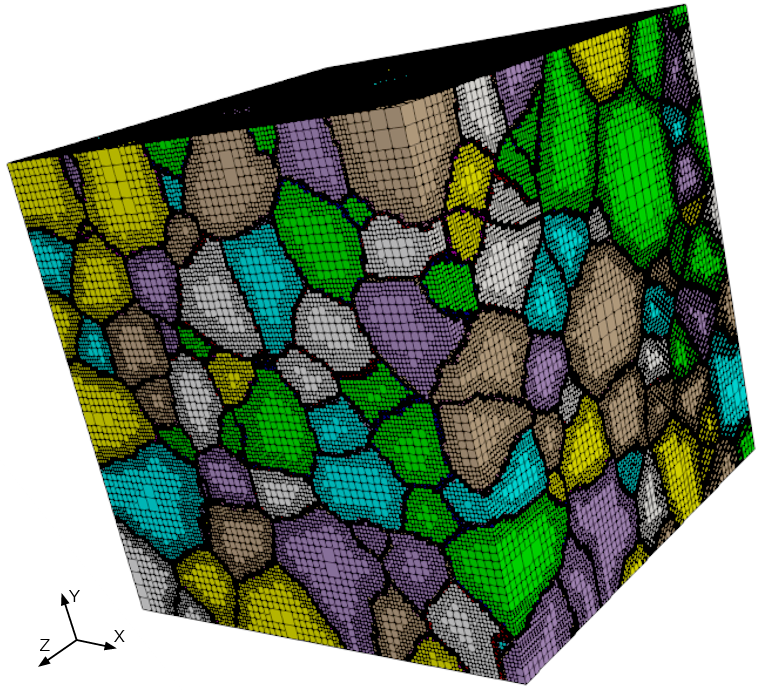}
	\caption{Locally refined XFEM analysis model of a polycrystalline micro-structure }
	\label{fig_Sketch_Example3_1}
\end{figure}

From the 3D image, we generate a level-set field for each phase such that the grain geometries are represented by the zero isocontours of the level-set fields. The grain geometries are immersed into a cubic domain, as illustrated in Fig.~\ref{fig_Sketch_Example3_1}. Note that the voxel-based grain shapes could be preprocessed to obtain smoother grain interfaces. However, this option is not utilized here to demonstrate the ability of the proposed analysis framework to operate directly on complex non-smooth voxel-based geometries. 

\begin{figure*}[h]\centering
	\includegraphics[width=0.98\columnwidth]{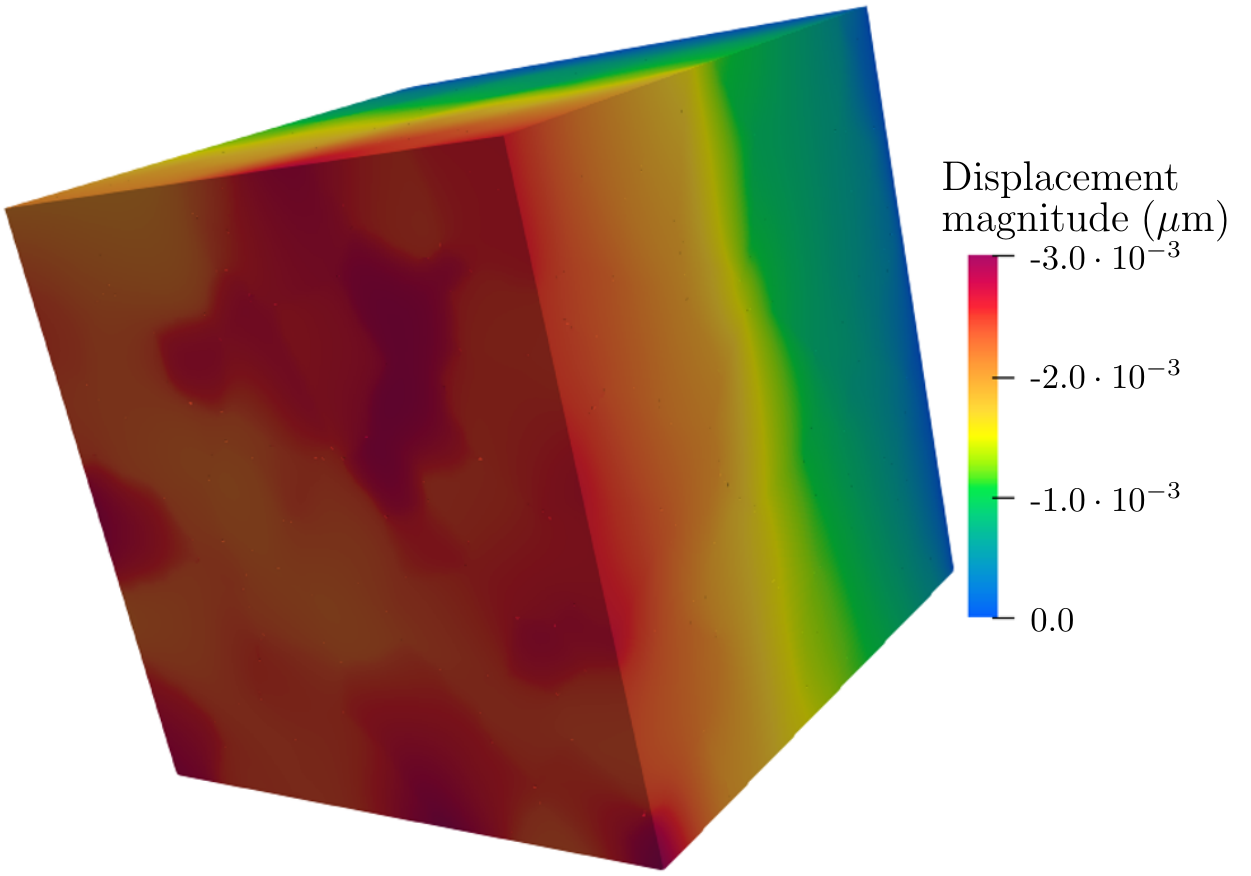}
	\includegraphics[width=0.98\columnwidth]{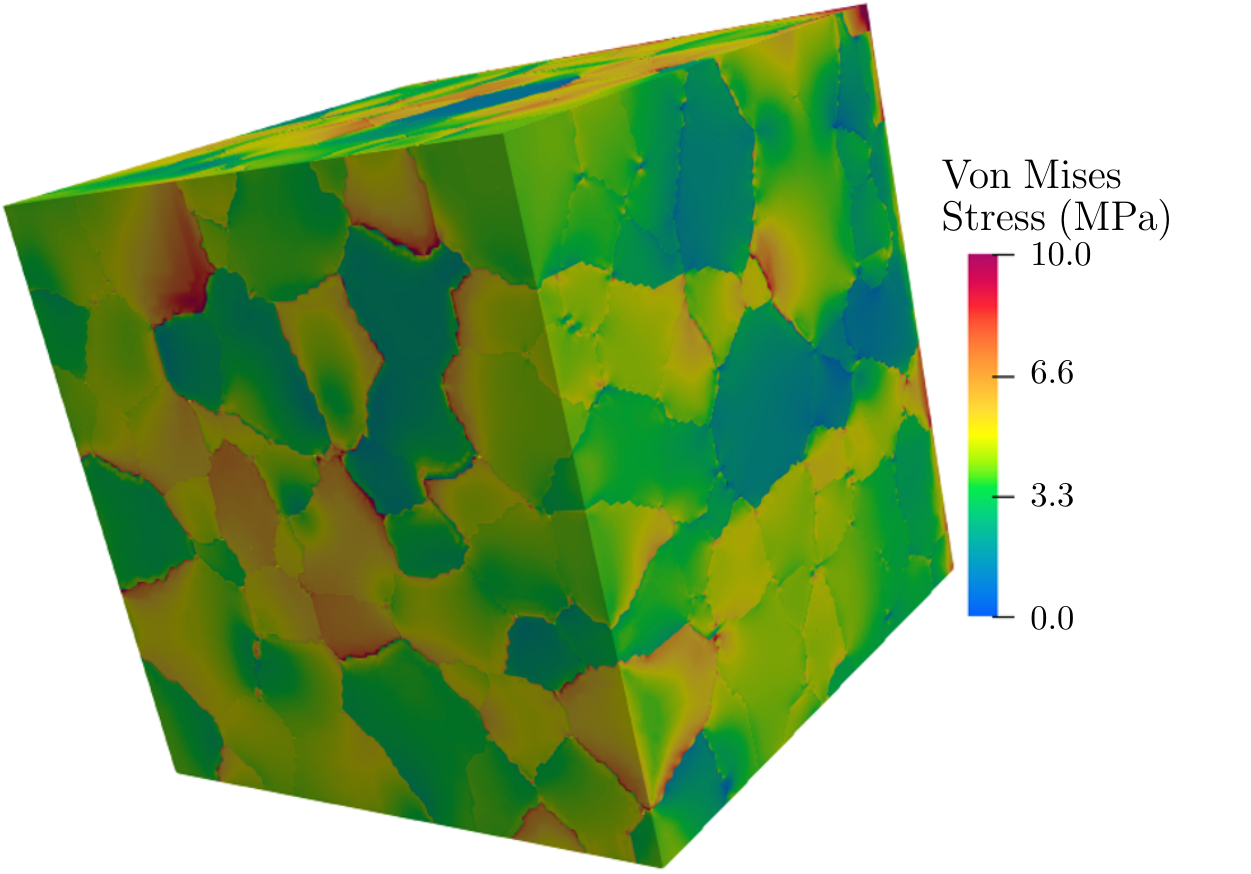}
	\caption{XFEM analysis of polycrystalline micro-structure: Displacement magnitude contours (left); Von Mises stress contours (right)}
	\label{fig_Sketch_Example3_2}
\end{figure*}

The mechanical behavior of the polycrystalline micro-structure is modeled by linear elasticity with an isotropic material behavior. Perfect bonding is assumed at the grain interfaces. To demonstrate the ability of performing an XFEM analysis for this micro-structure, we apply a pressure load of $F= \text{-} 4.0 \, \mathrm{MPa}$ on the $z = 0$ face and a zero-displacement Dirichlet boundary condition on the opposite face. The initial uniform B-spline background mesh has $32 \times 32 \times 32$ elements. In addition, two local refinement steps are performed at all grain boundaries. The displacement magnitude and Von Mises stress contours are shown in Fig.~ \ref{fig_Sketch_Example3_2}.

\subsubsection{Scalability study}\label{ScalibilityStudy}

The PT data structure and domain decomposition strategy presented in Subsection \ref{Paralle_Considerations} enable an efficient, parallel implementation of the proposed discretization framework. With the polycrystalline example, we study the overall performance and parallel scalability of this implementation.  

For the scalability study, we start with $200 \times 200 \times 200$ PT cells at refinement level 0 and refine locally twice around the grain boundaries. For the demonstration of the scalability, we purposefully exploit a finer mesh than as in the previous analysis to guarantee a large ratio of interior to aura PT cells on each processor local subdomain for a large processor count. A small ratio of interior to aura elements negatively impacts the scalability due to increased inter-processor communication.

The resulting THB and union background meshes have each a total of $105,558,692$ elements. The creation of the PT data structure and the derived discretizations can be subdivided into five distinct steps, see also Subsection \ref{sec:details}. These steps include the refinement of the PT cells, the construction of the THB and union background meshes, the construction of the extraction operators, and the construction of facets needed for ghost stabilization. 

\begin{figure}[h]\centering
	\includegraphics[width=0.96\columnwidth]{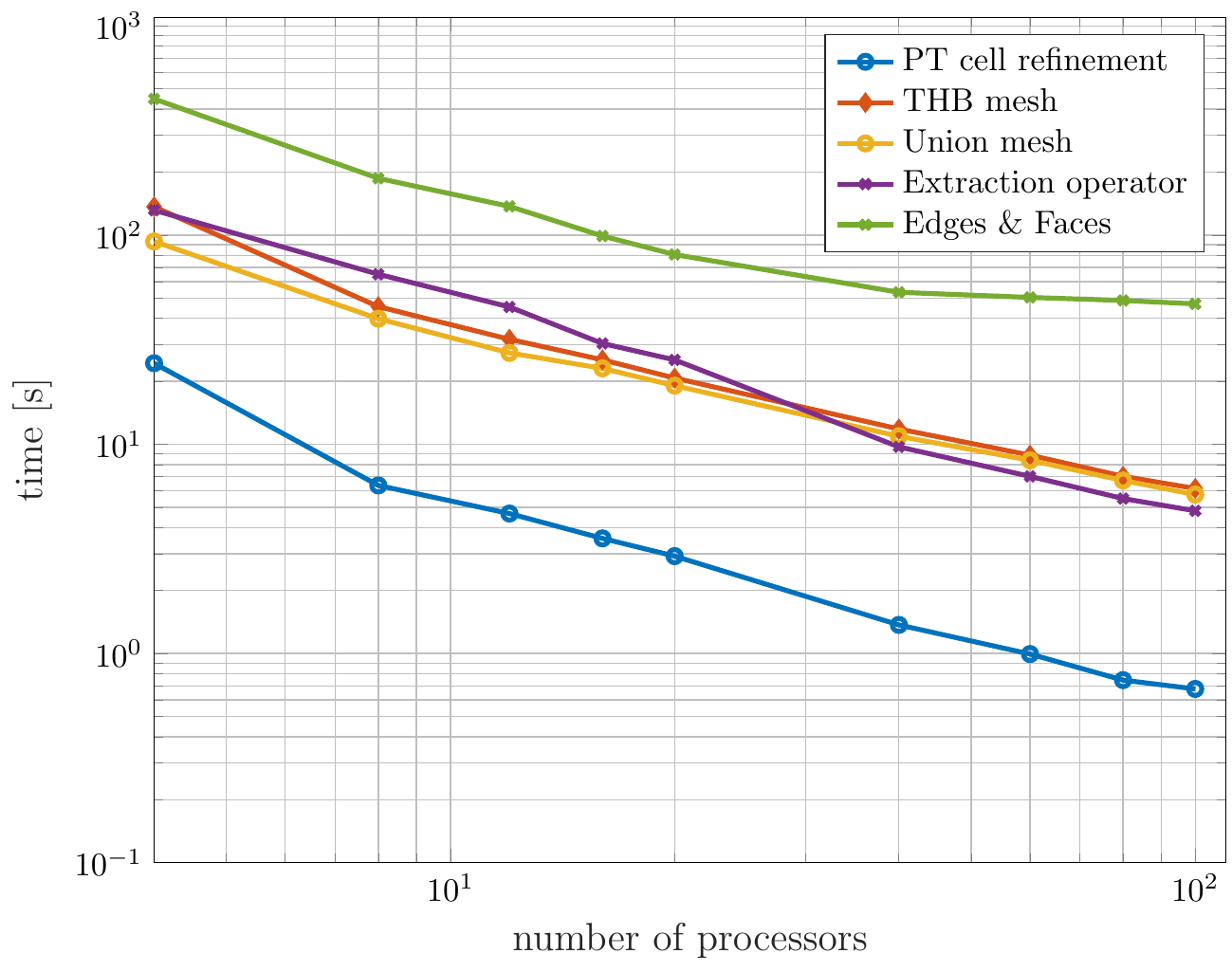}
	\caption{Scalability study with processor count of the hierarchical meshing library, using $200 \times 200 \times 200$ PT cells at refinement level 0 and two refinement steps at grain boundaries}
	\label{fig_Example_3}
\end{figure}

To characterize the scalability of our implementation, we generate the THB and union background meshes in parallel, varying the number of subdomains from $4$ to $160$. Fig.~\ref{fig_Example_3} shows the execution time for the mesh generation only, i.e., the time needed for XFEM analyses is omitted. The computations are preformed on four Intel Xeon Platinum 8160 "Skylake" nodes with 24 cores each, distributing the subdomains equally across all cores. 

Mostly linear scaling with the processor count is observed for the refinement operation and the construction of the projection operators. Such behavior is expected as these steps do not need any inter-process communication. The construction of THB and union background meshes, as well as the construction of facets, show a linear scaling for a low processor count. They start to plateau with a high processor count, as these steps require communication of IDs in the aura. The aura size is based on refinement level $l^0$ and the buffer size as outlined in Section ~\ref{Paralle_Considerations}. Increasing the processor count while keeping the total domain size constant increases the ratio of aura PT cells to interior cells. This consequently increases communication and affects scalability. 

The scalability study demonstrates that a mesh with over $10^8$ elements can be generated on just four nodes with a total of 96 processors in less than 50 seconds. Furthermore, when creating a higher order discretizations, the computational time for the refinement of the background meshes and the construction of faces stays the same as neither are affected by the interpolation order. The construction of the discretizations, as well as the calculation of the extraction operators, needs slightly more memory and computational time as the support of higher order basis functions is increased.

\subsection{Thermo-elastic Plate with Elliptic Hole}\label{elllipticHole}

\begin{figure*}[h]\center
	\includegraphics[width=0.68\columnwidth]{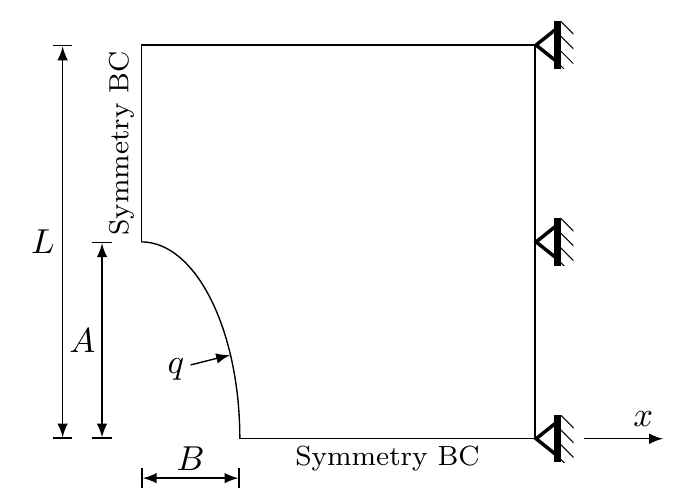}
	\includegraphics[width=0.62\columnwidth]{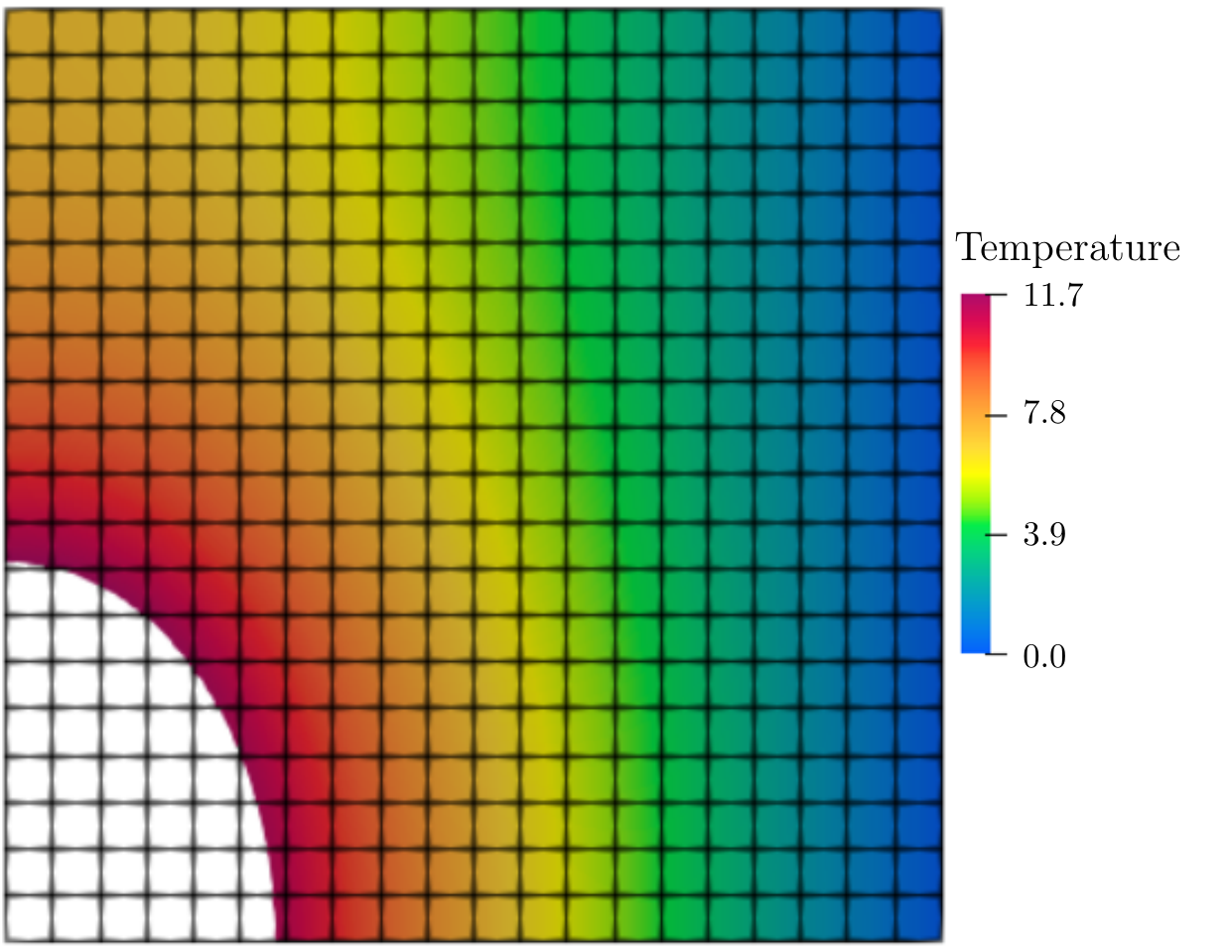}
	\includegraphics[width=0.62\columnwidth]{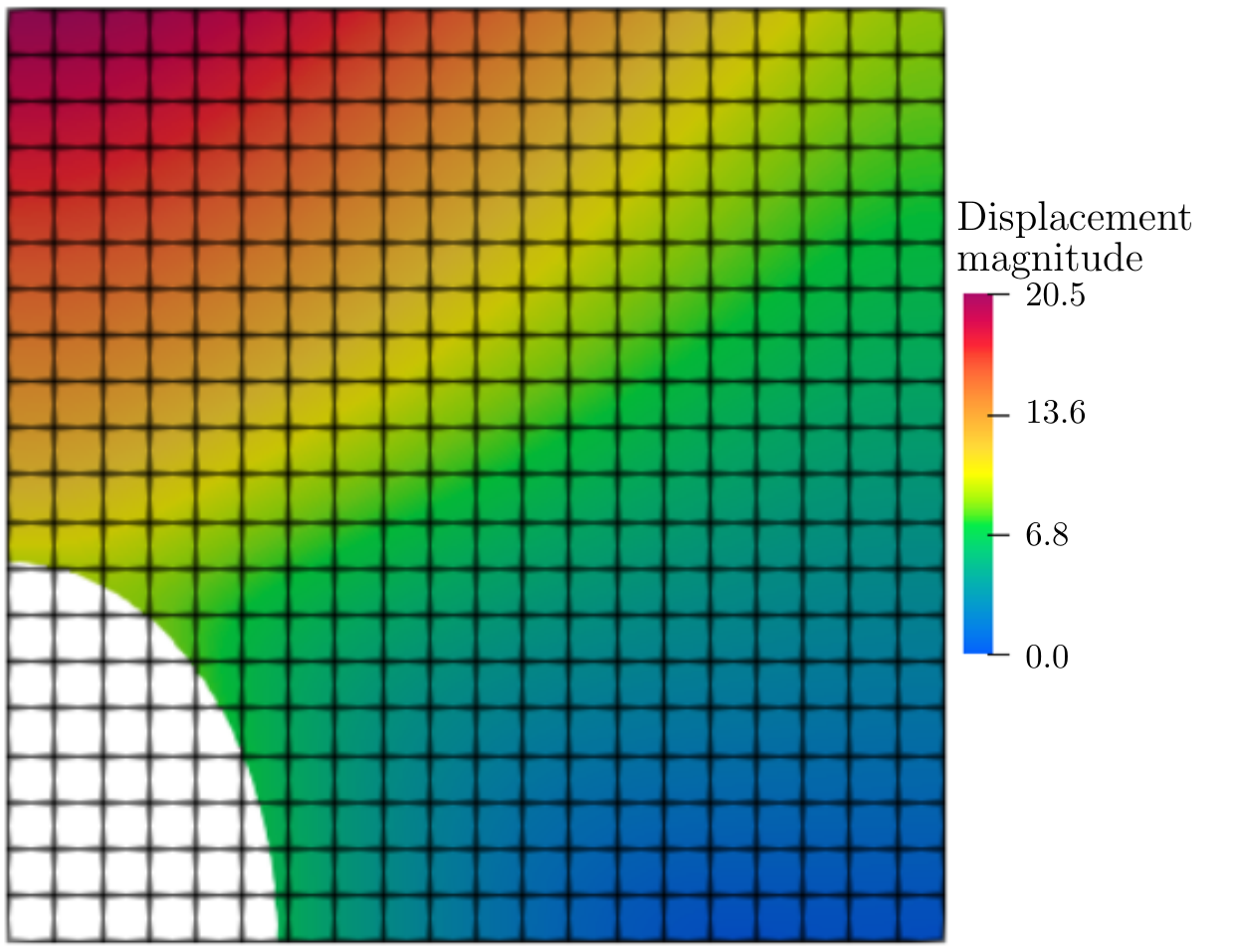}
	\caption{Thermo-elastic plate with elliptic hole under thermal load: problem setup (left); temperature contours (center); displacement magnitude contours (right)}
	\label{fig_Sketch_Example2_1}
\end{figure*}

This example considers a single-material but multi-field configuration. We study the problem of a thermo-elastic plate with an elliptical inclusion under in-plane tension due to a thermal load as illustrated in Fig.~\ref{fig_Sketch_Example2_1}. By exploiting symmetry, we only model a quarter of the domain. The setup consists of a solid, two-dimensional, square domain with a length of $L = 2.0$ and an elliptical inclusion at the origin with a semi-major axis of $A=0.8136$ and a semi-minor axis of $B=0.5753$. These values for the semi-axes are chosen such that the immersed interface does not align with the THB background mesh for all refinement levels. The temperature field is due to a uniform heat flux $q =10.0$ applied along the elliptical hole and a prescribed temperature $T = 1.0$ at the right domain boundary. The matrial conductivity is $\kappa=1.0$. The structural response is described by a linear thermo-elastic model and an isotropic constitutive behavior, with a Young's modulus $E = 1.0$, a Poisson's ratio $\nu = 0.3$, and a coefficient of thermal expansion $\alpha = 1.0$. The reference temperature is set to $T_0 = 0.0$. 

To study the convergence rate of the XFEM solution with mesh refinement, we evaluate the thermal and structural response on a series of B-spline discretizations. Starting from the coarsest THB background mesh with $10 \times 10$ elements, uniformly refined B-spline discretizations with $20 \times 20$, $40 \times 40$, $80 \times 80$, $160 \times 160$ and $320 \times 320$ elements are considered. The geometry of the elliptical inclusion is defined by an analytical level set function. A reference solution is computed numerically by calculating both fields on a uniformly refined B-spline mesh that is 7 times refined, $1280 \times 1280$  compared to the coarsest background mesh. Numerical studies lead by the authors have shown that the computation of the reference solution on a finer mesh does not noticeably improve the quality of the error computation. 

For all B-spline discretizations, the union background mesh is 7 times uniformly refined relative to the coarsest background mesh. This eliminates the geometrical error between the coarse and the reference solution. Moreover, it simplifies the integration of the $L^2$ and $H^1$ errors as both the coarse and the reference solutions are represented on equally refined XFEM background meshes. 

\begin{figure*}[h]\centering
	\includegraphics[width=1.8\columnwidth]{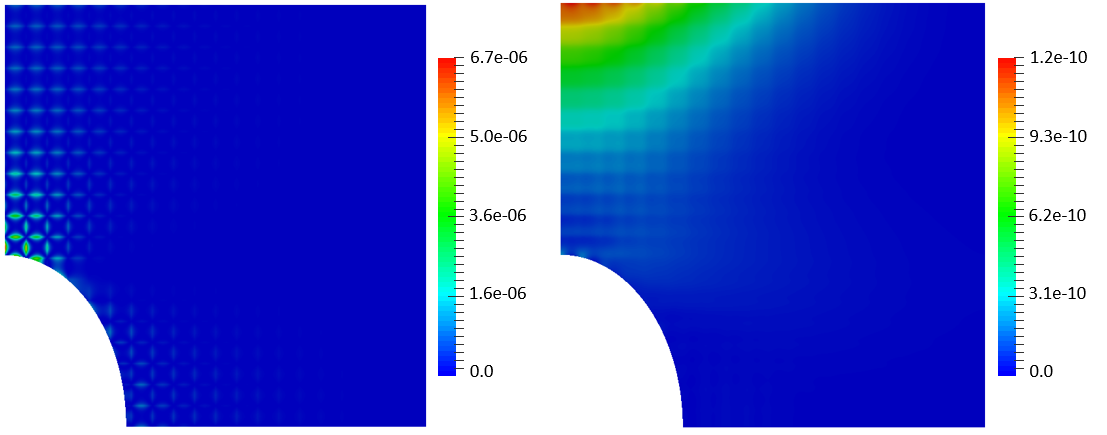}\hspace{0.25cm}
	\caption{Visualization of absolute error of temperature (left) and displacement field (right)}
	\label{fig_Example_2_Error}
\end{figure*}

To illustrate the interdependence of discretization errors in multi-physics problems, we vary the refinement for the temperature discretization while the displacement field is always evaluated on the finest B-spline discretization, i.e., $320 \times 320$ elements. The contours of the absolute errors of the temperature and displacement fields are visualized in Fig.~\ref{fig_Example_2_Error} for the case where the temperature field is discretized by linear B-splines on a $20 \times 20$ element mesh. This coarse temperature field discretization leads to errors in the interior of the elements. Therefore, the error visualization presents the characteristic tiled error patterns. 

\begin{figure*}[h]\centering
	\includegraphics[width=0.96\columnwidth]{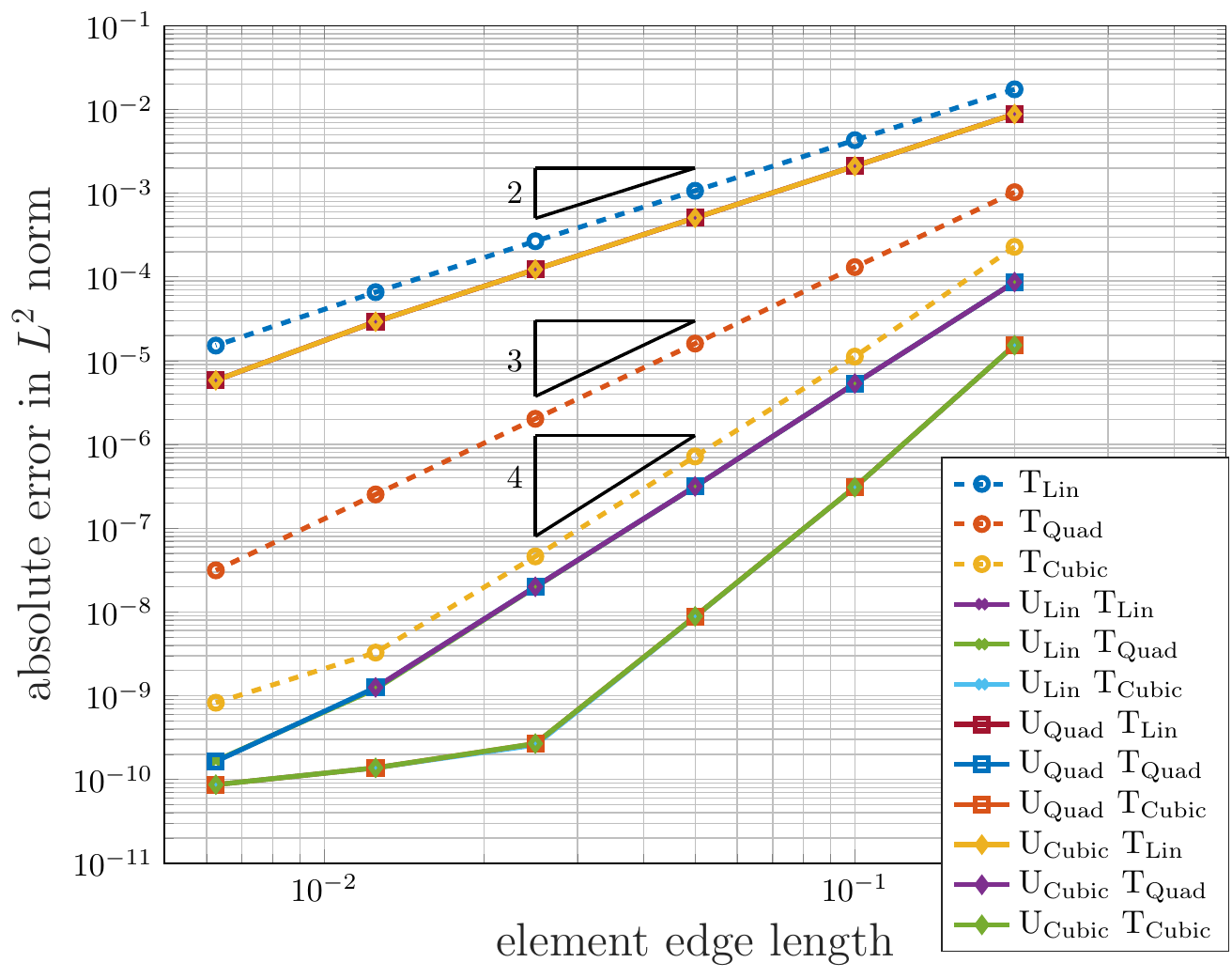}\hspace{0.25cm}
	\includegraphics[width=0.96\columnwidth]{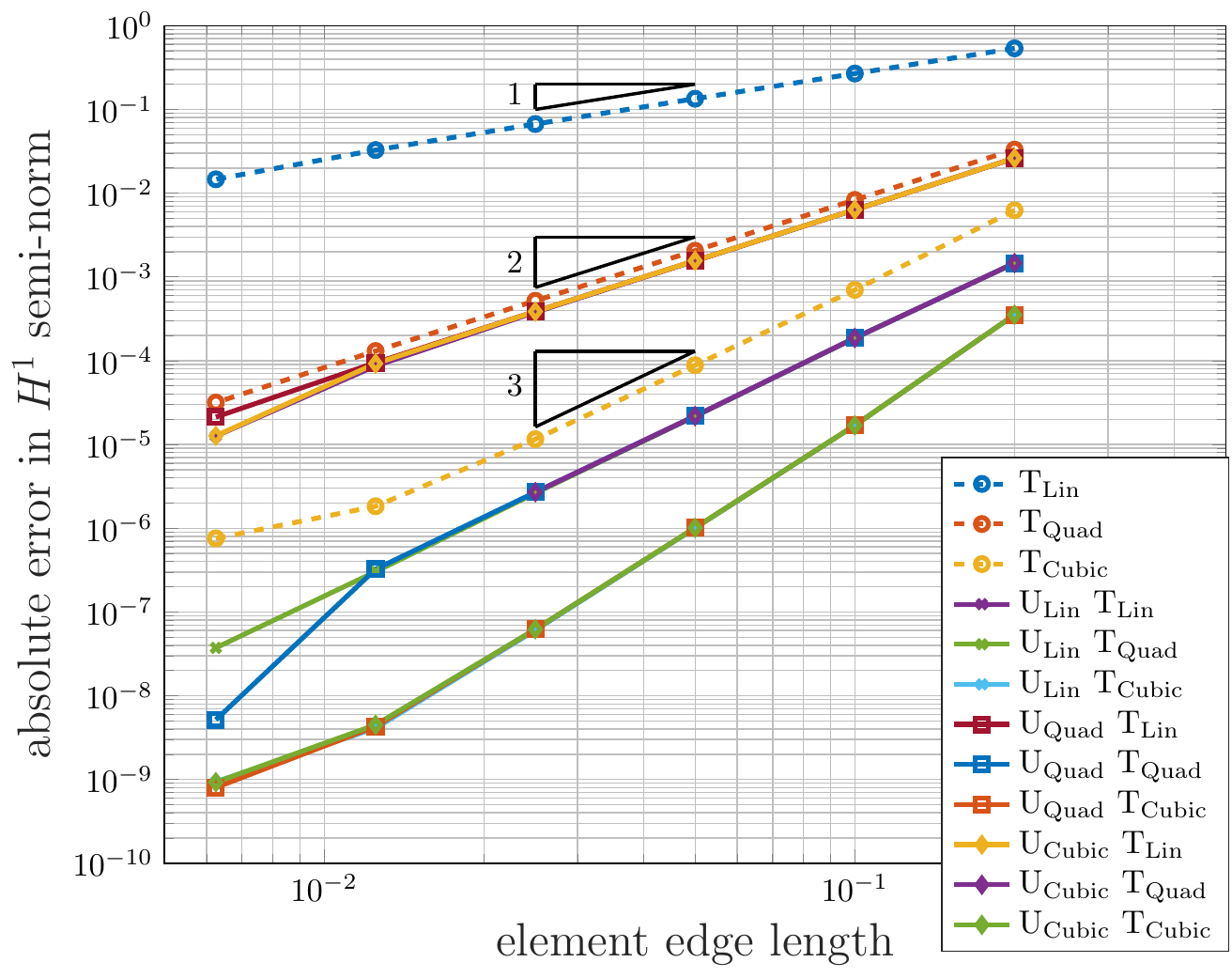}
	\caption{Error convergence rate of temperature and displacement fields with mesh refinement in the $L^2$ norm (left) and $H^1$ semi-norm (right) }
	\label{fig_Example_2}
\end{figure*}

The $L^2$ and $H^1$ semi-norm errors of the temperature and the displacement fields are presented in Fig.~\ref{fig_Example_2}. First, we consider the error of the temperature field. The convergence of the $L^2$ error norm for linear, quadratic, and cubic B-spline discretizations shows optimal convergence rate of $p+1$. For the $H^1$ semi-norm, we observe optimal convergence rates of $p$. As expected, using a higher-order basis functions results in smaller errors for the temperature field. 

Next, we consider the error of the displacement field, which is evaluated with linear, quadratic, and cubic basis functions and depends on the temperature field that is discretized by either linear, quadratic, and cubic B-splines. The errors of the displacement field in the $L^2$ norm and $H^1$ semi-norm are presented in Fig.~\ref{fig_Example_2}. The convergence rate of the error with mesh refinement in the $L^2$ norm is independent of the order of the displacement field basis functions and instead only depends on the order of the temperature field. Since the displacement field is discretized on a much finer mesh than the thermal field, the error of the displacement field is dominated by the error of the thermal field. The same observation can be made for the convergence rate in the $H^1$ semi-norm.

\begin{figure*}[h]\centering
	\includegraphics[width=0.96\columnwidth]{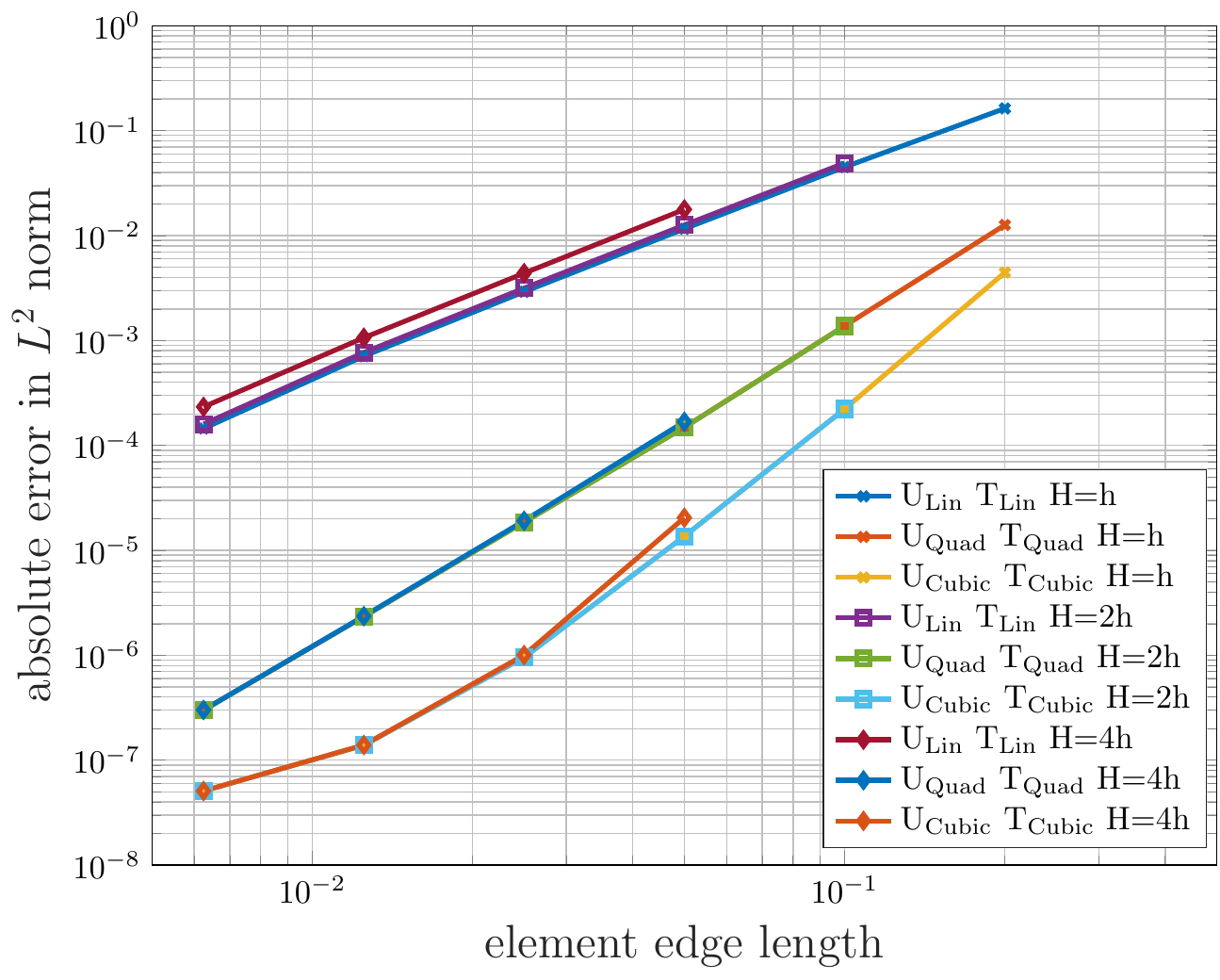}\hspace{0.25cm}
	\includegraphics[width=0.96\columnwidth]{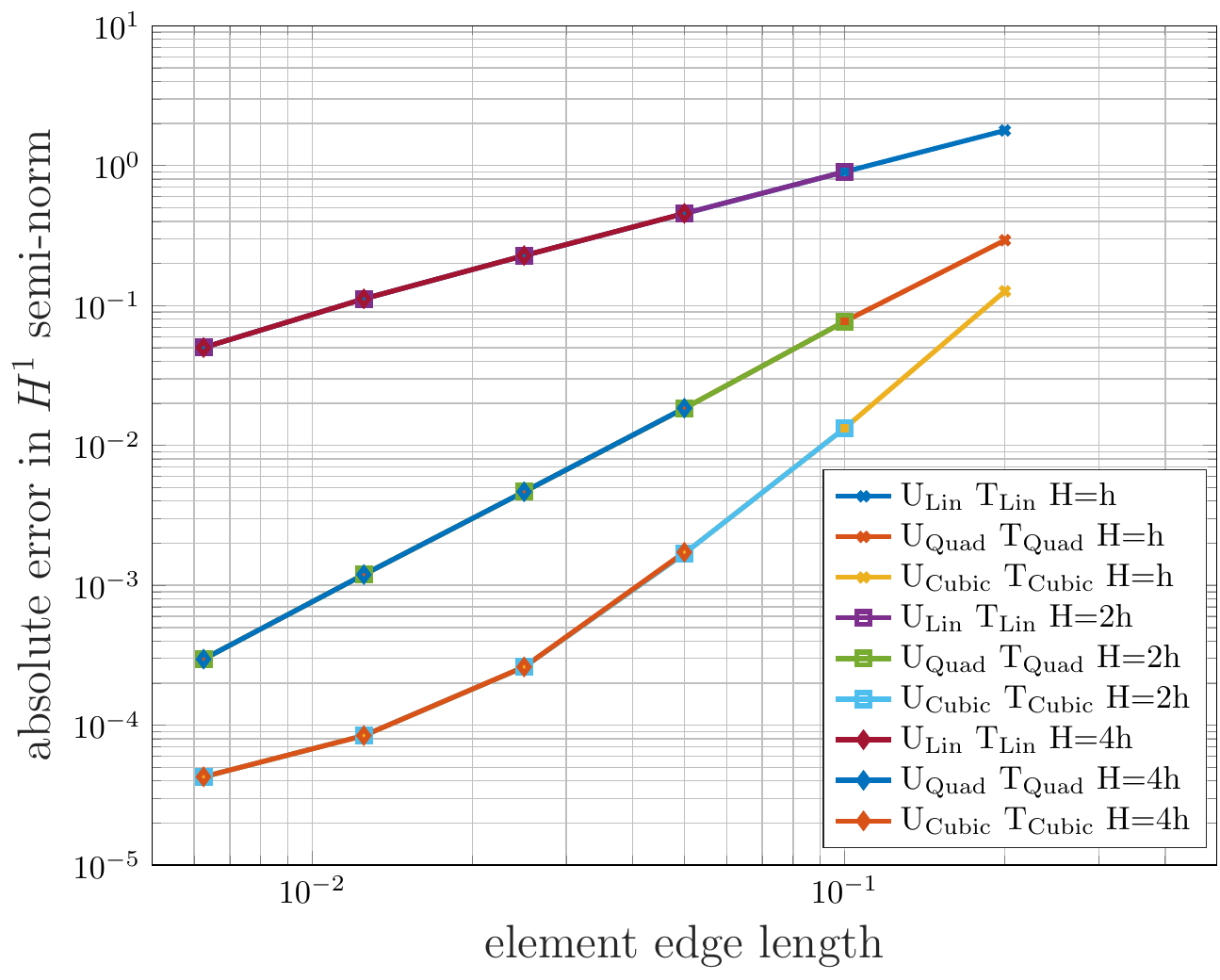}
	\caption{Error convergence rate the displacement field with mesh refinement. $L^2$ norm (left) and $H^1$ semi-norm (right). Temperature and displacement field are refined simultaneously }
	\label{fig_Example_2_hRef}
\end{figure*}

To gain further insight, we refine simultaneously the thermal and displacement fields with a specific difference in refinement level between both fields. The THB background element edge lengths of the thermal and displacement fields are denoted by $H$ and $h$, respectively. The study is performed for a maximal difference in element size of $H= 4h$. This study is performed for linear, quadratic, and cubic interpolation orders for the thermal and the displacement fields. The setup of the thermal problem is identical the one presented above, see Fig.~\ref{fig_Example_2}. The convergence rate of the displacement field  with mesh refinement in the $L^2$ and $H^1$ semi-norm is presented in Fig.~\ref{fig_Example_2_hRef}. We observe that the absolute error and convergence rate of the displacement problem only shows minimal differences when choosing an up to two times coarser thermal field than the displacement field. 

\begin{figure*}[h]\centering
	 \includegraphics[width=0.96\columnwidth]{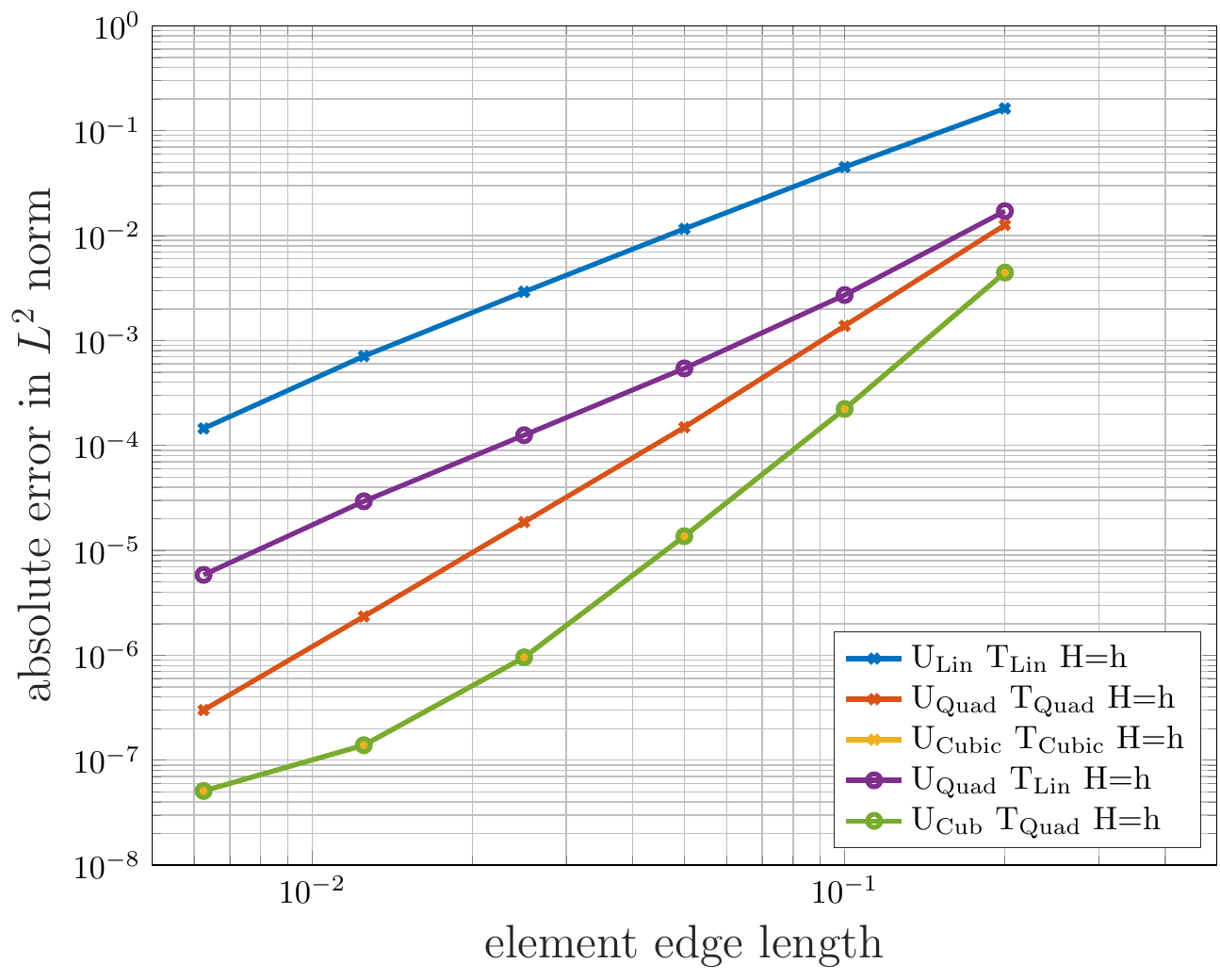}\hspace{0.25cm}
	\includegraphics[width=0.96\columnwidth]{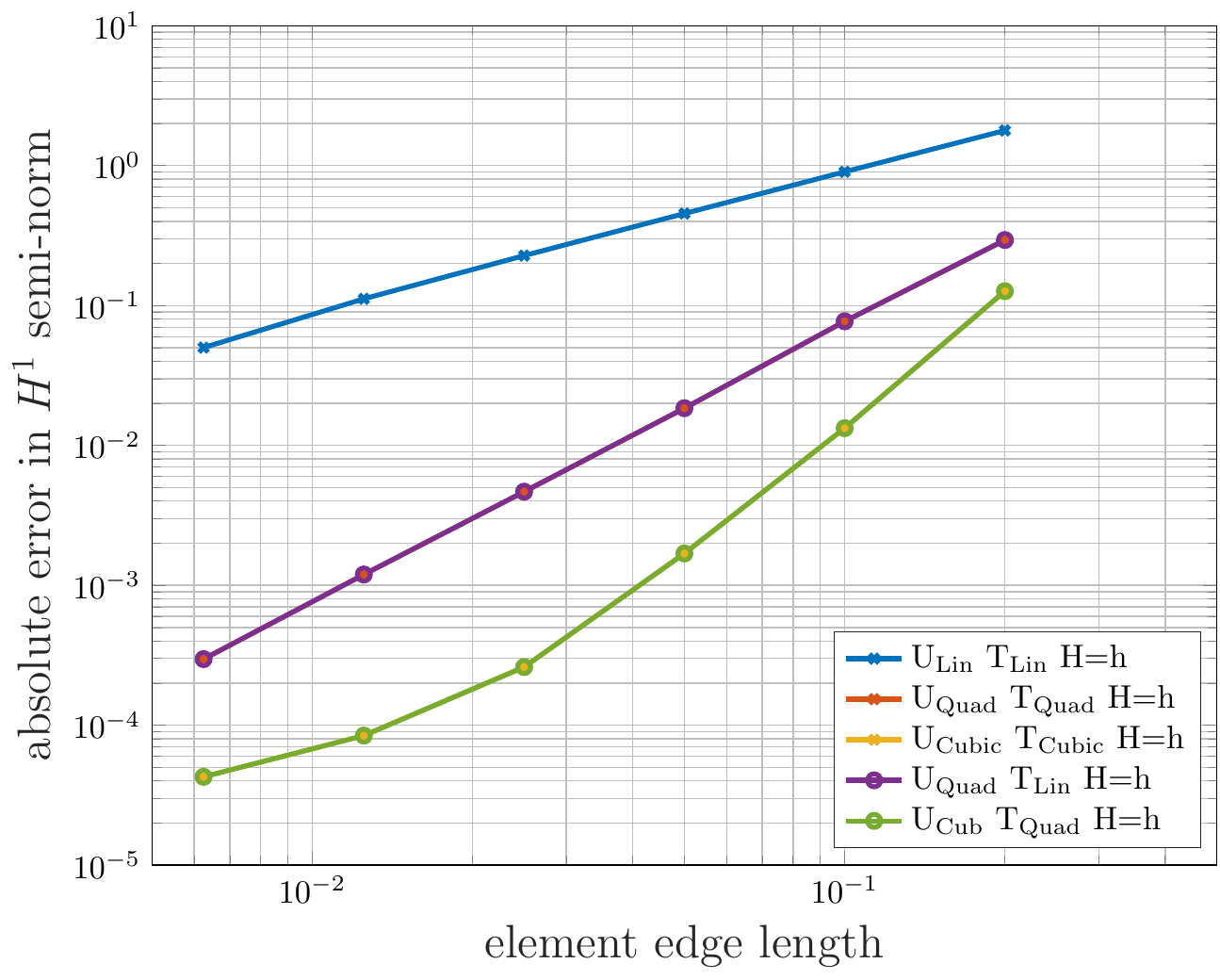}
	\caption{Error convergence rate the displacement field with mesh refinement in the $L^2$ norm (left) and $H^1$ semi-norm (right) using identical refinements for the temperature and displacement fields}
	\label{fig_Example_2_pRef}
\end{figure*}

Lastly, we present a convergence study where we allow for a difference in polynomial order. The thermal and displacement field are refined simultaneously such that both fields are on the same refinement level, i.e., $H=h$. Results are presented for a polynomial order of the displacement field $p_U$ for $p_U = p_T$ and $p_U = p_T + 1$, where  $p_T$ is the order of the thermal field. The convergence rates with mesh refinement of the displacement field in the $L^2$ and $H^1$ semi-norm are presented in Fig.~\ref{fig_Example_2_pRef}. We observe that choosing the thermal field one polynomial order lower than the displacement field results in the same convergence rate of the displacement field in the $H^1$ semi-norm. For this particular problem, theoretical results are not available for the rate of the convergence for the displacement field in the $L^2$ norm, as the well-known Aubin-Nitsche technique cannot be applied in the thermo-elastic setting.


This example illustrates that for multi-physics problems, the discretization of individual fields needs to be performed in concert. The proposed discretization framework provides an efficient and flexible tool to select the interpolation order and mesh refinement level for individual fields to obtain a numerical solution that meets accuracy requirements with minimal computational costs, i.e., with minimal number of DOFs.

\subsection{Thermo-elastic Multi-material Problem}\label{MultiMaterialProblem}

This final example studies a multi-material, multi-physics configuration. For such problems, areas with large spatial gradients can vary based on the type of physics and material. Differently locally  refined discretizations for each field reduce the computational cost while simultaneously enabling for an accurate evaluation of the physical responses. To demonstrate this aspect, we consider the two-material, thermo-elastic problem of an expanding circle embedded in a non-expanding plate, as presented in the introduction in Fig.~\ref{subfig-2:1}. The circular inclusion is occupied by a material A and the plate by a material B. The temperature field is due to a spatially varying heat load of $q= 100.0 \cdot \sin (10.0 \cdot y)+110.0$ applied along the left domain boundary and a prescribed temperature $T = 0.0$ at the right domain boundary. The conductivities of materials A and B are identical and are $\kappa_{A}=\kappa_{B}=1.0$. The structural response of both materials is described by a linear thermo-elastic model and an isotropic, constitutive behavior, with Young's moduli $E_{A}=E_{B} = 1.0$, Poisson's ratios $\nu_{A}= \nu_{B}= 0.3$ and coefficients of thermal expansion $\alpha_A = 1.0 \cdot 10^{-5}$ and $\alpha_B = 10.0 \cdot 10^{-5}$. The reference temperature is set to $T_0 = 0.0$. Perfect bonding is assumed at the interface between the two materials $\Gamma^{AB}$. 

\begin{figure*}[ht]\centering
	\begin{tabular}{cc}
		 \textbf{Local Refinement} & \textbf{Uniform Refinement}\\	
		 \includegraphics[width=0.96\columnwidth]{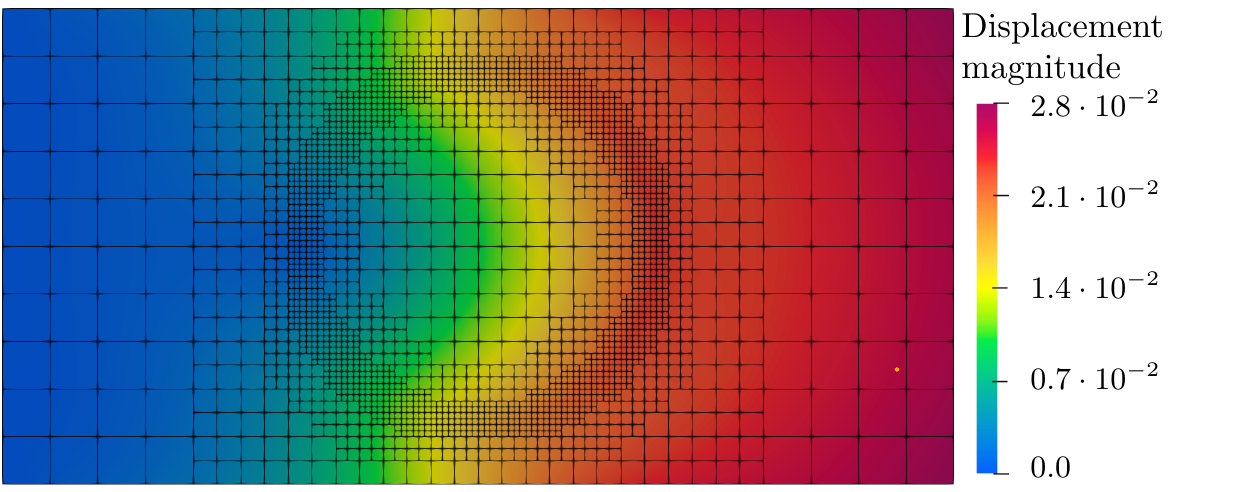}
		&\includegraphics[width=0.96\columnwidth]{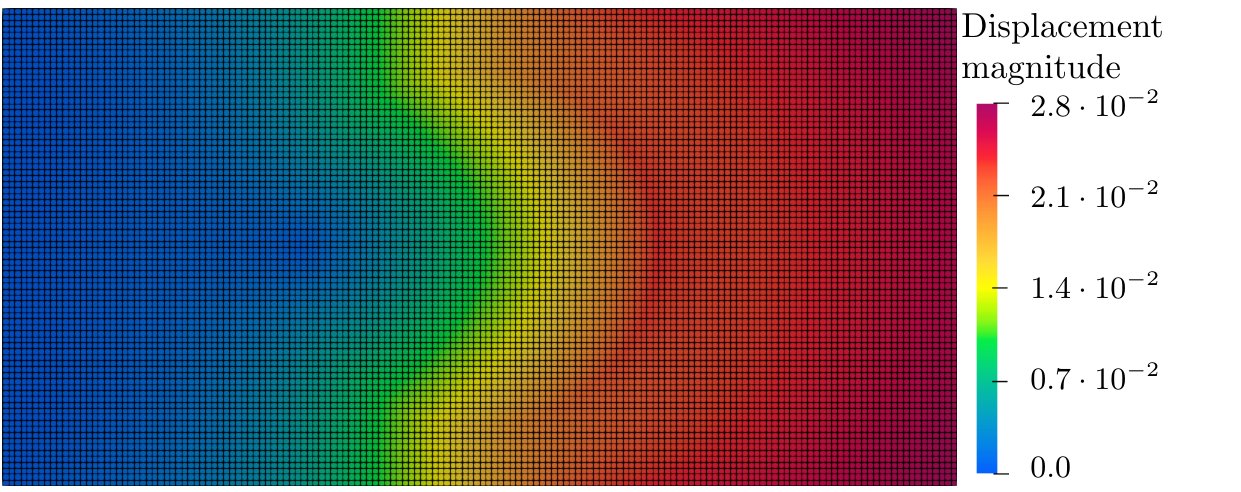}\\
		
		\includegraphics[width=0.96\columnwidth]{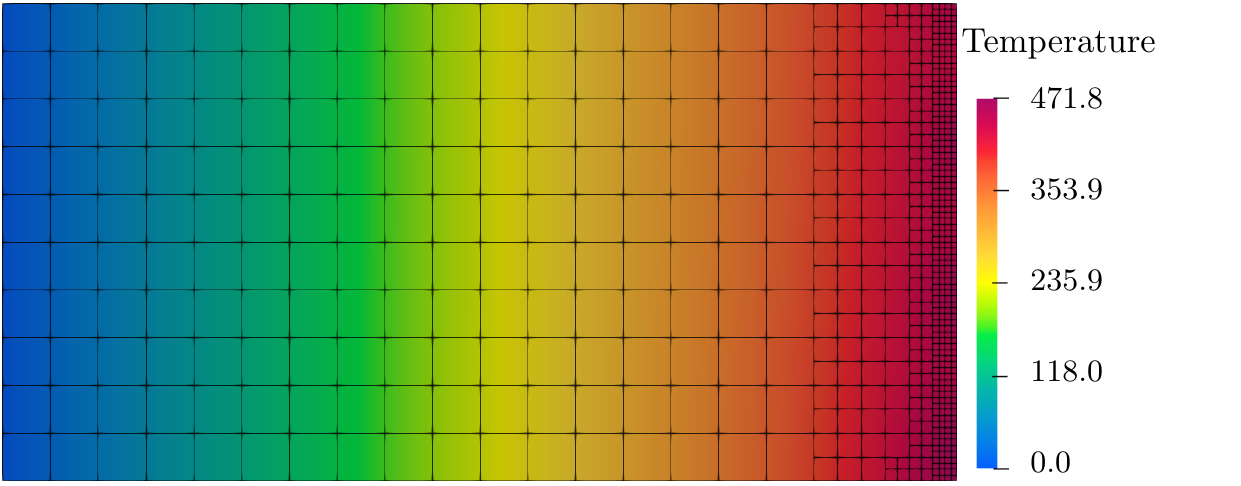}
		&\includegraphics[width=0.96\columnwidth]{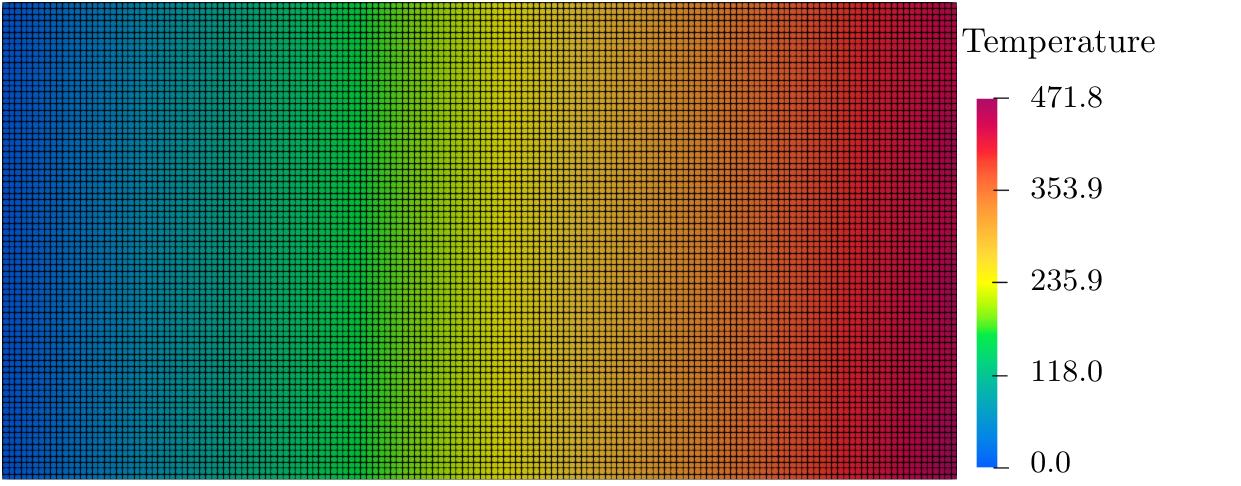}\\
		
		\includegraphics[width=0.96\columnwidth]{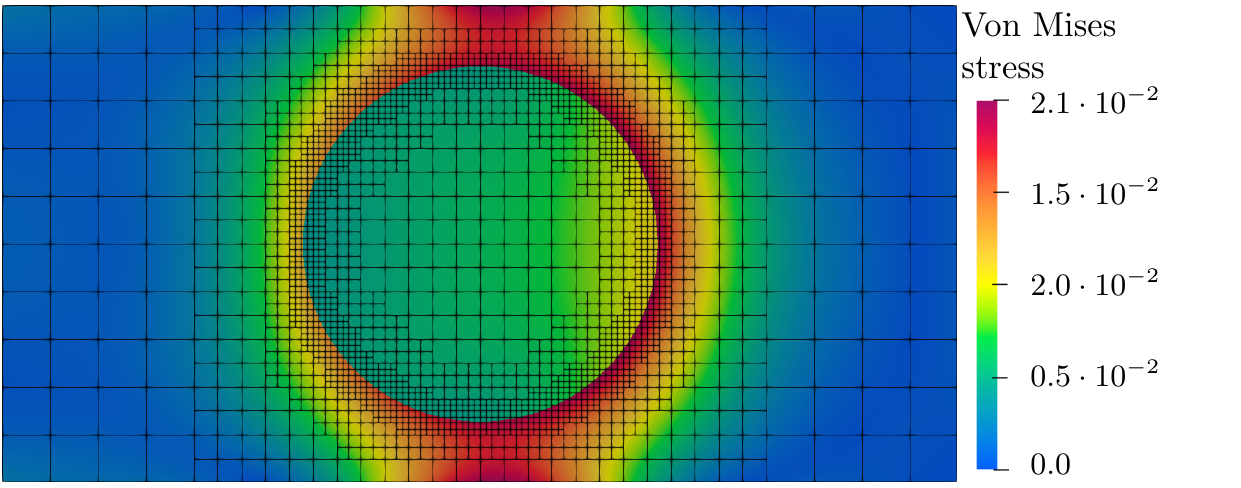}
		&\includegraphics[width=0.96\columnwidth]{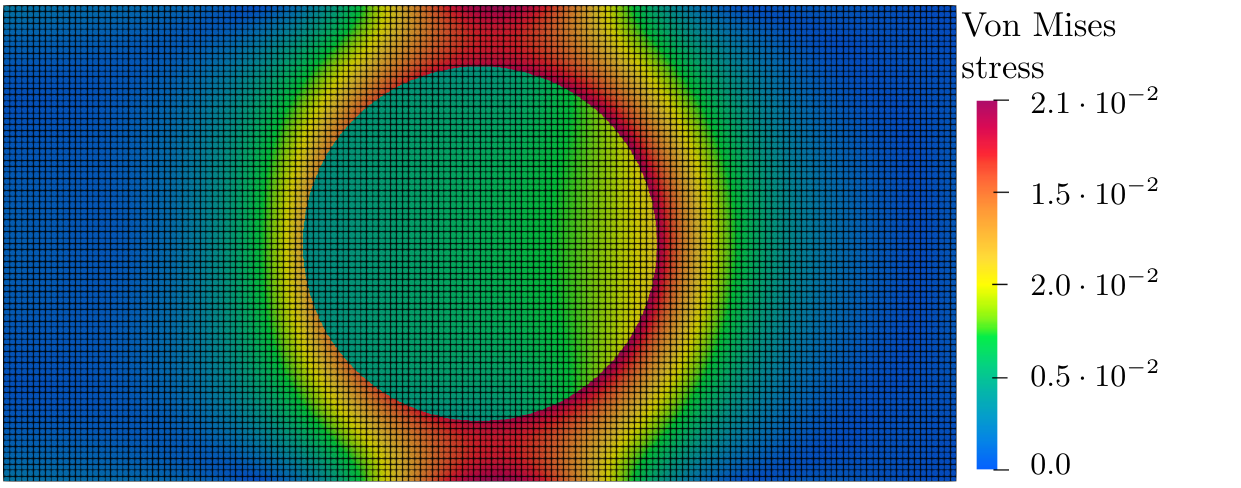}\\
		
		\includegraphics[width=0.96\columnwidth]{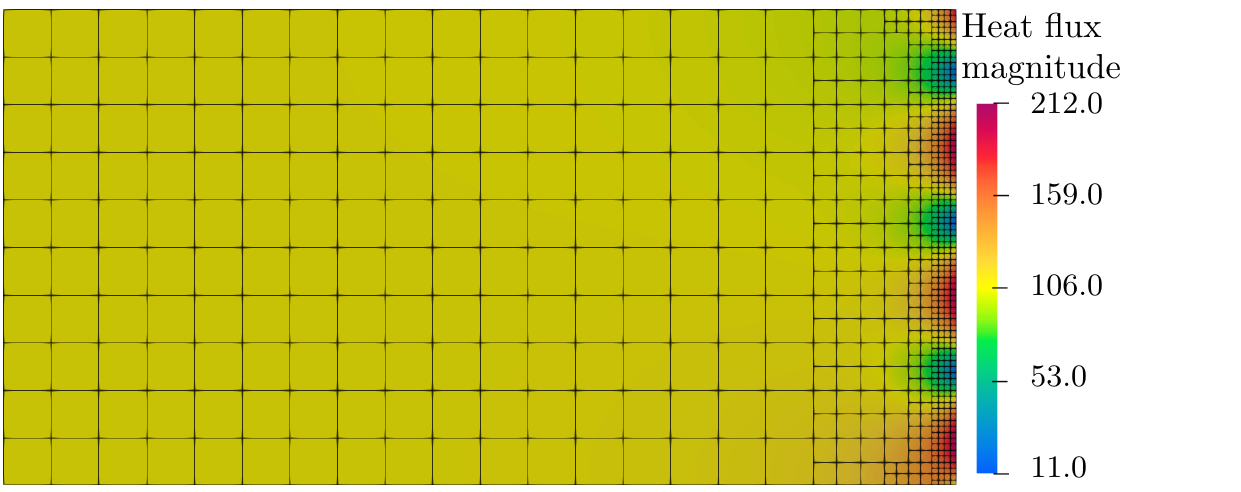}
		&\includegraphics[width=0.96\columnwidth]{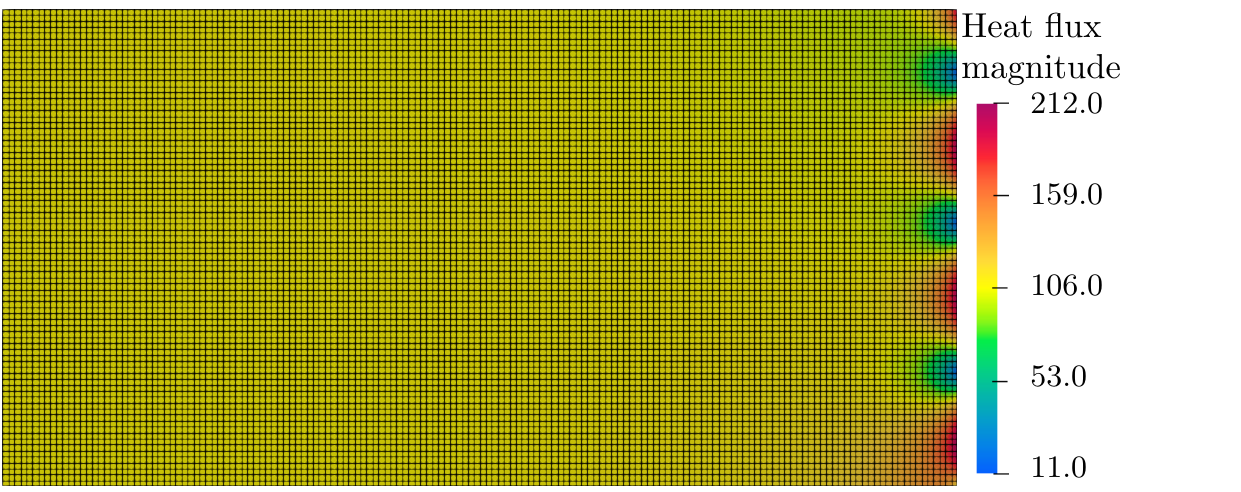}\\
		
	\end{tabular}
	\caption{Comparison of XFEM analysis results for locally and uniformly refined discretizations}
	\label{fig_result_multimaterial}
\end{figure*}

An XFEM analysis is performed to approximate the stress field at the material interface with high accuracy. We compare the temperature, displacement, Von Mises stress, and heat flux using uniformly and locally refined B-spline discretizations. A quadratic polynomial order is used for the temperature and displacement fields. The coarsest THB background mesh has $20 \times 10$ elements.  For the uniform refinement case, the discretizations of the temperature and displacement fields are three times uniformly refined. In the local refinement case, the discretization for the temperature field is three times locally  refined at the left domain boundary to accurately represent the spatially varying load. In contrast, the discretization of the displacement field is three times locally  refined around the circular material interface. The obtained heat flux and stress distributions are presented in in Fig.~\ref{fig_result_multimaterial} and show that the temperature and displacement fields exhibit large spatial gradients in these regions. 

The number of DOFs resulting from the uniformly and locally refined THB background meshes are presented in Table \ref{table_multimaterial}. The locally refined discretization of the thermal field reduces the number of DOFs by a factor of $\sim17$ when compared to the uniform discretization. When locally refining the displacement field, the number of DOFs is $\sim4.6$ times smaller compared to the linear system associated with a uniform discretization. To accurately represent the geometry, the union background mesh is four times uniformly refined, irrespective of the B-spline discretizations of the temperature and displacement fields.

\begin{table}[t!]\centering
	\caption{Comparison of the number of DOFs for temperature and displacement field for and locally uniformly refined discretizations }
	\label{table_multimaterial}
	\renewcommand{\arraystretch}{1.5}
	\begin{tabular}{p{0.20\columnwidth}p{0.30\columnwidth}p{0.34\columnwidth}}
		\hline
		 & Local Refinement & Uniform Refinement\\ \toprule
		$\# \mbox{DOFs}_{\mathrm{Temp}}$ & 818 & 13992\\
		$\# \mbox{DOFs}_{\mathrm{Disp}}$ & 6020 & 27981\\
		\hline
	\end{tabular}
\end{table}

The displacement and temperature fields, as well as the Von Mises stress and the heat flux magnitude contours, are shown in Fig.~\ref{fig_result_multimaterial}. Qualitatively, the resulting fields for uniform and local refinement are equivalent. 

The following four locally refined configurations are examined: (a) refinement for both the structural and thermal field at $\Gamma^{AB}$ and at the Neumann boundary $\Gamma^{N}$; (b) refinement of the structural field at $\Gamma^{AB}$ and refinement of the thermal field at $\Gamma^{N}$; (c) refinement of the structural field at $\Gamma^{AB}$; and (d) refinement of the thermal field at $\Gamma^{N}$. In all cases the union background mesh is chosen to be uniformly refined.

Table \ref{table_error} presents the $L^2$ error norm of the Von-Mises stress at the material interface $\Gamma^{AB}$ for different local refinement configurations together with the size of the linear system. Local refinement for the first two Configurations, (a) and (b), yields the smallest error in the $L^2$ error norm. However, Configuration (b) results in a significant reduction of the computational cost as measured by the number of DOFs.

In addition, Configuration (c) yields a similar small error as Configuration (a) and (b). This behavior might be unexpected but can be explained with the fine uniform refinement of the union background mesh which allows for an accurate integration of the spatial varying load even with a coarse thermal discretization. Configuration (d) results in the highest error as only the thermal field is locally  refined at the Neumann boundary.

\begin{table}[h]\centering
	\caption{Comparison of error of the stress field at the interface $\Gamma^{AB}$ in the $L^2$ error norm for different refinement configurations}
	\label{table_error}
	\renewcommand{\arraystretch}{1.5}
	\begin{tabular}{p{0.52\columnwidth}p{0.17\columnwidth}p{0.14\columnwidth}}
		\hline
		Local Refinement & $L^2$ error norm & \# DOFs\\ \toprule
		$\text{(a) Both fields at }\Gamma^{AB} \text{ and } \Gamma^{N}$  & 1.2023e-6 & 10626\\
		$\text{(b) Disp at }\Gamma^{AB} \text{ and Temp at } \Gamma^{N}$  & 2.8881e-6 & 6838\\
		$\text{(c) Disp at }\Gamma^{AB}$  & 2.8913e-6 & 5764\\
		$\text{(d) Temp at } \Gamma^{N}$  & 3.8764e-4 & 2126\\
		\hline
	\end{tabular}
\end{table}

\section{Conclusion}\label{sec:concl}

This paper presents an immersed isogeometric finite element analysis framework with local mesh refinement based on a Heaviside enriched XFEM. Hierarchical, locally refined discretizations allow for refinement of the finite element approximations in regions of interest, balancing accuracy and computational cost. THB functions are utilized as they provide an elegant way to construct suitable, locally refined discretization spaces. Moreover, B-spline basis functions are an appealing choice over Lagrange basis functions because of their higher inter-element continuity and their increased computational efficiency. In multi-material, multi-physics problems, the resolution requirements may be different for individual state variable fields. The proposed framework allows for separate discretizations with different polynomial orders for each physical field. Furthermore, each discretization can be refined individually, both globally and locally, to meet field-specific accuracy requirements. In contrast to using the same polynomial order and refinement for all state variable fields, the proposed framework may lower significantly the computational cost. In this paper, THB background meshes are refined based on geometric refinement indicators. However, the framework permits any refinement indicators and can be used for adaptive mesh refinement strategies. 

A PT data structure and mesh generation algorithms are presented for the efficient construction of differently refined meshes, both in terms of run time performance and memory needs. The concept of PT cell activation states enables using the same data structure to construct a set of different hierarchically refined THB background meshes. Using the PT data structure and the activation state concept, a union background mesh is constructed such that elements of the union background mesh are at the highest (or higher) refinement level of all corresponding elements of the THB background meshes. The union mesh supports the discretizations of all THB background meshes. 

The union background mesh serves as the XFEM background mesh in which the geometry is immersed. In this paper, the geometry is represented by level set functions, and intersected elements are cut recursively by the zero isocontours of the level set fields. This process yields a single integration mesh which is aligned with the boundaries and interfaces defined by the level set functions. In this paper, the weak form of the governing equations is integrated by standard quadrature rules on the integration mesh. However, other quadrature schemes can be applied, such as the ones proposed by \cite{Thiagarajan2016} and \cite{Gunderman2021}. The B-spline basis functions are enriched using the generalized Heaviside enrichment strategy of \cite{Noel2022}. To facilitate operations performed on the union background mesh and on the integration mesh in the XFEM analysis, the THB basis functions are represented by Lagrange basis functions defined on the union background mesh via extraction operators. The union background mesh supports extraction operators for each THB discretization.

Numerical examples suggest that the proposed immersed isogeometric finite element framework generates discretizations that converge at theoretical, optimal convergence rates with mesh refinement. The application of our framework to a multi-material polycrystalline micro-structure shows that it is well suited to discretize complex multi-material problems in 3D. A scalability study demonstrates that the proposed implementation scales with an increasing number of processors. The coupled thermo-elastic examples highlight the benefits of tailoring the discretization of individual state variable fields to their field-specific accuracy requirements. 

Future work will focus on extending the proposed framework to utilize additional refinement criteria such as finite element error estimators. Moreover, more complex physics for which local mesh refinement is crucial, such as fluid flow at high Reynolds numbers, will be addressed.

\section*{Acknowledgment}\label{Acknowledgment}

The first, second, fourth, and fifth authors acknowledge the support for this work from the Defense Advanced Research Projects Agency (DARPA) under the TRADES program (agreement HR0011-17-2-0022). The first author acknowledges partial auspice of the U.S. Department of Energy by Lawrence Livermore National Laboratory under Contract DE-AC52-07NA27344 (LLNL-JRNL-842737). The third and fifth authors acknowledge the support of Sandia National Laboratories under PO 2120843. The fourth and fifth authors acknowledge the support of the National Science foundation under Grant 2104106. The opinions and conclusions presented in this paper are those of the authors and do not necessarily reflect the views of the sponsoring organization.




\bibliographystyle{abbrvnat}
\bibliography{Literature.bib}

\end{document}